\documentclass[12pt]{article}
\usepackage{fullpage,amsmath,amssymb,amsthm}
\usepackage{subcaption}
\usepackage{graphicx}
\usepackage{float}
\usepackage[utf8]{inputenc}
\usepackage[section]{placeins}
\usepackage{pgffor}
\usepackage{url}
\graphicspath{{./Bickley/}{./Bickley/KMeans_Partitions/}{./Bickley/NegativeMass/}{./Bickley/SPCA_Partitions/}{./Bickley/SpEigFunc/}{./Bickley/Spectrum/}
{./Turbulence/}{./Turbulence/Diagnostics/}{./Turbulence/Eigenvectors/}{./Turbulence/Examples/}{./Turbulence/KMeans_Partitions/}{./Turbulence/NegativeMass/}{./Turbulence/SPCA_Partitions/}{./Turbulence/SparseEigFunc/}{./Turbulence/SparseVectorSum/}{./Turbulence/Spectrum/}{./Turbulence/Minimum_Value/}{./Turbulence/SecondEigFunc/}{./Turbulence/Cheeger/}
{./NorthAtlantic/}{./NorthAtlantic/Diagnostics/}{./NorthAtlantic/Spectrum/}{./NorthAtlantic/Eigenvectors/}{./NorthAtlantic/KMeans_Partitions/}{./NorthAtlantic/NegativeMass/}{./NorthAtlantic/Minimum_Value/}{./NorthAtlantic/SPCA_Partitions/} {./NorthAtlantic/SparseEigFunc/}{./NorthAtlantic/Cheeger/}{./NorthAtlantic/SparseVectorSum/}{./NorthAtlantic/ThirdEigFunc/}
{./OneDim/}
{./Introfigs/}}
\title{Sparse eigenbasis approximation: multiple feature extraction across spatiotemporal scales with application to coherent set identification}
\author{Gary Froyland, Christopher P.\ Rock, and Konstantinos Sakellariou \\ School of Mathematics and Statistics \\
University of New South Wales \\
Sydney NSW 2052, Australia}

\newcommand{\R}{\mathbb{R}}

\DeclareMathOperator{\Min}{Min}
\DeclareMathOperator{\Polar}{Polar}
\DeclareMathOperator{\sech}{sech}
\DeclareMathOperator{\sign}{sign}

\DeclareMathOperator{\Span}{span}
\DeclareMathOperator{\Trace}{Trace}

\DeclareMathOperator*{\argmin}{arg\,min}
\DeclareMathOperator*{\argmax}{arg\,max}

\newtheorem{theorem}{Theorem}[section]
\newtheorem{algorithm}[theorem]{Algorithm}

\newtheorem{assumption}{Assumption}
\newcounter{heuristicstep}

\usepackage[framed]{matlab-prettifier}
\newtheoremstyle{straightremark}
  {\topsep}   
  {\topsep}   
  {}          
  {0pt}       
  {\bfseries} 
  {.}         
  {5pt plus 1pt minus 1pt} 
  {}          

\theoremstyle{straightremark}
\newtheorem{remark}{Remark}[section]

\begin{document}
\maketitle
\begin{abstract}
The output of spectral clustering is a collection of eigenvalues and eigenvectors that encode important connectivity information about a graph or a manifold.
This connectivity information is often not cleanly represented in the eigenvectors and must be disentangled by some secondary procedure.
We propose the use of an approximate sparse basis for the space spanned by the leading eigenvectors as a natural, robust, and efficient means of performing this separation.
The use of sparsity yields a natural cutoff in this disentanglement procedure and is particularly useful in practical situations when there is no clear eigengap.
In order to select a suitable collection of vectors we develop a new Weyl-inspired eigengap heuristic and heuristics based on the sparse basis vectors.
We develop an automated eigenvector separation procedure and illustrate its efficacy on examples from time-dependent dynamics on manifolds.
In this context, transfer operator approaches are extensively used to find dynamically disconnected regions of phase space, known as almost-invariant sets or coherent sets.
The dominant eigenvectors of transfer operators or related operators, such as the dynamic Laplacian, encode dynamic connectivity information.
Our sparse eigenbasis approximation (SEBA) methodology streamlines the final stage of transfer operator methods, namely the extraction of almost-invariant or coherent sets from the eigenvectors.
It is particularly useful when used on domains with large numbers of coherent sets, and when the coherent sets do not exhaust the phase space, such as in large geophysical datasets.
\end{abstract}

\newpage
\tableofcontents

\newpage
\section{Introduction}
Spectral clustering has found broad application in areas such as network analysis, manifold learning (e.g.\ diffusion maps), Lagrangian dynamics, and stochastic processes.
The output of a spectral method is the spectrum and eigenvectors of some matrix or linear operator, typically either of Laplace type such as a graph Laplacian or Laplace-Beltrami operator on a manifold, or of Markov type\footnote{often there is a ``spectral mapping'' relationship between these two types, of the form $M=\exp(L)$, where $M$ is Markov type and $L$ is Laplace type.}.
The eigenvectors carry geometric structure to e.g.\ help discern the topology of an unknown manifold, cluster point data, or analyse nonlinear dynamics.

\begin{figure}[hbt]
  \centering
  \includegraphics[width=\textwidth]{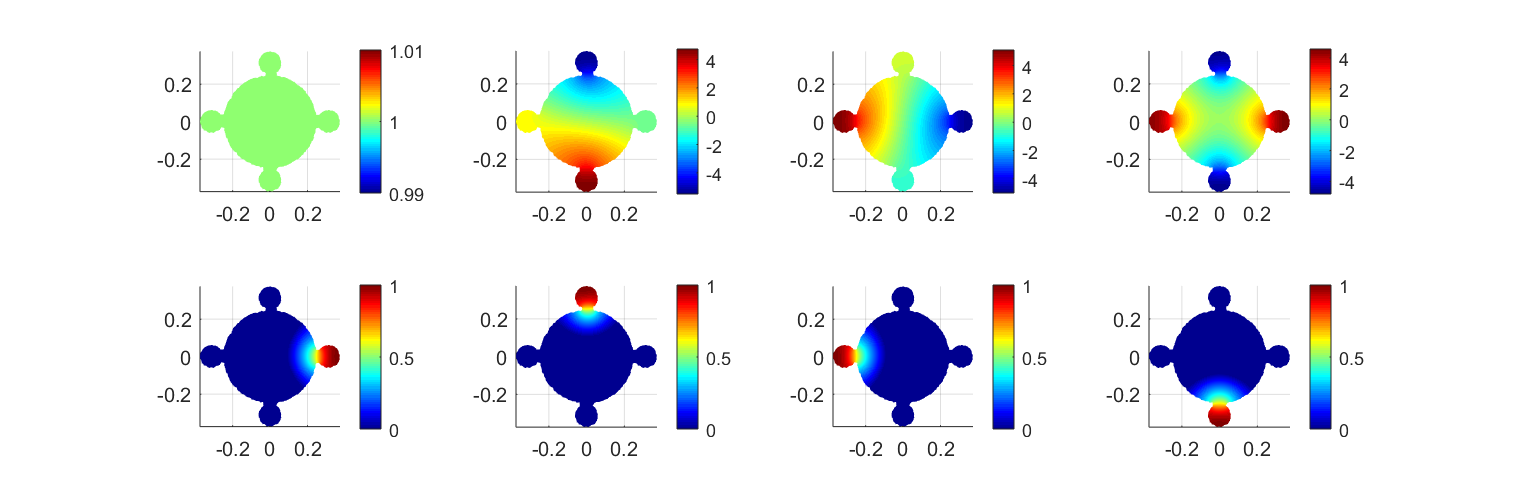}
  \caption{Upper:  the leading four eigenfunctions of the Laplace-Beltrami operator on the manifold shown. Lower: a sparse basis output by Algorithm \ref{alg1}, separating the four ``blobs'' at the periphery of the central disk.  Dark blue represents zero value.  The eigenfunctions were computed using a finite-element mesh with 10774 nodes.}\label{fig:manif}
\end{figure}

For example, the Laplace-Beltrami operator arises as the generator of heat flow on a manifold.
The four leading eigenfunctions of this operator\footnote{with homogeneous Neumann boundary conditions.} (corresponding to the largest four (real) eigenvalues $0=\lambda_1\ge \cdots\ge \lambda_4$) are shown in the upper row of Figure \ref{fig:manif}.
These eigenfunctions -- beyond the first ``trivial'' constant eigenfunction corresponding to $\lambda_1=0$ and shown in Figure \ref{fig:manif} (upper left) --  describe signed ``heat modes'' that decay most slowly under heat flow, at rates $\sim e^{\lambda_it}$, $i=2,\ldots,4$.
Because the four ``blobs'' at the periphery of the main disk have narrow channels connecting them to the main disk, the slowest decaying mode (second from the left in the upper row of Figure \ref{fig:manif}) has a lot of heat (red) in the bottom blob and lower part of the main disk and is cool (blue) in the upper blob and upper part of the main disk.
The four ``blobs'' are correctly detected in the upper row of Figure \ref{fig:manif}, but they are mixed together in the leading four eigenfunctions.
We wish to separate these four features and this is what is achieved 
in the lower row of Figure \ref{fig:manif}, by finding an approximate sparse basis through orthogonal rotation of the four leading eigenfunctions.

\begin{figure}
  \centering
  \includegraphics[width=10cm]{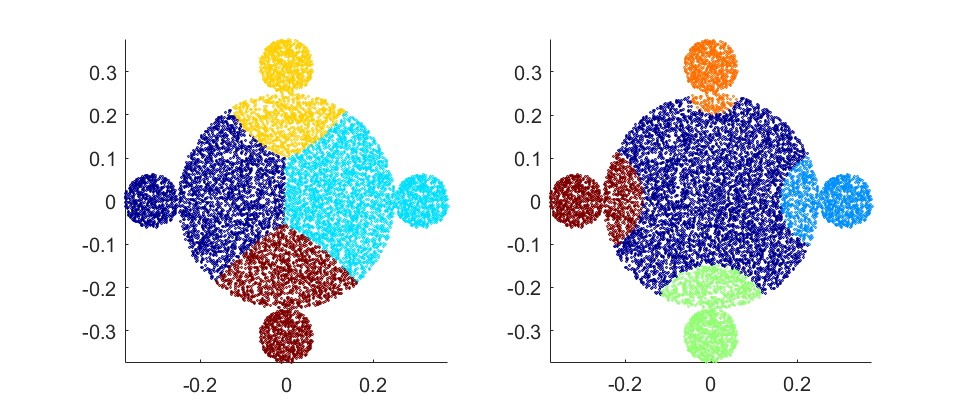}
  \caption{Left: clustering an embedding of the second to fourth eigenvectors from Figure \ref{fig:manif} in $\mathbb{R}^3$ (10774 points in $\mathbb{R}^3$) and applying k-means clustering, requesting four clusters. Right: As per left, but requesting five clusters.}\label{fig:kmeans}
\end{figure}

Arguably the most popular method of performing this separation is the embedding of (discrete versions of) the second to fourth eigenfunctions in $\mathbb{R}^3$ and applying a hard clustering algorithm such as k-means;  see e.g.\ \cite{jain10} for a detailed description and the historical development of k-means.
In this example, for the k-means approach to at least partially isolate the peripheral blobs in the example of Figure \ref{fig:manif}, one has to request a clustering of \emph{five} clusters because k-means \emph{partitions} the dataset/domain, and requesting four clusters groups each of the peripheral blobs with approximately one quarter of the large central disk; see Figure \ref{fig:kmeans}.
A key advantage of our approach is that we do \emph{not necessarily} separate by partitioning, but instead may classify a large part of the dataset/domain as ``unclustered''.
Such an idea has been partly implemented in \cite{hadjighasem16} by adding one more cluster to k-means as a ``background'' cluster, but we show in Sections \ref{sec:turbulence} and \ref{sec:northatlantic} that this is not a good approach for more challenging dynamics.
Other separation approaches are discussed in Section \ref{sec:hardcluster}.

Perhaps the most similar precursor to our approach in the dynamical systems literature, applied in the context of Markov chains, is the search for a linear transformation of the leading eigenvectors \cite{deuflhard} to make them as close as possible to a collection of indicator functions (vectors taking a small number of distinct values) spanning a similar subspace, under constraints of non-negativity and small intersection of support.
The optimisation in \cite{deuflhard} can be expensive for large vectors and large numbers of vectors, and will always return a full partition of the domain.
An LU decomposition approach \cite{berrysauer} has been effective in the case where the eigenvectors are already close to linear combinations of indicator vectors.
This is a somewhat idealised situation and, in such a setting, all of the above methods should perform very well.

In the context of general Markov processes and Lagrangian dynamics, there are many situations where it is not instructive to \emph{partition} all nodes or the entire manifold into almost-disconnected pieces according to the relevant eigenvectors.
Rather, it is more informative to know that some parts of the network or manifold have a natural almost-disconnected structure, while the remainder (possibly a large remainder) of the network or manifold is well connected.
The sparsity built in to our new sparse basis approach will automatically identify those nodes or subsets of a manifold that do not belong to any almost-disconnected subset;  such nodes or subdomains will be zeroed out by the sparsification.


In Section \ref{sec:examples} we illustrate our approach on spectral algorithms designed to identify almost-invariant or coherent sets in dynamical systems on manifolds.
The term \emph{almost-invariant set} refers to a set that is approximately invariant under autonomous dynamics.
The term \emph{coherent set} refers to a time-dependent family of sets that are approximately equivariant under time-dependent or nonautonomous dynamics.
In contrast to an almost-invariant set, which is fixed in phase space, the family of coherent sets are in general mobile in phase space as time evolves;  see \cite{FPG14} for a short overview.
Almost-invariant and coherent sets in the phase space of complicated dynamics are important because they are the most predictable features in the medium term.
These concepts have been used in molecular dynamics to identify conformations to aid drug design \cite{deuflhard_schuette}, in atmospheric dynamics to identify atmospheric vortices \cite{FSM10}, in oceanography to identify ocean gyres \cite{FPET07} and eddies \cite{FHRSS,FHRvS}, 
and in fluid dynamics \cite{stremler_etal}.
Spectral approaches to the identification of coherent features rely on the dominant eigenvectors of transfer operators or related operators, such as the dynamic Laplacian \cite{F15,FK17}.
There are various categories of constructions, based on the problem to be solved;  see   Table \ref{tab:TO}.
The final stage of the extraction of almost-invariant or coherent sets from eigenvectors or singular vectors has previously proceeded either by time-consuming hand-selected thresholds or by common clustering methods such as k-means.


A main goal of this research is to describe  simple and robust procedures for automatically separating features represented in the eigenvectors or singular vectors arising from transfer operator methods.
Our separation procedures streamlines the use of these methods  on domains with large numbers of coherent sets, particularly when the coherent sets do not exhaust the phase space.
Although the dynamical context is our primary motivation,  we believe our approach will translate well to other application domains  where disentangling the output of spectral clustering is important.
We have thoroughly tested our approach on a variety of nonlinear dynamics and have  proposed heuristics that have consistently performed well in our tests.
We illustrate our new approach by separating coherent features in idealised fluid flow dynamics and in large-scale ocean currents in the North Atlantic derived from satellite altimetry.

Our main contributions include a refinement of the standard eigengap criterion for determining suitable numbers of features.
Our Weyl-inspired refinement strongly highlights and appropriately scales gaps in the spectrum, emphasising natural time scales.
We also develop new vector-based heuristics for determining suitable numbers of features;  these help to determine good choices of dimension for the subspaces and therefore the corresponding spatial scales of features.
We develop a ``reliability'' heuristic for the sparse basis vectors, including feature rejection criteria.
These vector-based heuristics supplement the eigengap criteria and we have found them to be particularly useful in realistic situations where there is no clear eigengap.
Finally, we propose automatic methods of hard thresholding of the sparse basis to provide hard separation of features.

An outline of the paper is as follows.
Section \ref{sec:background} provides a brief background on hard separation of spectral clustering output, sparse principal components analysis, and transfer operator and dynamic Laplace operator constructions in nonlinear dynamics.
Section \ref{sec:sparsebasis} describes in detail the variant of sparse basis approximation we will use, based on the SPCArt algorithm \cite{hu16}.
Section \ref{sec:techniques} discusses the use of our sparse basis algorithm, including new eigenvalue and vector based heuristics for selecting the size of the basis and assigning reliability to the sparse basis output, and simple thresholding procedures.
We present our numerical examples in Section \ref{sec:examples} and conclude in Section \ref{sec:conclusion}.

\section{Background} \label{sec:background}
\subsection{Separation of spectral clustering output by hard clustering}
\label{sec:hardcluster}
Denote the eigenvalue (resp.\ eigenvector) output of a spectral clustering method by ordered lists of eigenvalues $\lambda_1\ge\lambda_2\ge\cdots$ and their corresponding eigenvectors $v_1, v_2,\ldots$.
Suppose that we wish to separate the features encoded in the vectors $v_1,\ldots,v_r$.
One of the most popular separation approaches is to perform  hard clustering on the $p\times r$ data array $V=[v_1|\cdots|v_r]$, treating each row of this array as a Euclidean point in $\mathbb{R}^r$.
Figure \ref{fig:kmeans} illustrates this procedure with $r=3$ using k-means for the hard clustering.
In the context of nonlinear dynamics, the k-means algorithm \cite{jain10} has been applied to eigenvector data embedded in Euclidean space, for example \cite{hadjighasem16,banisch17,schlueter-kuck17,padberg-gehle17,FJ18} to identify coherent sets.
A ``soft'' variant, fuzzy c-means, has earlier been used in the dynamics context on eigenvector data \cite{FD03,F05} to identify almost-invariant sets.
Other variants include k-medians \cite{jain88}, k-medoids \cite{kaufman87}, and kernel k-means \cite{scholkopf98}.

Other hard separation approaches for spectral clustering include \cite{yu_shi03}, who try to find a rotation of the eigenvectors that brings them as close as possible to a collection of binary vectors, which represent a feature partition.
The work of \cite{leeleelee15} builds on \cite{yu_shi03} by proposing rotating the eigenvectors to approximately achieve a collection of nonnonegative vectors, followed by a maximum likelihood assignment.
Related to this is \cite{han_xiong17} who search for nonnegative vectors ``nearby'' the  eigenvectors via a penalty term.
Another recent contribution includes \cite{lu_yan16}, who search for sparse orthogonal vectors that maximise the trace of the usual bilinear form involving the Laplace matrix, followed by an application of k-means.
All of the above approaches produce a \emph{partition} of the dataset/domain (as does k-means), and this is something we wish to avoid through the use of sparse bases.
Moreover, the explicit imposition of nonnegative vectors removes what we have found to be an important indicator of the quality of the separation, as we show in Section \ref{sec:techniques}.
There are several other related methods that involve manipulating the weights in the affinity matrix, but we omit discussing these because in our dynamical systems applications the Markov and Laplace operators that arise have additional special properties we wish to preserve.

\subsection{Sparse principal component analysis}
Principal component analysis (PCA) is an unsupervised learning technique often used for dimension reduction and feature selection.
PCA may be viewed as the process of finding a linear projection from a high-dimensional space (occupied by the dataset) onto a low-dimensional space that preserves as much of the variance of the original dataset as possible.
One drawback of PCA is that every variable in the lower-dimensional space generally depends on \emph{all} of the variables in the high-dimensional space.

This shortcoming motivated sparse principal component analysis (SPCA), which modifies PCA by imposing constraints (resp.\ adding a penalty) to force (resp.\ encourage) each variable in the lower-dimensional space to depend on only a subset of the original variables.
In 2003, Joliffe et al. proposed the first sparse principal component analysis technique, called SCoTLASS \cite{jolliffe03}, and several versions of SPCA have appeared since then.
Many such versions, termed \emph{deflation methods}, find one sparse principal component at a time, then deflate the original data in that direction and find the next orthogonal principal component, e.g.  \cite{jolliffe03,d'aspremont08}.
To avoid sequential deflation finding a non-global optimum, so-called \emph{block methods} have been proposed, which optimise all the subspaces at once \cite{journee10a,d'aspremont07,zou_hastie,hu16}.

\begin{figure}[hbt]
\centering
\includegraphics[width=0.75\textwidth]{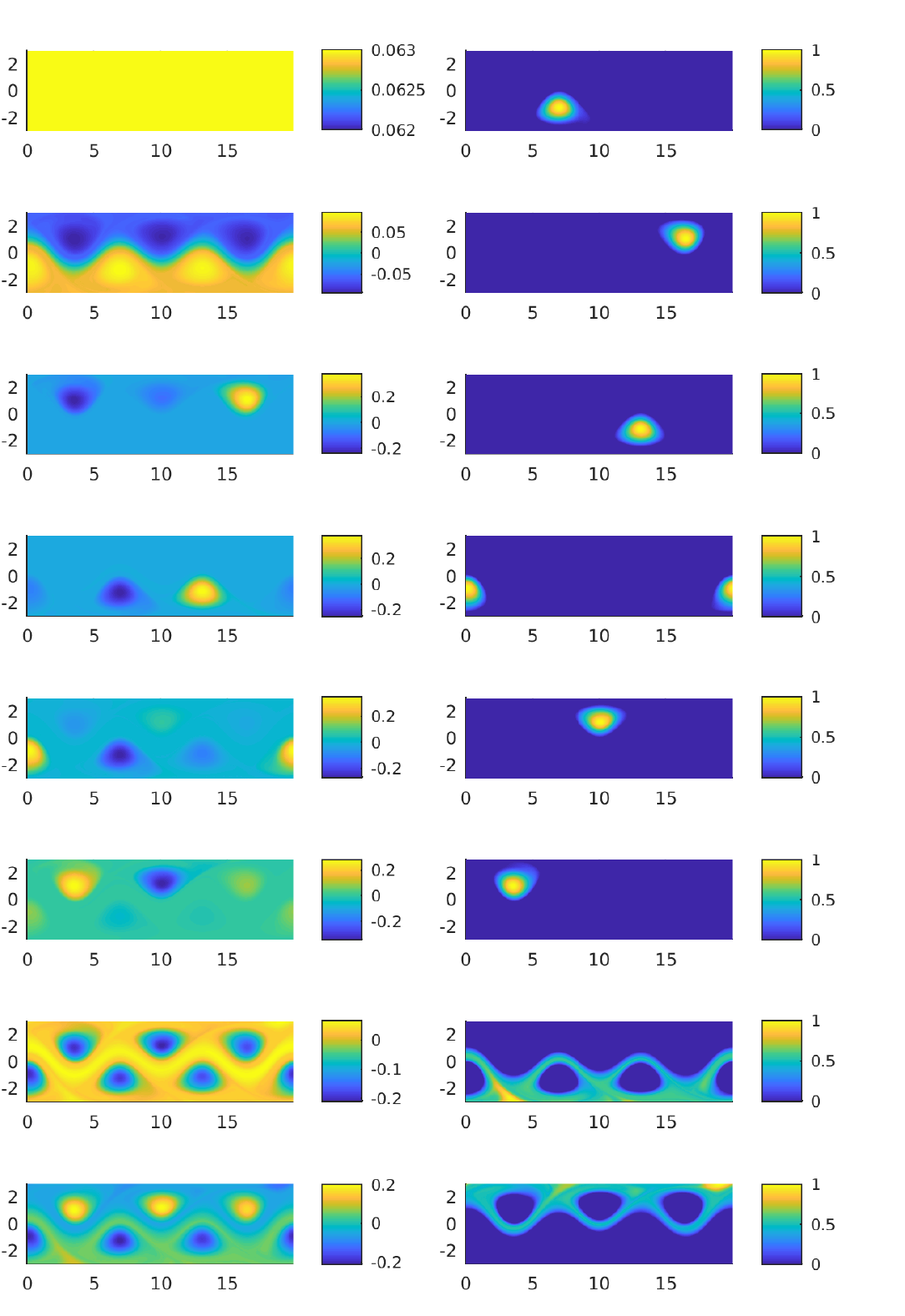}
\caption{Left: the 8 leading eigenvectors $v_1,\ldots,v_8$ (upper to lower) of the dynamic Laplacian constructed from the Bickley jet (see Section \ref{sec:bickley}).  Right: an approximate sparse basis $s_1,\ldots,s_8$ output by Algorithm \ref{alg1} separating the 8 main coherent sets;  small negative values of have been removed from the vectors $s_7, s_8$.} \label{fig:bickley_example}
\end{figure}

Because of its simplicity, efficiency, and robustness, in this paper we will use a tailored version of the sparse principal component analysis by rotation and truncation (SPCArt) algorithm \cite{hu16}, which is a block method.
Letting $\mathcal{V}=\Span\{v_1,\ldots,v_r\}$ we implicitly assume that the eigenvectors from our spectral clustering satisfy the following assumption.
\begin{assumption} \label{assumptionA}
There is a $\mathcal{V}$-approximating, $r$-dimensional subspace $\mathcal{S}\subset \mathbb{R}^p$ with a nonnegative sparse basis $s_1,\ldots,s_r$ whose supports have little or no intersection.
\end{assumption}
Assumption \ref{assumptionA} should hold when there exist pairwise disjoint subsets of the graph or manifold that are sufficiently strongly disconnected, provided suitable eigenvectors $\{v_1,\ldots,v_r\}$ of a Laplace-type or Markov-type operator or matrix are chosen.
Figure \ref{fig:manif} (lower row) and Figure \ref{fig:bickley_example} (right column) illustrate the type of sparse basis vectors we will obtain.


\subsection{Markov- and Laplace-type operators arising in nonlinear dynamics}

Our examples in Section \ref{sec:examples} will be drawn from nonlinear dynamics and we provide a brief overview here, referring the reader to original papers for further details.
One has a smooth domain $M$, typically a smooth $d$-dimensional manifold or subset of $\mathbb{R}^d$, on which the nonlinear dynamics acts.

There are several related linear operators that are intimately connected with the nonlinear evolution on phase space.
These operators fall into two broad classes, which we will call \emph{Markov-type} and \emph{Laplace-type}.
The former have their spectrum contained in $\{z\in\mathbb{C}: |z|\le 1\}$, the closed unit disk in the complex plane, and  the spectrum of the latter is contained in $\{z\in\mathbb{C}:\Re(z)\le 0\}$, the closed half plane with non-positive real part.

In the most straightforward case, the dynamics is generated by repeated application of a transformation $T:M\to M$ or by solving a differential equation $\dot{x}=F(x)$.
Both of these mechanisms generate \emph{autonomous} dynamics because the underlying governing dynamics does not depend on time.
In the former case one uses left\footnote{the description of left and/or right eigenvectors has been chosen to match the original papers.  For Markov-type operators, multiplying vectors on the left of matrices is consistent with most Markov chain texts.} eigenvectors of the transfer operator (Markov type) to identify subsets of the domain that are approximately invariant under repeated application of the dynamics;  see Construction 1 in Table \ref{tab:TO}.
\begin{table}[hbt]
  \centering
\makebox[\textwidth][c]{
  \begin{tabular}{|l|c|c|c|c|c|}
  \hline\hline
  Construction & Autonomous/ & Finite/ & Closed/ & Uses time & Objects  \\
     & Nonauton. & Infinite & Open & derivative & Identified  \\
         & Dynamics & Time & Dynamics &  &  \\ \hline\hline
1.\  Left eigenvectors of& Auton.& $\infty$ & closed & no & almost- \\
  transfer operator \cite{DJ99}&&&&&invariant sets \\ \hline
2.\   Left eigenvectors of& Auton.\ or& $\infty$ & closed & yes & almost- \\
    generator \cite{FJK13,FK17}&Periodic&&&&invariant sets \\ \hline
3.\ Right eigenvectors of transfer& Auton. & $\infty$ & open & no  & Basins of\\
  operator (resp.\ generator) \cite{koltai11,FSvS14}&&&&(resp.\ yes)&  attraction \\ \hline
4.\   Right eigenvectors of& Auton. & finite & closed & no & finite-time \\
    symmetrised transf.\ op. \cite{F05,FPG14}&&&&&almost-inv.\ sets  \\ \hline
5.\  Singular vectors of& Nonauton. & finite & closed & no & finite-time  \\
   transfer operator \cite{FSM10,F13}&&&&&coherent sets \\ \hline
6.\   Eigenvectors of dynamic& Nonauton. & finite & closed & no & finite-time  \\
   Laplacian \cite{F15,FK17a}&&&&&coherent sets \\ \hline
7.\   Oseledets vectors of& Nonauton. & $\infty$ & closed & no & coherent sets \\
   transfer operator cocycle \cite{FLS10}&&&&& \\ \hline
  \hline
\end{tabular}
}
  \caption{Summary of transfer operator and dynamic Laplace operator constructions to identify almost-invariant and coherent sets, and the dynamical settings in which each construction is used.}\label{tab:TO}
\end{table}
In the latter case, one uses left eigenvectors of the generator (Laplace type) to detect almost-invariant sets in continuous time without trajectory integration;  this is Construction 2 in Table \ref{tab:TO}.
When there is open dynamics created by the existence of a ``hole(s)'' in phase space, one is interested in the basins of attraction;  this is Construction 3 in Table \ref{tab:TO}, covering both discrete and continuous time.
In many instances, there is a natural timescale of interest and one wishes to find subsets that are approximately invariant over this \emph{finite} timescale;  this is addressed in Construction~4 in Table \ref{tab:TO}.

\emph{Nonautonomous} dynamics arises when the generating law varies with time.
In discrete time one applies a sequence of different transformations $T_{n}\circ T_{n-1}\circ\cdots\circ T_1$, and in continuous time, one solves a differential equation with time-dependent right-hand-side $\dot{x}=F(x,t)$.
Because the underlying dynamics changes over time, one wishes to find subsets that remain ``coherent'' over the finite time duration;  this is addressed in Constructions 5 and 6 in Table \ref{tab:TO}.
Finally, identifying coherent sets in time-dependent dynamics over an infinite time period is addressed in Construction 7 in Table \ref{tab:TO}.
While we have described the dynamics above as deterministic, each of the Constructions 1--7 can also be applied to stochastic dynamics arising from transformations with additional noise, or stochastic differential equations.

Each of the Constructions 1--7 in Table \ref{tab:TO} is a type of spectral clustering;  Constructions 1, 3, 4, 5, and 7 produce Markov-type operators and Constructions 2, 3, and 6 produce Laplace-type operators.
The eigenvectors arising from Constructions 4, 5, and 6 are orthogonal in a suitable inner product space, while the vectors produced by Constructions 1, 2, 3, and 7 need not be.
For the latter constructions, one applies e.g.\ QR factorisation to orthogonalise the vectors before applying Algorithm \ref{alg1}.

There are various ways in which the constructions in Table \ref{tab:TO} can be numerically implemented.
A thin spline implementation of Construction 5 is introduced in \cite{williams15}, and a spectral implementation in continuous time of Construction 5 is presented in \cite{denner15}.
A radial basis function implementation of Construction 6 is contained in \cite{FJ15} and a finite element implementation is introduced in \cite{FJ18};  the latter is the numerical implementation we will use for all dynamic Laplacian eigenvector computations in this paper.
A graph-based method similar to Construction 6 is described in \cite{hadjighasem16}.
Diffusion map implementations of Constructions 5 and 6 are described in \cite{banisch17}.

In section \ref{sec:examples} we will apply our method to output from Constructions 5 and 6;  the former is of Markov type and the latter is of Laplace type.
While the different constructions in Table \ref{tab:TO} solve different dynamical problems, the properties of the eigenvectors from all constructions are strongly related because they arise from spectral approaches.
We expect our methods to perform equally well on all constructions in Table \ref{tab:TO}.

\section{Sparse eigenbasis approximation} \label{sec:sparsebasis}

Let $v_1,\ldots,v_r\in \mathbb{R}^p$ be a set of linearly independent eigenvectors or singular vectors;  typically $r\ll p$.
These vectors form a basis for a subspace $\mathcal{V}\subset\mathbb{R}^p$.
We wish to transform these vectors to a basis of sparse vectors $s_1,\ldots,s_r\in \mathbb{R}^p$ for a subspace $\mathcal{S}\subset\mathbb{R}^p$ with $\mathcal{V}\approx \mathcal{S}$.
Without loss of generality, we can assume that the vectors $v_1,\ldots,v_r$ are an orthonormal basis;  if not, we perform a QR factorisation of the matrix $V=[v_1|\cdots|v_r]$ and extract the orthonormal columns.
If the subspace $\mathcal{V}$ satisfies Assumption A, we can expect the sparse vectors $s_1,\ldots,s_r$ to have supports with small overlaps,  and therefore  be close to orthogonal.
Thus, we wish to find an  orthogonal ``rotation'' matrix $R$ that transforms $v_1,\ldots,v_r$ into $s_1,\ldots,s_r$ such that (i) the $r$-dimensional space $\mathcal{S}\subset \mathbb{R}^p$ spanned by $s_1,\ldots,s_r$ is close to $\mathcal{V}\subset\mathbb{R}^p$ (approximate subspace preservation) and (ii) the vectors $s_1,\ldots,s_r$ are sparse.

The approach we take is based on the SPCArt algorithm (sparse PCA via truncation and rotation \cite{hu16}). We present below a simplified ``full rank'' version of SPCArt, tuned to the setting of Assumption \ref{assumptionA}, which we call \emph{Sparse Eigenbasis Approximation} (SEBA).
Let $\mathfrak{S}^{r}$ denote the Stiefel manifold $\{A \in \R^{r \times r} : A^\top A = I_r\}$ of $r\times r$ orthogonal matrices, and let $\mathfrak{U}^{p,r}=\{A\in \mathbb{R}^{p\times r}: \mbox{each column of $A$ has $\ell_2$ norm 1}\}$.
Define $\|A\|_F:=\sqrt{\sum_{i=1}^p\sum_{j=1}^r A_{ij}^2}$ to be the Frobenius norm and $\|A\|_{1,1}:=\sum_{i=1}^p\sum_{j=1}^r |A_{ij}|$ to be the $\ell_{1,1}$ matrix norm.
Given some small positive sparsity parameter
$\mu$, we wish to solve the following nonconvex optimisation problem for $S\in \mathfrak{U}^{p,r}$ and $R\in \mathfrak{S}^{r}$:
\begin{align}
\argmin_{\substack{S\in\mathfrak{U}^{p,r} \\ R \in \mathfrak{S}^{r}}} \frac{1}{2} \|V-SR\|_F^2 + \mu \|S\|_{1,1}, \label{eq:spcart}
\end{align}
The $r$ columns of $S$ define a sparse basis $\{s_1,\ldots,s_r\}$ for a subspace $\mathcal{S}$ close to $\mathcal{V}$.
The first term in (\ref{eq:spcart}) measures how close\footnote{Note that if we define $S:=VR^\top$ for some $r\times r$ orthogonal (in fact, $R$ nonsingular is sufficient) matrix $R$, then the columns of $S$ will span the same subspace as the columns of $V$.} the rotated columns of $V$, namely $VR^T$, are to the columns of $S$ (recall $\|VR^\top-S\|_F=\|V-SR\|_F$ for orthogonal $R$).
The second term measures the sparsity of the basis formed by the columns of $S$ using the common sparsity-inducing \cite{donoho} $\ell_1$ penalty.



Because of the nonconvexity of the problem (\ref{eq:spcart}), finding a global optimum $(S,R)\in\mathfrak{U}^{p,r}\times\mathfrak{S}^{r}$ is difficult.
Hu \emph{et al.} \cite{hu16} proposed alternately fixing $R$ and optimising $S$, and fixing $S$ and optimising $R$ to find a local minimum $(S,R)$.
Each of these individual optimisation problems are fast to solve exactly.
\begin{enumerate}
\item \textbf{Fixed $R$:} The optimisation for $S$ in \eqref{eq:spcart} is exactly solved by ``soft thresholding'' \cite{journee10a}. Define a thresholding transformation $C_\mu:\mathbb{R}\to\mathbb{R}$ by $C_\mu(z)=\sign(z)\max\{|z|-\mu,0\}$;  $C_\mu$ is applied to vectors elementwise. For $j=1,\ldots,r$, set the $j^{th}$ column of $S$ to $S_j=C_\mu((VR^\top)_j)/\|C_\mu((VR^\top)_j)\|$.
\item \textbf{Fixed $S$:} The optimisation for $R$ becomes a Procrustes problem $\min_{R \in \mathfrak{S}^{r}} \frac{1}{2} \|V-SR\|_F^2$.
This problem may be efficiently exactly solved by polar decomposition \cite{zou_hastie}. Set $R=\Polar(S^\top V)$, where $\Polar(\cdot)$ is the orthonormal component of the polar decomposition;  that is, $S^\top V=RH$, where $R$ is orthogonal and $H$ is symmetric and positive definite. If $S^\top V$ has singular value decomposition $S^\top V=PDQ^\top$, then $R=PQ^\top$.   In MATLAB, we use \lstinline[style=Matlab-editor]"[P,~,Q]=svd(S'*V,0);"
\end{enumerate}
One initialises with $R=I_r$ and alternately applies steps 1 and 2 above until the change in $R$ is below a specified tolerance\footnote{\label{foot:tol}We used a tolerance of $10^{-14}$ in our experiments, and we observed there is little, if any, difference if the tolerance is increased by a few orders of magnitude.}.
The parameter $\mu$ should be chosen less than $1/\sqrt{p}$ as $C_\mu$ sends the constant unit vector to the zero vector for $\mu\ge 1/\sqrt{p}$.
Hu \emph{et al.} \cite{hu16} suggest using $\mu=1/\sqrt{p}$ and with extensive experimentation we have found that this works well, as do values between $75\%$ and $100\%$ of $1/\sqrt{p}$.
In our experiments we have used $\mu=0.99/\sqrt{p}$ because in many cases the constant unit vector is part of our initial basis $\mathcal{V}$.
Alternative thresholding methods were presented in \cite{hu16}, but in our experiments we found soft thresholding was the most robust and we use this in all computations.
In summary, our sparse eigenbasis approximation (SEBA) algorithm is:
\begin{algorithm}[(SEBA) -- Input orthonormal $p\times r$ matrix $V$; output sparse $p\times r$ matrix $S$]
\label{alg1}\
\begin{enumerate}
\item \label{itm:setmu} Set $\mu=0.99/\sqrt{p}$ and $R=I_r$.
\item \label{itm:setS} For $j=1,\ldots,r$, set the $j^{th}$ column of $S$ to $S_j=C_\mu((VR^\top)_j)/\|C_\mu((VR^\top)_j)\|$.
\item \label{itm:setR} Set $R=\Polar(S^\top V)$.
\item \label{itm:iter} If the matrix 2-norm of the difference between the revised $R$ in step 3 and the previous $R$ is larger than some tolerance\textsuperscript{\ref{foot:tol}}, go to step 2;  otherwise go to step 5.
\item \label{itm:withsign} For $j=1,\ldots,r$, set $S_j\to\sign(\sum_{i=1}^p S_{ij})S_j$.
\item \label{itm:setmax} For $j=1,\ldots,r$ set $S_j\to S_j/\max_{1\le i\le p}S_{ij}$.
\item \label{itm:negsort} Reorder the columns of $S$ so that $m_j:=\min_{i}S_{ij}$ is in decreasing order, $j=1,\ldots,r$.
\end{enumerate}
\end{algorithm}
In Step 5 we possibly change the sign of the columns of $S$ to ensure they are predominantly nonnegative and in Step 6 we scale the columns of S so that their maximum value is 1;  we will later interpret $S_{ij}$ as a likelihood of membership of index $i$ in the $j^{th}$ feature.
Step 7 orders the columns of $S$ in terms of ``reliability''; see Section \ref{sec:reliability} for a discussion.

As noted in \cite{hu16}, SPCArt (steps 1--4 of Algorithm \ref{alg1}) produces vectors $s_1,\ldots,s_r$ that should span a similar subspace to $\mathcal{V}$ because of the first term of (\ref{eq:spcart}).
Because the input vectors $v_1,\ldots,v_r$ are orthonormal and the vectors $s_1,\ldots,s_r$ are arrived at by orthogonal rotation and truncation, the latter should also be close to orthogonal.
The sparsity should be reasonably balanced because of the global optimisation of $R$ in step 3 and uniform thresholding across each vector in step 2.
A short MATLAB code listing to implement Algorithm \ref{alg1} is included in Appendix \ref{code:alg1}.
%
%

\section{Using the SEBA Algorithm \ref{alg1}}\label{sec:techniques}
In this section we describe the ``reliability'' ordering in Step 7 of Algorithm \ref{alg1}, and introduce some rules of thumb for selecting an appropriate basis size for input to Algorithm \ref{alg1} and how to extract a (sub)partition from the output of Algorithm \ref{alg1}.  We also detail how to (optionally) apply weights to vector entries.
Code is provided in Appendices \ref{code:thresh1} and \ref{code:thresh2}.

\subsection{Ordering the sparse vectors in terms of ``reliability''}
\label{sec:reliability}

Under Assumption \ref{assumptionA}, if Algorithm \ref{alg1} is performing well we should be able to find  a subspace $\mathcal{S}$ that is a good approximation of $\mathcal{V}$ and which is spanned by a nonnegative  sparse basis $s_1,\ldots,s_r$.
Nonnegativity of $s_1,\ldots,s_r$ is not guaranteed by Algorithm \ref{alg1} and we have found that if a sparse vector $s_j$ has one or more large negative entries, this is an indication that the feature $s_j$ has not been cleanly separated from other features.
We propose a ``reliability ordering'' of the sparse output vectors:  in Step 7 of Algorithm \ref{alg1} we order the vectors $s_1,\ldots,s_r$ so that $m_j:=\min_{1\le i\le p}s_{ij}$ is in descending order.
That is, $j$ for which $m_j=0$ means that $s_j$ should highlight a very reliably coherent feature and appear early in the ordering, while those $s_j$ toward the end of the ordering are potentially increasingly spurious.
In the case of the Bickley jet with $r=8$, we found $m_j=0$, $j=1,\ldots,6$, with the remaining $m_7,m_8>-0.02$;  see the eight sparse vectors shown in Figure \ref{fig:bickley_example} (right column), ordered by decreasing $m_j$.
In Section \ref{sec:negmass} we will adapt this heuristic to select the number of input eigenvectors. 

\subsection{Selecting an appropriate number of input vectors}
\label{sec:scales}
Consider a set $v_1,\ldots,v_r\in\mathbb{R}^p$ of linearly independent eigenvectors as in Section \ref{sec:sparsebasis}.
As one increases $r$, the $r^{\rm th}$ eigenvector becomes increasingly oscillatory, corresponding to increasingly rapidly decaying modes.
The ``spatial scale'' of the features encoded in the $r^{\rm th}$ eigenvector therefore decreases with increasing $r$.
Similarly, the $r^{\rm th}$ eigenvalue describes the exponential decay rate of the $r^{\rm th}$ eigenvector, and therefore a ``temporal scale'' for the corresponding features.
In the dynamical systems context, the number $r$ controls the smallest \emph{spatial scale} at which eigenvectors identify almost-invariant or coherent sets, while the $r^{\rm th}$ eigenvalue provides a lower bound on the \emph{temporal scale} at which each of the almost-invariant or coherent sets mix with the rest of phase space.
In particular, gaps in the spectrum indicate a jump in temporal scales, rather than a jump in spatial scales.
The reliability ordering of Section \ref{sec:reliability} will guide the choice of a number $k\le r$ to indicate that Algorithm \ref{alg1} has produced $k$ reliable features at the spatiotemporal scale governed by $r$.

\subsubsection{Heuristics based on the spectrum}
\label{sec:eigengap}
\begin{figure}[hbt]
\centering
\includegraphics[scale=0.5]{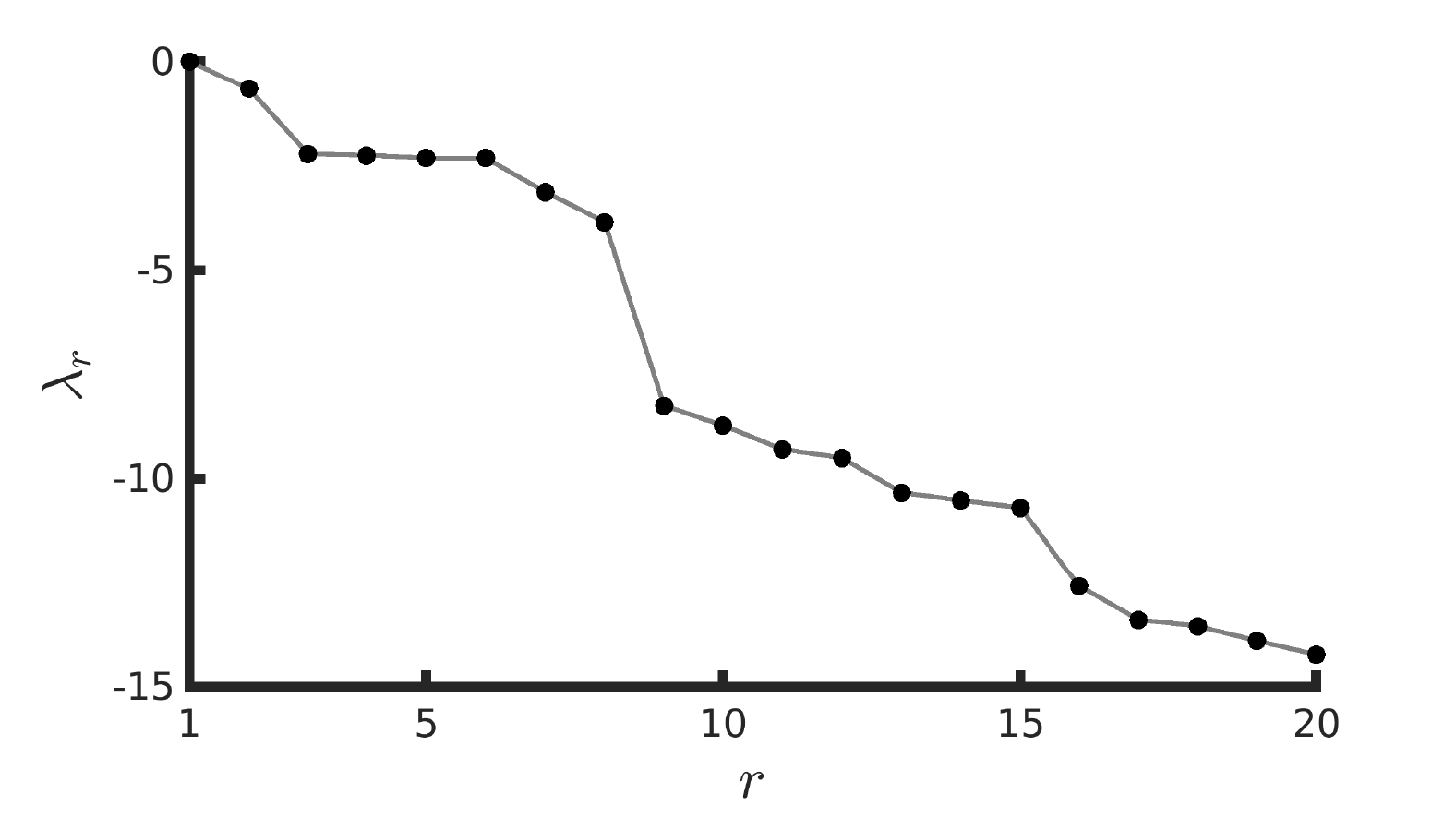}
\includegraphics[scale=0.5]{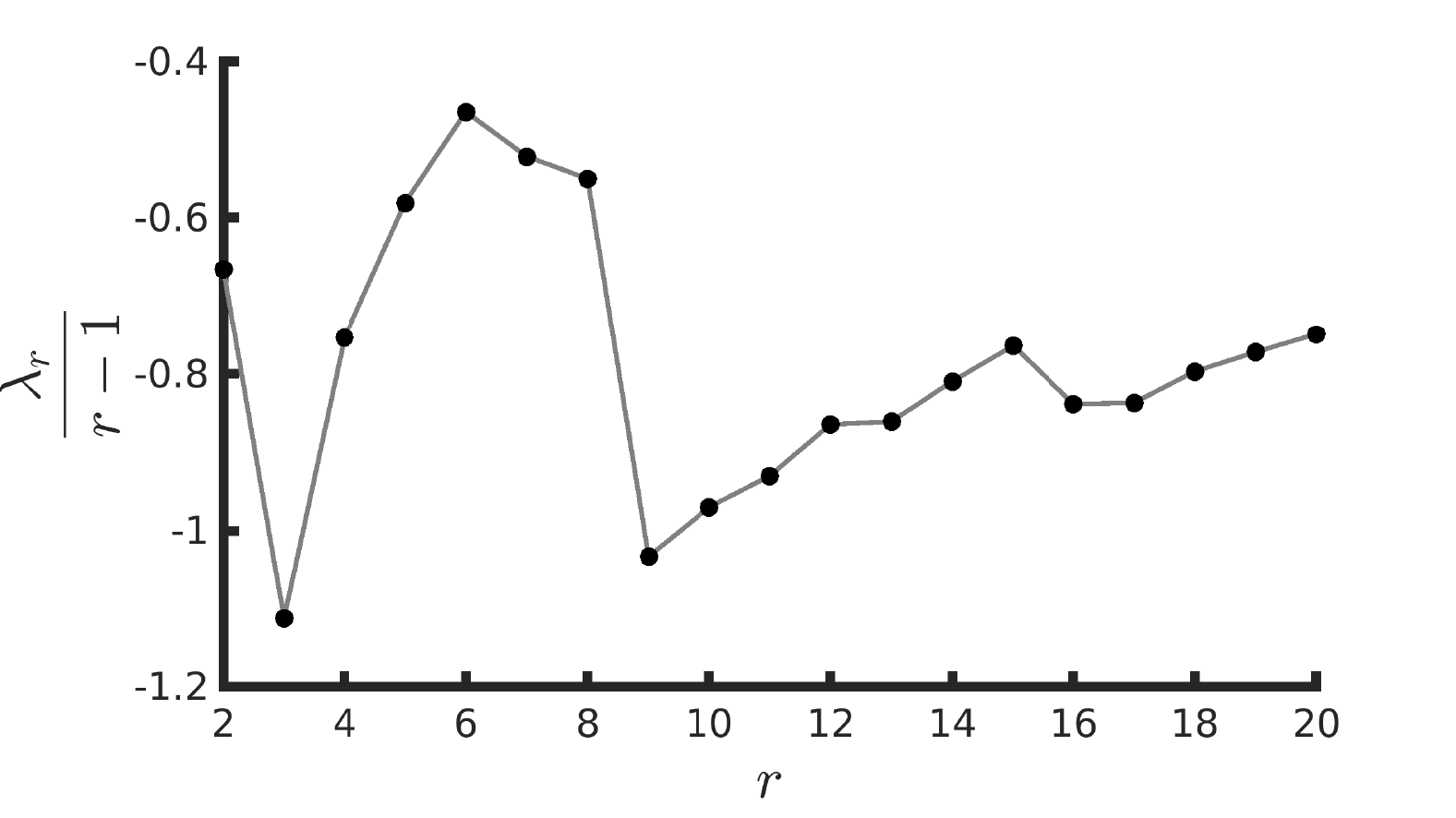}
\caption{Left: plot of $\lambda_r$ vs $r$ for $r=1,\ldots,20$ for the Bickley jet of Section \ref{sec:bickley}.
Right: plot of $\lambda_{r}/(r-1)$ vs $r$ for $r=2,\ldots,20$.  The eigengaps in the left image at $r=2$ and $r=8$ are highlighted even more strongly in the Weyl rescaling in the right image.  The eigengap at $r=15$ in the left image has been de-emphasised in the right image relative to $r=2$ and $r=8$.}
\label{eigVal_Bickley}
\end{figure}

In spectral clustering, the eigengap heuristic (see \cite{von_luxburg07} and references therein) is a common method for selecting an appropriate number of vectors for further analysis.
The eigengap heuristic suggests to look for large gaps in the spectrum, i.e.\ if $|\lambda_{r+1}-\lambda_r|$ is large compared to $|\lambda_{i+1}-\lambda_i|$ for $i=1,\ldots,r-1$, one should truncate the collection of eigenvectors to $v_1,\ldots,v_r$.
A potential issue with this rule (apart from the absence of an unambiguous spectral gap) is that depending on the dimension of the underlying manifold, the natural asymptotic behaviour of $\lambda_k$ as a function of $k$ given by Weyl's law \cite{weyl} is not necessarily linear.
We therefore propose a new variant of the eigengap heuristic by taking this asymptotic into account.
If our phase space is a $d$-dimensional Riemannian manifold, denote by $N(\lambda)$ the number of eigenvalues of the Laplace(-Beltrami) operator $\Delta$ on $M$ no greater than $\lambda$ (applying either homogeneous Neumann or Dirichlet boundary conditions).
Weyl's law states that $\lim_{\lambda\to\infty}N(\lambda)/\lambda^{d/2}=-(2\pi)^{-d}\omega_d|M|$ where $\omega_d$ is the volume of the $d$-dimensional unit ball and $|M|$ is the volume of $M$.
Thus, ordering the eigenvalues of $\Delta$ as $0\ge \lambda_1\ge \lambda_2\ge\cdots$, we have $\lambda_r\sim -Cr^{2/d}$ for a constant $C$ as $r\to\infty$.
\paragraph{Laplace-type matrix/operator, Neumann boundary conditions:}
For example, a graph Laplacian, a Laplace-Beltrami operator, or the dynamic Laplacian\footnote{which is not necessarily a Laplace-Beltrami operator \cite{karrasch17}.} \cite{F15,FK17} of Sections \ref{sec:bickley} and \ref{sec:northatlantic} with homogeneous Neumann boundary conditions.
In this setting, the leading eigenvalue is $\lambda_1=0$ and fitting the asymptotic $\lambda_r\sim -Cr^{2/d}$ to this condition suggests plotting $\lambda_{r}/(r-1)^{2/d}$ vs.\ $r$ for $r=2,\ldots,$ and looking for the largest drops from $r$ to $r+1$; see Figure \ref{eigVal_Bickley}.

\paragraph{Laplace-type operator, Dirichlet boundary conditions:}
For example, a Laplace-Beltrami operator or the dynamic Laplacian with Dirichlet boundary conditions.
The leading eigenvalue $\lambda_1$ is strictly negative (and unknown) so we plot $\lambda_{r}/r^{2/d}$ vs.\ $r$ for $r=1,\ldots,$ and look for the largest drops from $r$ to $r+1$.

\paragraph{Markov-type operator:} For example, a matrix or operator arising from diffusion maps \cite{banisch17}, the various transfer operator constructions in Table \ref{tab:TO}, or the normalised finite-time transfer operator \cite{FSM10,F13} of Section \ref{sec:turbulence}.  In these examples $\lambda_1=1$ and the rest of the spectrum is contained in the unit circle in the complex plane.  Taking logs and retaining only the real part, we obtain a Laplace-type spectrum with Neumann boundary conditions.  Thus, we plot $\Re(\log\lambda_{r})/(r-1)^{2/d}$ vs.\ $r$ for $r=2,\ldots,$ and look for the largest drops from $r$ to $r+1$. In the case of transfer operators arising from continuous-time flows, this logarithmic transformation can be made rigorous via the spectral mapping theorem, see \cite{FJK13,FK17}.

The above heuristics are applied to spectra arising from both transfer operators and the dynamic Laplacian in Section \ref{sec:examples}.


\subsubsection{A heuristic based on sparse vectors}
\label{sec:negmass}

Following the arguments of Section \ref{sec:reliability}, for the $p\times r$ matrix $S$ produced by Algorithm \ref{alg1}, we create a cumulative minimum value quantity $\Min(S):=\sum_{j=1}^r -m_j=-\sum_{j=1}^r \min_{1\le i\le p}S_{ij}$.
Denoting by $S^{(r)}$ the sparse array produced by Algorithm \ref{alg1} with $r$ input vectors, one may then plot $\Min(S^{(r)})$ vs.\ $r$ and select those $r$ for which there has been a large drop from $r-1$;  i.e.\ $\Min(S^{(r)})-\Min(S^{(r-1)})$ is negative.
The rationale for this choice is that $\Min(S^{(r)})$ is typically increasing with $r$ because it is a sum of non-negative values.
If there is a drop from $r-1$ to $r$ it means that despite adding an extra term to the sum, the sum decreases, and thus overall the new $m_j$ values (based on $r$ input vectors) are smaller than the old $m_j$ (based on $r-1$ input vectors), indicating that the new sparse basis $s_1,\ldots,s_r$ has better separated the $r$ features.
Note that such a heuristic implicitly assumes that \emph{all} $r$ sparse vectors are reliable, and no sparse vector will be rejected from the collection.
In Section \ref{sec:turbulence} we will combine this heuristic with the reliability heuristic of Section \ref{sec:reliability} to suggest a way to simultaneously select both $r$ (the number of input vectors), and $k$ (the number of reliable sparse vectors), with $k< r$ after some sparse vectors have been rejected.


\subsection{Extraction of a (sub)-partition from a sparse basis output}
\label{sec:subpart}
Recall that we wish to separate the features encoded in the eigenvectors $v_1,\ldots,v_r$ resulting from a spectral clustering.
Algorithm \ref{alg1} has already done most of the hard work by forming a sparse approximate basis of the space spanned by the eigenvectors.
Indeed one would often be very satisfied with this separation procedure.
A simple and effective way of summarising the separated features is to create a ``superposition vector''.
\paragraph{Superposition vector:}  $\mathfrak{s}:=\min\{1,\sum_{j=1}^r \max\{S_{ij},0\}\}$.  This vector combines membership likelihoods across features, and may be interpreted as the likelihood that the $i^{\rm th}$ index belongs to \emph{some} feature.
In MATLAB, this is realised as \verb"super=min(1,sum(max(S,0),2))";  superposition vectors $\mathfrak{s}$ are  illustrated in Figures \ref{fig:turbulence_eigenvector_sparse}, \ref{fig:turbulence_sparse_r30_k7_sparse}, \ref{fig:turbulence_sparse_r38_k21_sparse}, and \ref{fig:nat_sparse_r96_k36}.

If a hard separation is required, then a final step is to threshold the sparse vectors to fix which indices are in or out of a feature.
If there is a two-dimensional phase space associated to the sparse vectors,
then one may apply by inspection a manual threshold $\tau$ so that feature $j$ consists of grid cells or nodes with index $i$ satisfying $S_{ij}>\tau$.

For situations where visualisation is not possible, the main thing we desire is to ensure that distinct features are disjoint.
We consider three algorithms. Algorithm \ref{thresh1} applies the minimum threshold $\tau^{pu}$ to ensure that after thresholding, the columns $S_j$ form a sub-partition of unity: that is $\sum_{j=1}^r S_{ij} \le 1$ for each $i=1,\ldots,p$.
Algorithm \ref{thresh2} selects the minimum threshold $\tau^{dp}$ to ensure that after thresholding, the $s_j$ have disjoint support.
In each case, the principle of maximum likelihood is then applied to produce a feature vector $a \in \{0,\ldots,r\}^p$ with $a_i=j$ if element $i$ belongs to feature $j$ for $i=1\ldots,p$ and $j=1,\ldots,r$, and $a_i=0$ if feature $i$ is unassigned.
A third algorithm, described at the end of this section, simply applies the above maximum likelihood procedure without thresholding.
Define a thresholding transformation $H_\mu:\mathbb{R}\to\mathbb{R}$ by $H_\mu(z)=z$ if $|z|>\mu$ and $H_\mu(z)=0$ otherwise; $H_\mu$ is applied to vectors elementwise.

\begin{algorithm}[Input: vectors $\{s_1,\ldots,s_r\}\subset \mathbb{R}^p$ produced by Algorithm \ref{alg1}; output: thresholded $\{s_1,\ldots,s_r\}$ vectors and feature vector $a \in \{0,\ldots,r\}^p$]
\label{thresh1}\
\begin{enumerate}
\item For $j=1,\ldots,r$, set $s_j \to \max\{s_j, 0\}$ to make the vectors nonnegative.
\item For each row $i=1,\ldots,p$, let $\tilde{s}_{i1},\ldots,\tilde{s}_{ir}$ be the values of $s_{i1},\ldots,s_{ir}$ in decreasing order.
\item Set the threshold $\tau^{pu}:=\max_{1\le i\le p,1\le j \le r} \{\tilde{s}_{ij}:\sum_{l=1}^j \tilde{s}_{il}>1\}$. For each $j=1,\ldots,r$, hard threshold $s_j \to H_{\tau^{pu}}(s_j)$.
\item For each $i=1,\ldots,p$ let $j^*=\argmax_{1\le j\le r} s_{ij}$ (resolving ties arbitrarily).
If $s_{ij^*}>0$, set $a_i=j^*$,
otherwise set $a_i=0$.
\end{enumerate}
\end{algorithm}

\begin{algorithm}[Input: vectors $\{s_1,\ldots,s_r\}\subset \mathbb{R}^p$ produced by Algorithm \ref{alg1}; output:  thresholded $\{s_1,\ldots,s_r\}$ vectors and feature vector $a \in \{0,\ldots,r\}^p$] \label{thresh2}\
\begin{enumerate}
\item For $j=1,\ldots,r$, set $s_j \to \max\{s_j, 0\}$ to make the vectors nonnegative.
\item For each row $i=1,\ldots,p$, let $\tilde{s}_{i1},\ldots,\tilde{s}_{ir}$ be the values of $s_{i1},\ldots,s_{ir}$ in decreasing order.
\item Set the threshold $\tau^{dp}:=\max_{1\le i \le p} \tilde{s}_{i2}$, and for each $j=1,\ldots,r$ hard threshold $s_j \to H_{\tau^{dp}}(s_j)$.
\item For each $i=1,\ldots,p$, $j=1,\ldots,r$, if any $s_{ij}>0$ for $j=1,\ldots,r$, set $a_i=j$, otherwise, set $a_i=0$.
\end{enumerate}
\end{algorithm}

\begin{figure}[hbt]
\centering
\includegraphics[width=0.45\textwidth]{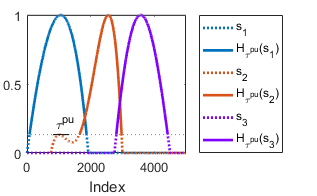}
\includegraphics[width=0.45\textwidth]{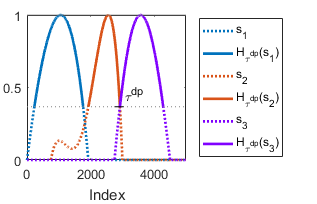}
\includegraphics[width=0.45\textwidth]{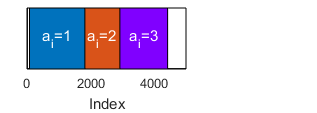}
\includegraphics[width=0.45\textwidth]{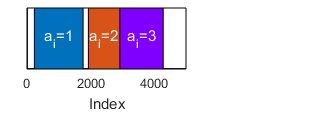}
\caption{
Upper left: $\tau^{pu}$ from Algorithm \ref{thresh1} is indicated by the horizontal dotted line.
$H_{\tau^{pu}}(s_1),H_{\tau^{pu}}(s_2),H_{\tau^{pu}}(s_3)$ (solid lines) are equal to $s_1,s_2,s_3$ wherever those vectors are greater than $\tau^{pu}$, and zero elsewhere.
Note that $H_{\tau^{pu}}(s_1)+H_{\tau^{pu}}(s_2)+H_{\tau^{pu}}(s_3)\le 1$ as required for a partition of unity.
Lower left: The resulting partition $a$, obtained by performing maximum likelihood on $H_{\tau^{pu}}(s_1),H_{\tau^{pu}}(s_2),H_{\tau^{pu}}(s_3)$.
Upper right: $\tau^{dp}$ from Algorithm \ref{thresh2} is indicated by the horizontal dotted line.
$H_{\tau^{dp}}(s_1),H_{\tau^{dp}}(s_2),H_{\tau^{dp}}(s_3)$ are  plotted as solid lines.
Lower right: The resulting partition $a$ obtained by  performing maximum likelihood on $H_{\tau^{dp}}(s_1),H_{\tau^{dp}}(s_2),H_{\tau^{dp}}(s_3)$.}
\label{fig:thres_example}
\end{figure}

Because a disjoint partition is a stronger requirement than a partition of unity, one has $\tau^{dp}\ge \tau^{pu}$.
If $\tau^{dp} \ge 0.5$, then $\tau^{dp}=\tau^{pu}$.
We illustrate both algorithms in Figure \ref{fig:thres_example}, where three synthetic sparse vectors $s_1,s_2,s_3\in \mathbb{R}^{5000}$ are plotted as dotted curves.
These vectors have overlapping supports, and $s_1+s_2>1$ on part of the domain.
The post-thresholding vectors are shown as solid curves in Figure \ref{fig:thres_example}.

One is free to apply other thresholdings to achieve specific outcomes;  the above two algorithms are of generic use for feature separation.
Thresholds may be applied vectorwise to individual vectors, or they may be computed from individual vectors and applied uniformly across the vectors as in the two algorithms above.
In Section \ref{sec:turbulence} we illustrate another option specifically related to coherent sets in nonlinear dynamical systems.

If no thresholding method is effective for a particular system (for example, when feature separations are not clear), maximum likelihood can be applied directly to the output of Algorithm \ref{alg1};  that is, one can apply the last step of Algorithm \ref{thresh1} to obtain a feature vector $a$.
This can be realised in MATLAB by \verb"[m,a]=max(S,[],2); a(m<=0)=0".
The support of $a$ coincides with the positive support of $\max_j s_{ij}$.

\subsection{Weighting entries of the data vectors}
\label{sec:weights}
In certain situations, one may wish to apply a weight vector $0< \nu \in \mathbb{R}^p$ to the calculation of the objective in \eqref{eq:spcart}. 
The weight vector $\nu$ allows one to emphasise (or deemphasise) specific indices $i=1,\ldots,p$.
One situation where this may be appropriate is if one is constructing a transfer operator from a grid with differently sized grid cells;  larger cells can be weighted according to their area or volume.

We define a weighted inner product $\langle v,w\rangle_\nu:=\sum_{i=1}^p \nu_iv_iw_i$, denote the corresponding norm by $\|\cdot\|_{\ell_{2,\nu}}$, and denote a weighted Frobenius norm by $\|A\|_{F,\nu}:=\sqrt{\sum_{i=1}^p\sum_{j=1}^r \nu_iA_{ij}^2}$.
Let $\mathfrak{U}^{p,r}_\nu=\{A\in \mathbb{R}^{p\times r}:{}$ each column of $A$ has $\ell_{2,\nu}$ norm $1\}$.
We define a weighted $\ell_{1,1}$ matrix norm by $\|A\|_{1,1,\nu}:=\sum_{i=1}^p\sum_{j=1}^r \nu_i|A_{ij}|$.
Let $V \in \mathbb{R}^{p \times r}$ be a matrix whose columns are pairwise orthonormal in the inner product $\langle\cdot,\cdot\rangle_\nu$.
Problem \eqref{eq:spcart} can be generalised to
\begin{align}
\argmin_{\substack{S\in\mathfrak{U}_\nu^{p,r} \\ R \in \mathfrak{S}^{r}}} \frac{1}{2} \|V-SR\|_{F,\nu}^2 + \mu \|S\|_{1,1,\nu}. \label{eq:weighted_spcart}
\end{align}
Denoting by $D_\nu$ the diagonal matrix with $\nu$ on the diagonal,
we have that $\|A\|_{1,1,\nu}=\|D_\nu A\|_{1,1}$ and $\|A\|_{F,\nu}=\|D_\nu^{1/2} A\|_F$.
Substituting $V'=D_\nu^{1/2}V$ and $S'=D_\nu^{1/2}S$,  problem \eqref{eq:weighted_spcart} is equivalent to
\begin{align}
\argmin_{\substack{S'\in\mathfrak{U}^{p,r} \\ R \in \mathfrak{S}^{r}}} \frac{1}{2} \|V'-S'R\|_F^2 + \mu \left\|D_\nu^{\frac{1}{2}}S'\right\|_{1,1}.
\label{eq:diagonal_spcart}
\end{align}

We generalise the soft threshold function $C_\mu$ to $C'_{\mu,\nu}:\mathbb{R}^p\to\mathbb{R}^p$, defined by $[C'_{\mu,\nu}(v)]_i=C_{\mu \nu_i^{1/2}}(v_i)$.
We prove in Appendix \ref{app:proof} that Problem \eqref{eq:diagonal_spcart} can be solved using Algorithm \ref{alg1} with $V$ replaced with $V'$, $\mu=0.99/\sqrt{p}$ replaced with $\mu=0.99/\sqrt{\|\nu\|_1}$ in Step \ref{itm:setmu}, $C_\mu$ replaced with $C'_{\mu,\nu}$ in Step 2, and Step \ref{itm:withsign} replaced with $S_j\to\sign(\nu^\top S_j) S_j$.
After running Algorithm \ref{alg1}, set $S=D_\nu^{-1/2}S'$ to obtain a solution for Problem \eqref{eq:weighted_spcart}. The parameter $\mu$ should be chosen less than $1/\sqrt{\|\nu\|_1}$ to ensure  $C'_{\mu,\nu}(D_\nu^{1/2}v)$ is nonzero for every $v \in \mathbb{R}^p$ with $\|v\|_{2,\nu}=1$; taking $\mu\ge 1/\sqrt{\|\nu\|_1}$ will result in $C'_{\mu,\nu}(D_\nu^{1/2}v)$ being the zero vector when $v$ is the constant $\ell_{2,\nu}$-unit vector.

\section{Examples} \label{sec:examples}
We perform numerical experiments on three systems: the Bickley jet, a turbulence simulation obtained from the Navier-Stokes equations, and North Atlantic ocean currents derived from satellite altimetry. We use the transfer operator (Construction 5 in Table \ref{tab:TO}) and the dynamic Laplacian (Construction 6 in Table \ref{tab:TO}) to detect coherent sets.

We can quantify the goodness of fit and sparsity of our solutions using the corresponding parts of the objective \eqref{eq:spcart} being minimised.
\begin{itemize}
\item \textbf{$\ell_2$ subspace error}: For subspace approximation, note that for $p\times r$ orthogonal $V$, one has $\|V\|_F^2=r$. Thus, reporting $\frac{1}{r}\|V-SR\|_F^2$ quantifies the subspace approximation error relative to the ``zero'' estimate $S=0$; this value should lie between 0 and 1 (lower value means better approximation).

\item \textbf{Absolute sparsity}: We report the absolute sparsity as $\|S\|_{0,1}/(p r)$, where $\|S\|_{0,1}:=\sum_{i=1}^p\|S_i\|_0$.
The absolute sparsity is the proportion of elements of $S$ that are nonzero, and therefore  
lies between zero and one (lower value means sparser).
\item \textbf{Relative  sparsity}: 
We report the relative ($\ell_1$) sparsity $\|S\|_{1,1}/\|V\|_{1,1}$. The $\|\cdot\|_{1,1}$ norm is used as the sparsity penalty term in Algorithm \ref{alg1}. The relative sparsity should lie between 0 and 1 (lower value means sparser).
\end{itemize}

In Sections \ref{sec:bickley} and \ref{sec:northatlantic} all time integrations are carried out using an explicit Runge-Kutta method of order (4,5) with adaptive step size, as implemented by MATLAB's \verb"ode45" integrator, using the absolute and relative error tolerance $10^{-3}$. The methodologies listed in Table \ref{tab:TO} apply to any finite dimension, we restrict ourselves to 2-dimensional examples to present our results with greater transparency.

\subsection{Bickley jet} \label{sec:bickley}
The Bickley jet is an idealised system introduced as a model of ``banded chaos'' by \cite{del-castillo-negrete93}, which has become a popular test case for coherent set detection methods \cite{rypina07,williams15,hadjighasem16,schlueter-kuck17,padberg-gehle17,Hadjighasem_etal,FJ18}. The Bickley jet is a Hamiltonian system modeling a meandering jet with vortices on either side. Its stream function is
\begin{align}
\psi(x,y,t) &= -U_0L_0\tanh(y/L_0) + \sum_{i=1}^3 A_i U_0 L_0 \sech^2(y/L)cos(k_i(x-c_it)).
\end{align}
We use parameter values $U_0 = 62.66$ m/s, $L_0 = 1770$ km, $A_1 = 0.0075$, $A_2 = 0.15$, $A_3 = 0.3$, $c_1 = 0.1446U_0$, $c_2 = 0.205U_0$, $c_3 = 0.461U_0$, $k_1 = 2/r_e$; $k_2 = 4/r_e$, $k_3 = 6/r_e$, $r_e = 6371$ km as in \cite{Hadjighasem_etal}. We consider the associated flow on the initial domain $\mathcal{M}=[0,20] \times [-3, 3]$, periodic in the $x$-direction, over the time interval $[0,40]$ days.
We calculate the dynamic Laplacian using the ``adaptive transfer operator method'' of \cite{FJ18} on a Delaunay triangulation (alpha complex with $\alpha=0.06$) of $\mathcal{M}$ with nodes on an initial grid of $600\times 180$ points, on the time set $\mathcal{T}=\{0, 40\}$.
There are clear eigengaps after the second and eighth eigenvalues, in both the unnormalised and Weyl-normalised plots (Figure \ref{eigVal_Bickley}). These two eigengaps characterise two natural spatiotemporal scales for the Bickley dynamics.

\begin{figure}[hbt]
\centering
\begin{subfigure}[b]{0.45\textwidth}
\includegraphics[width=\textwidth]{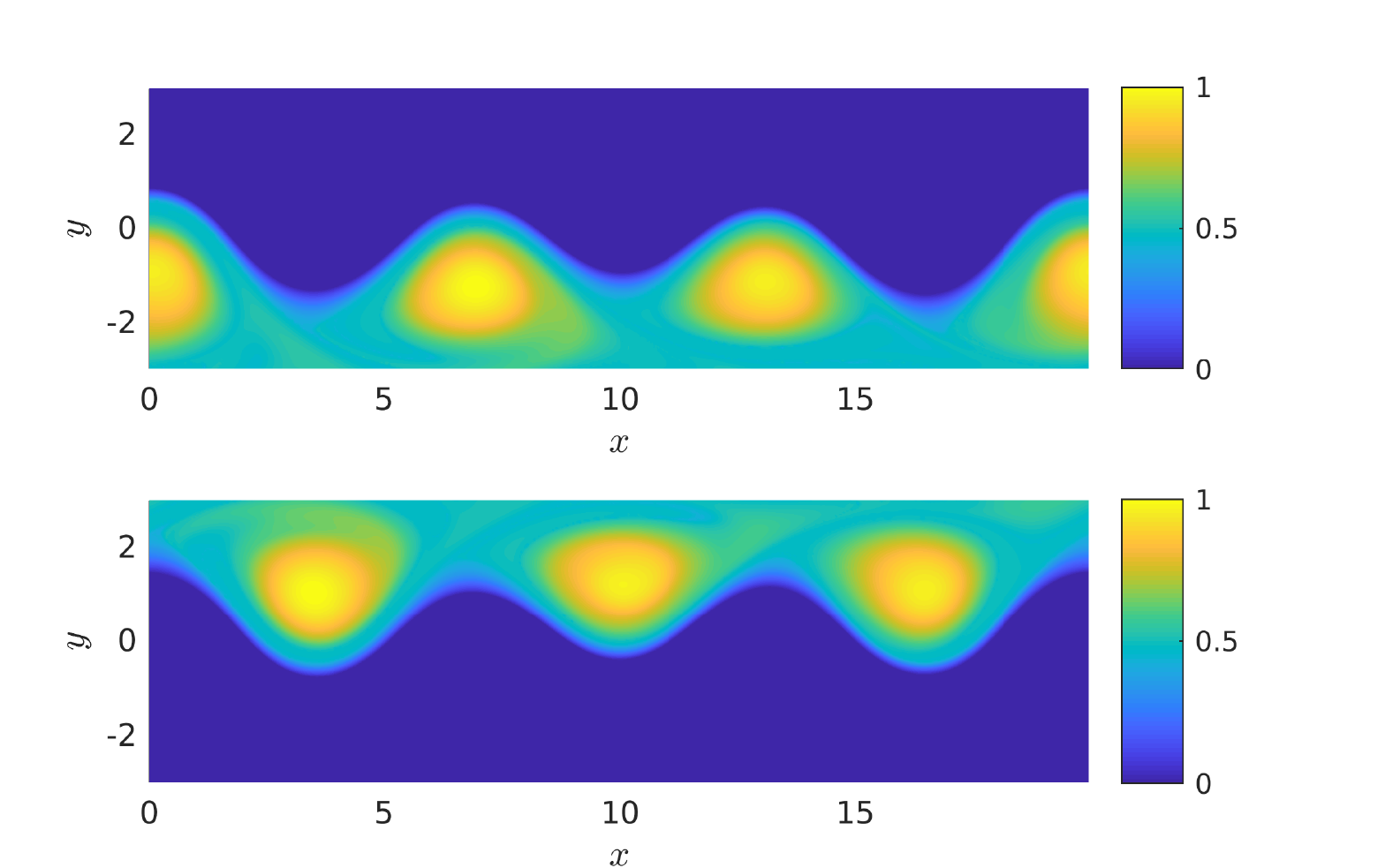}
\vspace{1.5\bigskipamount}
\caption{Sparse vectors approximating the subspace spanned by the two leading eigenvectors of the dynamic Laplacian ($r=2$).}
\label{fig:bickley_SpEigFunc_r2}
\end{subfigure}\qquad
\begin{subfigure}[b]{0.45\textwidth}
\includegraphics[width=\textwidth]{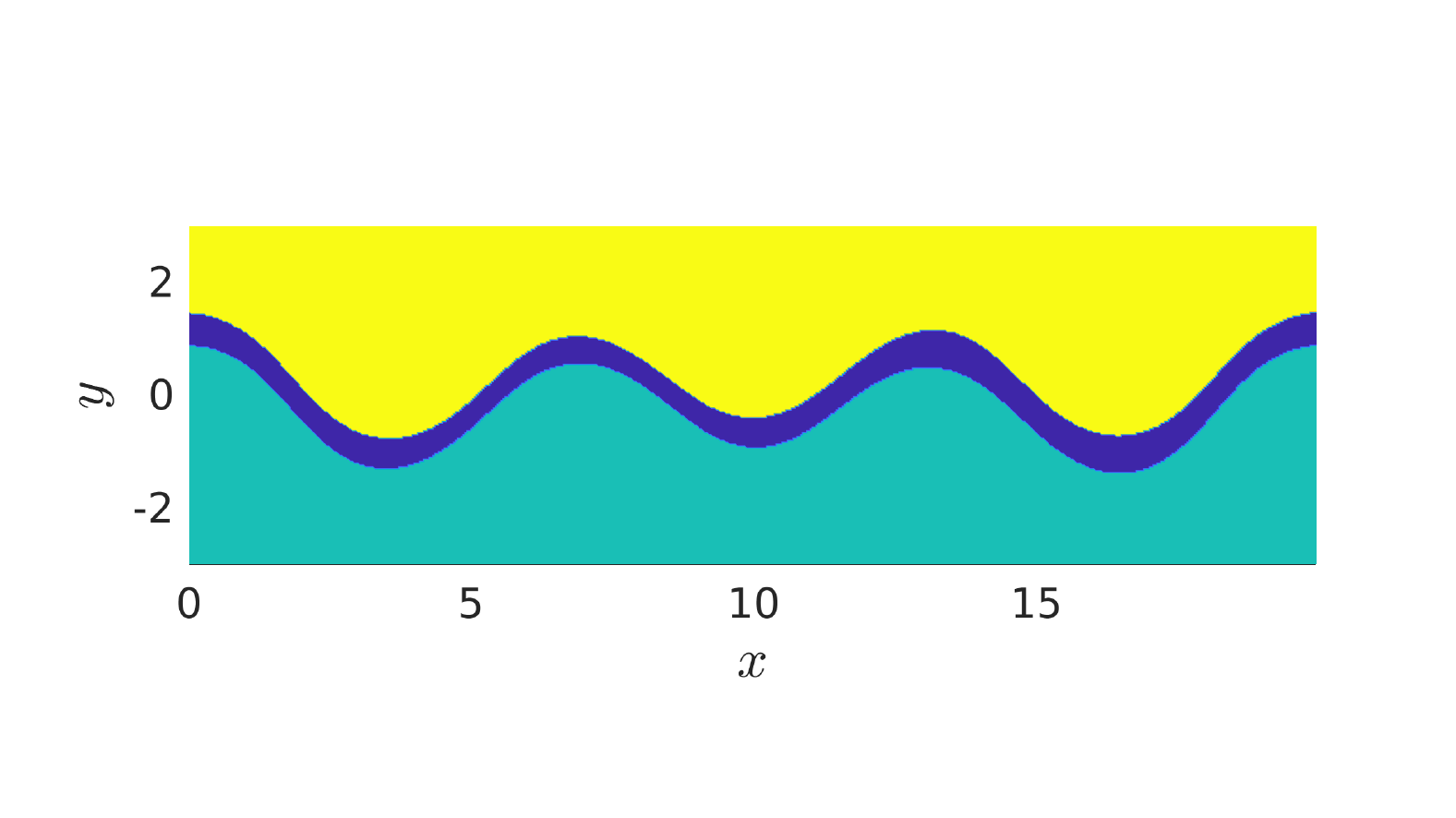}
\caption{Hard thresholding result from Algorithms \ref{thresh1} and \ref{thresh2}.}\label{fig:bickley_SPCA_Partition_r2}
\includegraphics[width=\textwidth]{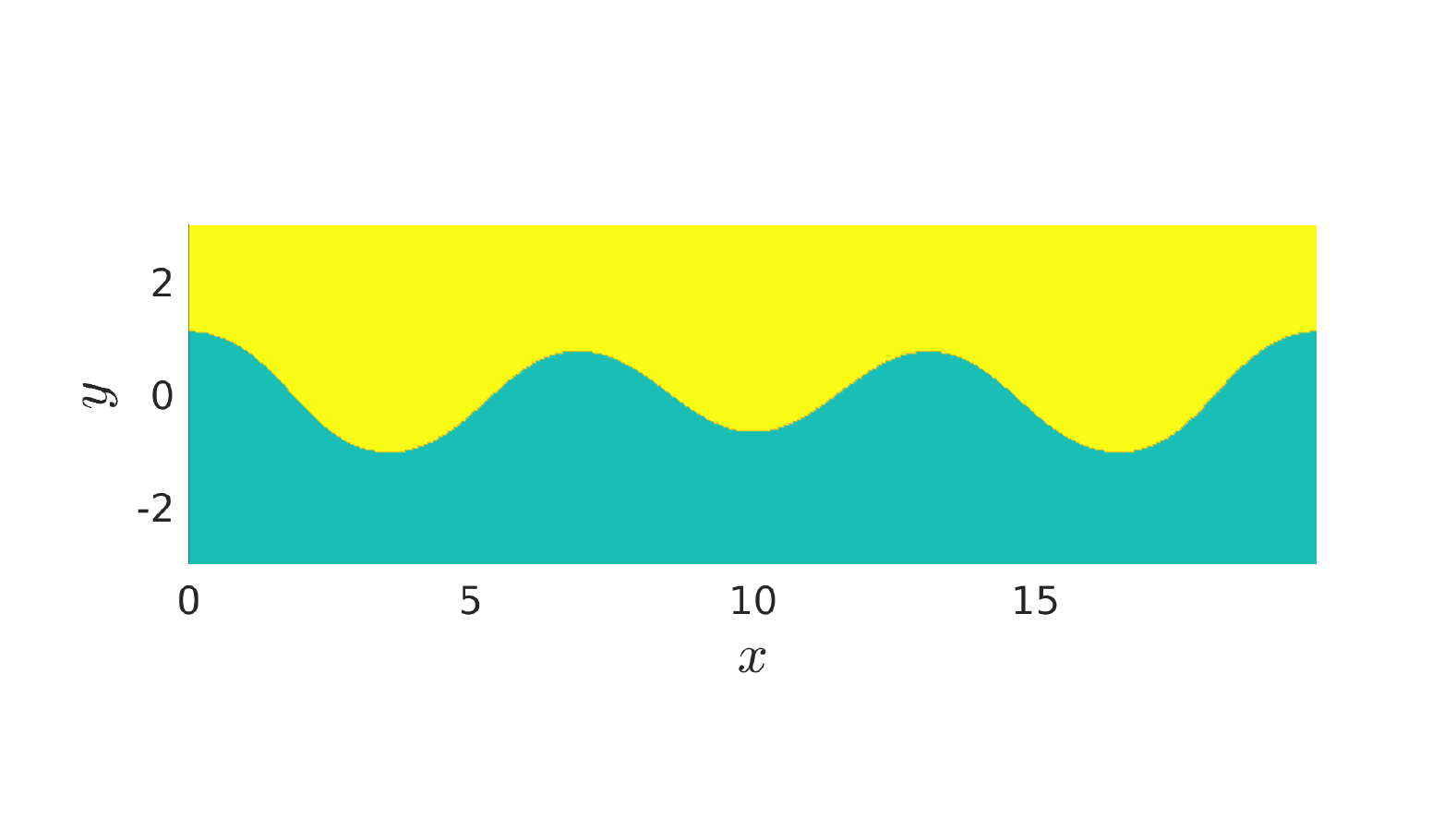}
\caption{Result of k-means (requesting two clusters) applied to the two leading eigenvectors of the dynamic Laplacian.}
\label{fig:bickley_KMeans_r2}
\end{subfigure}
\caption{Bickley jet, with $r=k=2$.}
\end{figure}

We applied Algorithm \ref{alg1} with $r=2$ and $r=8$ to the eigenvectors of the dynamic Laplacian.
With $r=2$, the two eigenvectors are the leading (constant) eigenvector and a second non-constant eigenvector.
Algorithm \ref{alg1} produces two sparse vectors whose supports separate the domain into upper and lower pieces;  see Figure \ref{fig:bickley_SpEigFunc_r2}.
The six vortices are also highlighted within the supports of the two sparse vectors.
The span of these two sparse vectors is a reasonable approximation of the span of the two leading eigenvectors (subspace error of 11.8\%).
The sparse vectors have an absolute sparsity of 45.2\% and a relative sparsity of 65.0\%.

We apply Algorithm \ref{thresh1} and \ref{thresh2} to obtain a partition of the domain $A_1 \sqcup A_2 \subset M$. Since the sparse vectors in Figure \ref{fig:bickley_SpEigFunc_r2} have disjoint support, both algorithms return the same sub-partition, formed by the supports of the two sparse vectors;  see Figure \ref{fig:bickley_SPCA_Partition_r2}.
Note that the central jet between the upper and lower halves is also visible (zeroed as dark blue), using only information from the first nontrivial dynamic Laplacian eigenvector.
k-means clustering creates a partition of the entire phase space; see
Figure \ref{fig:bickley_KMeans_r2}.

With $r=8$, Algorithm \ref{alg1} produces the sparse vectors in Figure \ref{fig:bickley_example} (right column).
The subspace error is low at 6.4\%, as are the absolute sparsity and relative sparsity at 14.1\%  and 44.6\%, respectively.
The low absolute sparsity indicates that we have effectively separated the domain into disjoint features, with little overlap between the features.

\begin{figure}[hbt]
\centering
\begin{subfigure}{0.3\textwidth}
\includegraphics[width=\textwidth]{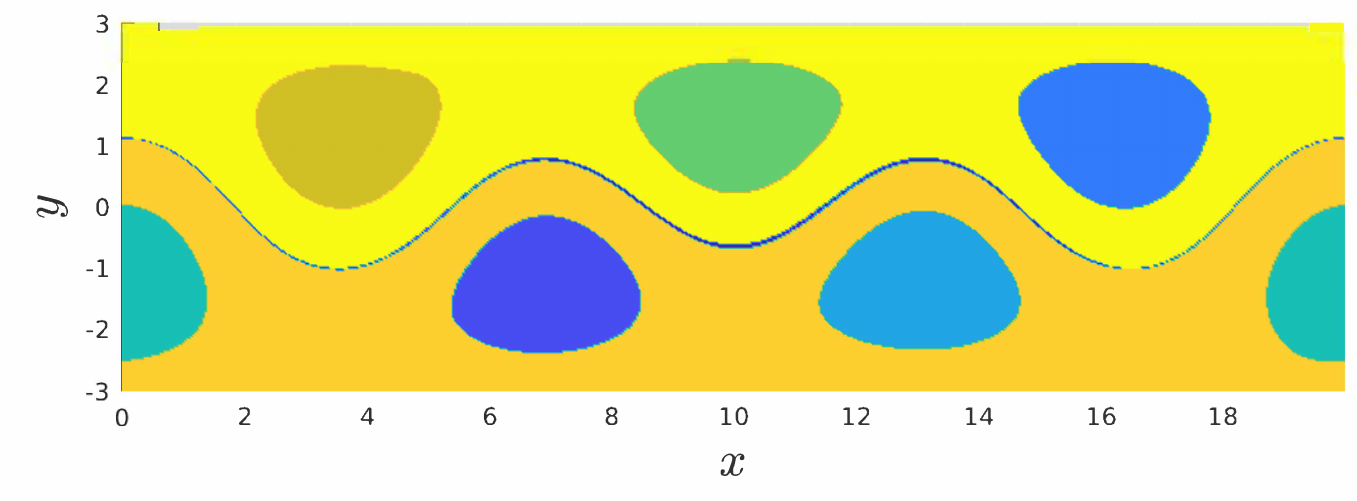}
\caption{Algorithm \ref{thresh1} } \label{fig:bickley_SPCA_Partition_r8_tau1}
\end{subfigure}
\begin{subfigure}{0.3\textwidth}
\includegraphics[width=\textwidth]{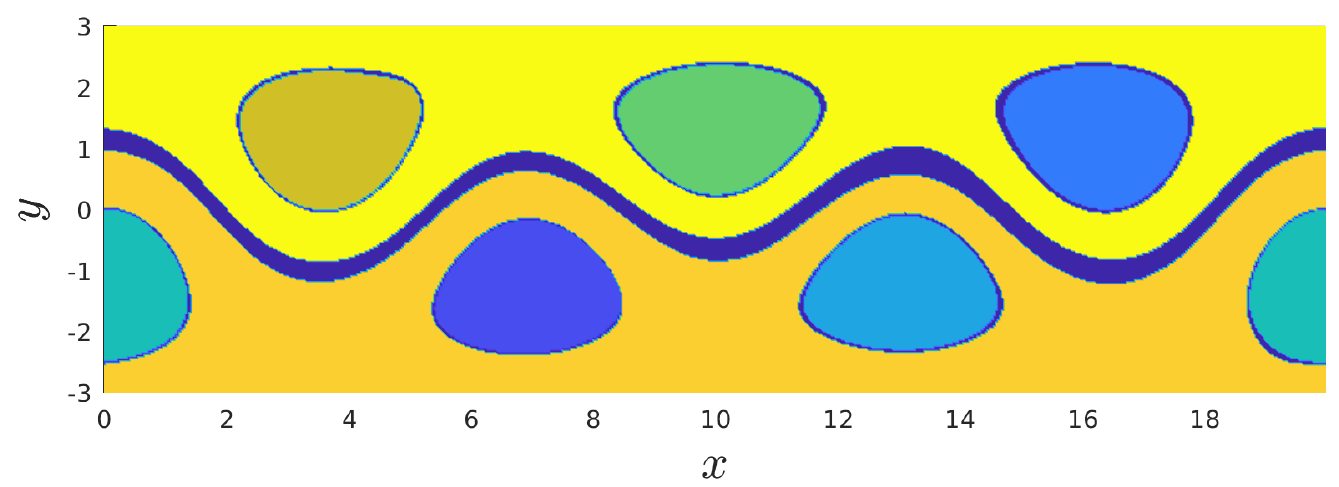}
\caption{Algorithm \ref{thresh2} } \label{fig:bickley_SPCA_Partition_r8_tau2}
\end{subfigure}
\begin{subfigure}{0.35\textwidth}
\includegraphics[width=\textwidth]{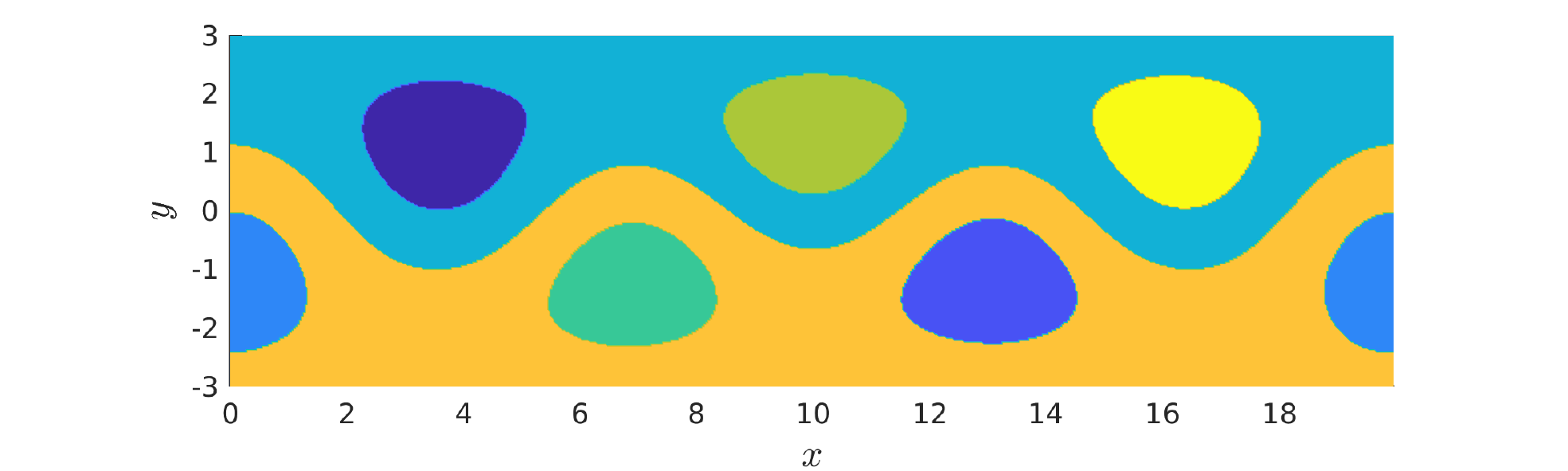}
\caption{k-means} \label{fig:bickley_KMeans_r8}
\end{subfigure}
\caption{State space partition of Bickley dynamics using dynamic Laplacian eigenvectors and Algorithm \ref{alg1} ($r=8$), followed by: (a) Algorithm \ref{thresh1};  (b)  Algorithm \ref{thresh2};  (c) k-means with 8 clusters.
} \label{fig:bickley_SPCA_Partition_r8}
\end{figure}

We can again use either Algorithm \ref{thresh1} or Algorithm \ref{thresh2} to obtain a partition of the domain;
see Figures \ref{fig:bickley_SPCA_Partition_r8_tau1} and \ref{fig:bickley_SPCA_Partition_r8_tau2} respectively.
Algorithm \ref{thresh1} produces a sub-partition of the domain consisting of almost the entire domain, while Algorithm \ref{thresh2} produces a sub-partition of the domain with a larger gap between the upper and lower regions corresponding to a central jet.
Both of these results are similar to the k-means clustering result with $k=8$ clusters, as shown in Figure \ref{fig:bickley_KMeans_r8}.

\begin{remark}
The k-means algorithm can be viewed as a coarse type of sparse eigenbasis approximation in the following way.
Apply k-means to the vectors $v_1,\ldots,v_r$ arising from a spectral clustering, requesting $r$ clusters.
For $j=1,\ldots,r$ define vectors $\hat{s}_{j}\in\mathbb{R}^p$ as $\hat{s}_{ij}=1$ if index $i$ is contained in cluster $j$, and $\hat{s}_{ij}=0$ otherwise.
In the ``ideal'' situation where data points are very clearly separated into distinct clusters, the eigenvectors $v_1,\ldots,v_r$ will be approximately constant, and the (sparse) vectors $\hat{s}_1,\ldots,\hat{s}_r$ will approximately span the same space.
By solving a Procrustes problem, one can also infer an orthogonal rotation implicitly performed by k-means to approximately rotate the vectors $v_1,\ldots,v_r$ to the vectors $\hat{s}_1,\ldots,\hat{s}_r$.
This relationship with sparse eigenbasis approximation can help explain why k-means can be very effective in some (closer to ``ideal'') cases where one expects to classify \emph{every} data point (i.e.\ one \emph{partitions} the dataset into features that exhaust the dataset).
\end{remark}

\subsection{Turbulent flow} \label{sec:turbulence}
\begin{figure}[hbt]
\centering
\includegraphics[width=0.325\textwidth]{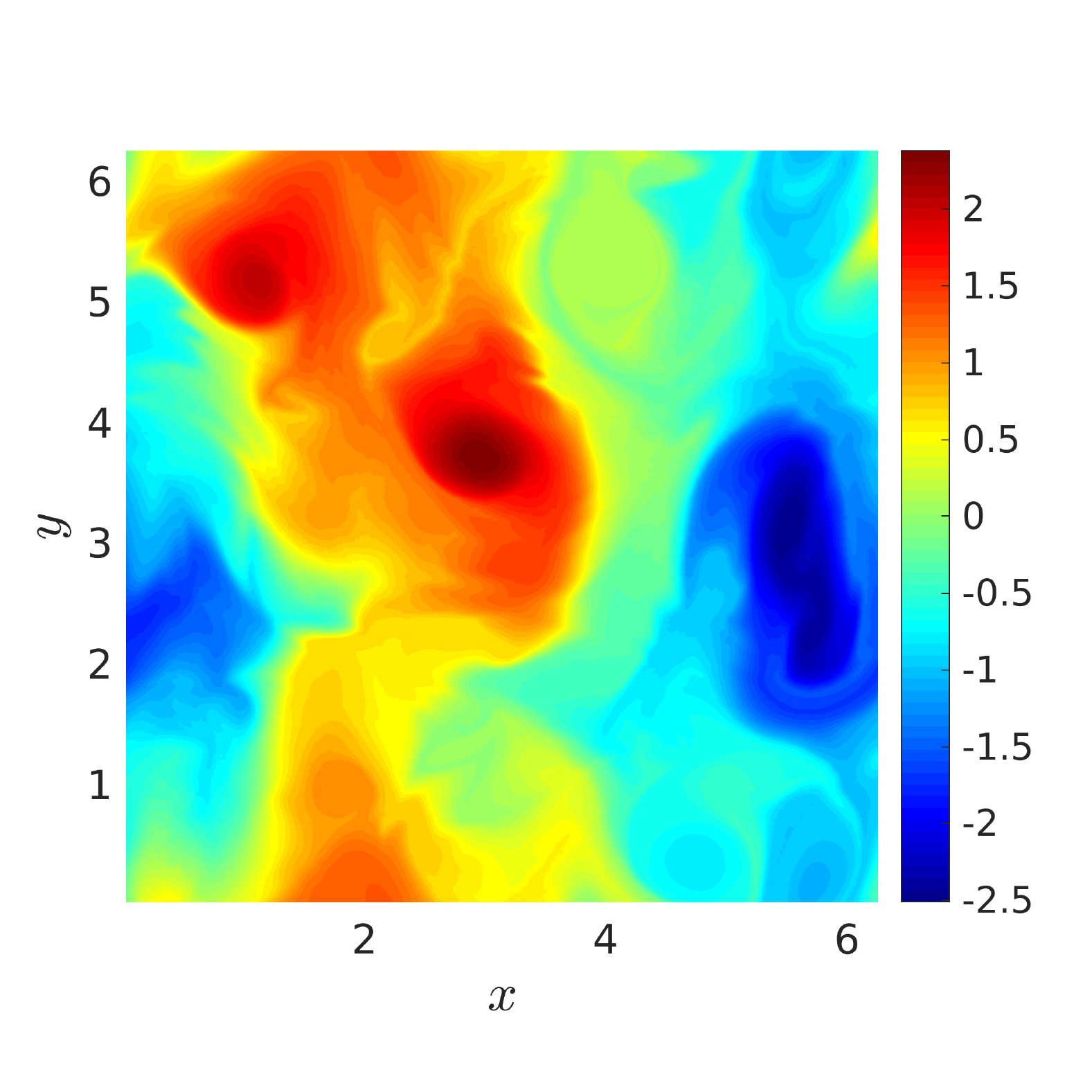} \hspace{0.1cm}
\includegraphics[width=0.325\textwidth]{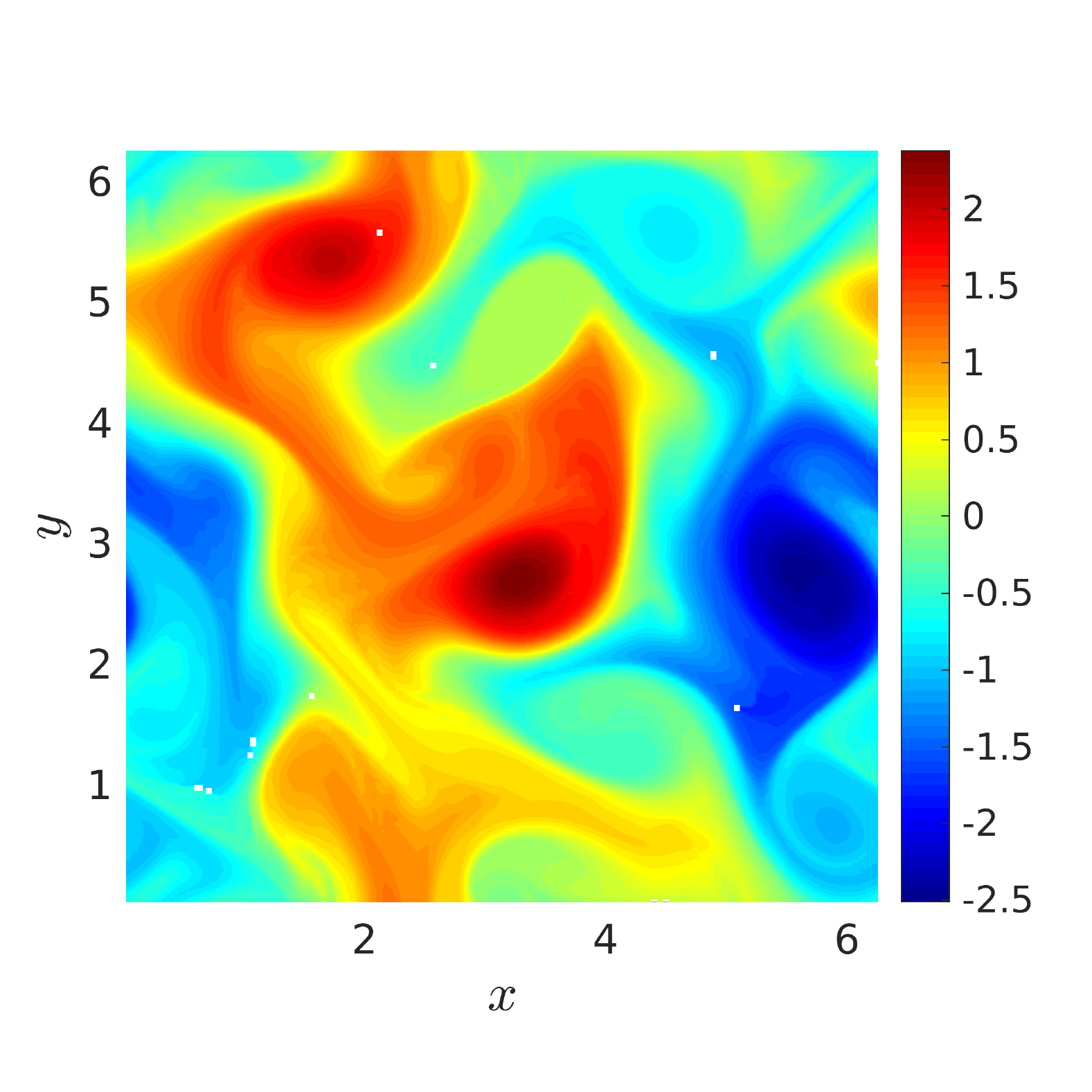}
\caption{Second left singular vector (left); second right singular vector (right).}
\label{fig:turbulence_eigenvector_orig}
\end{figure}

\begin{figure}[hbt]
\centering
\includegraphics[width=0.325\textwidth]{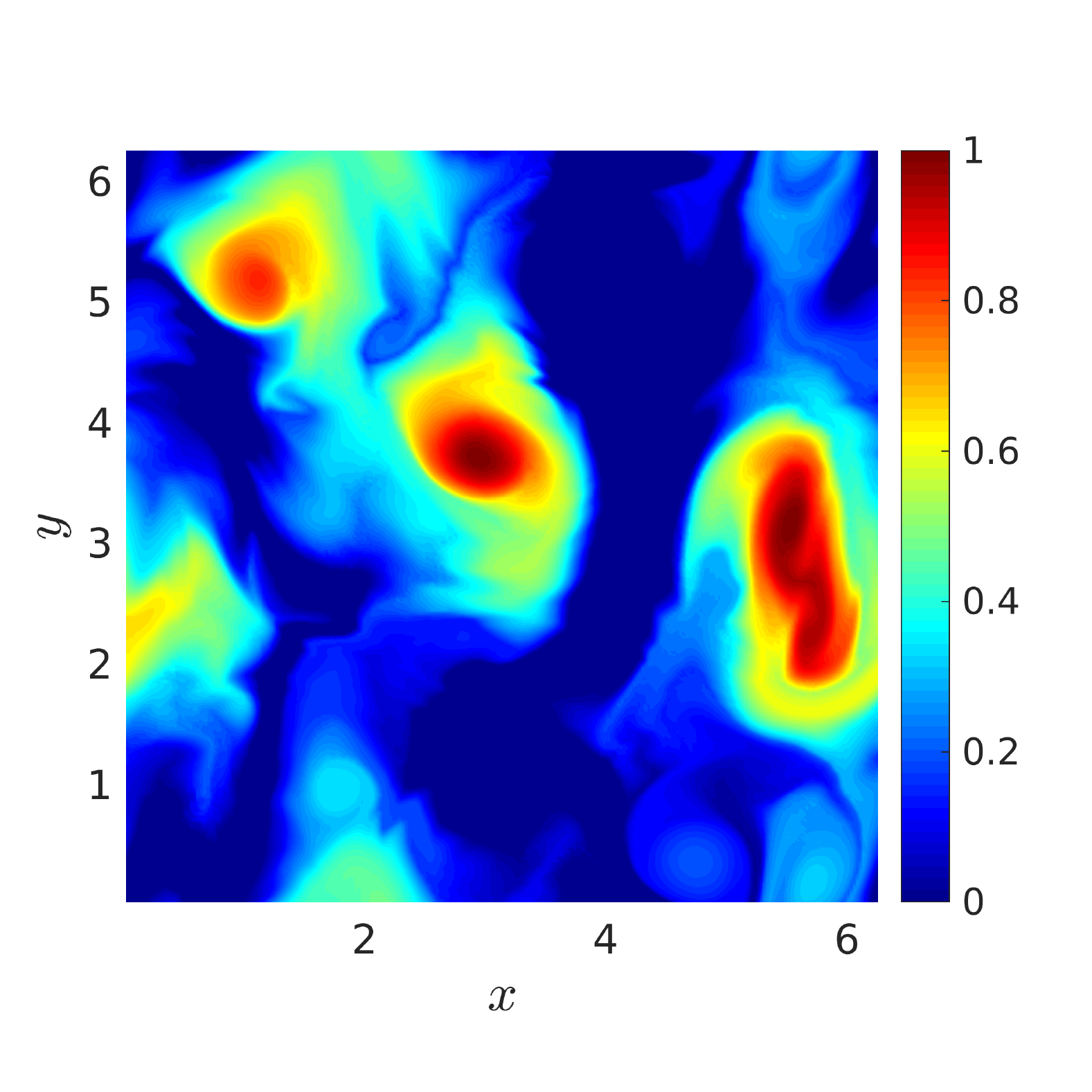} \hspace{0.1cm}
\includegraphics[width=0.325\textwidth]{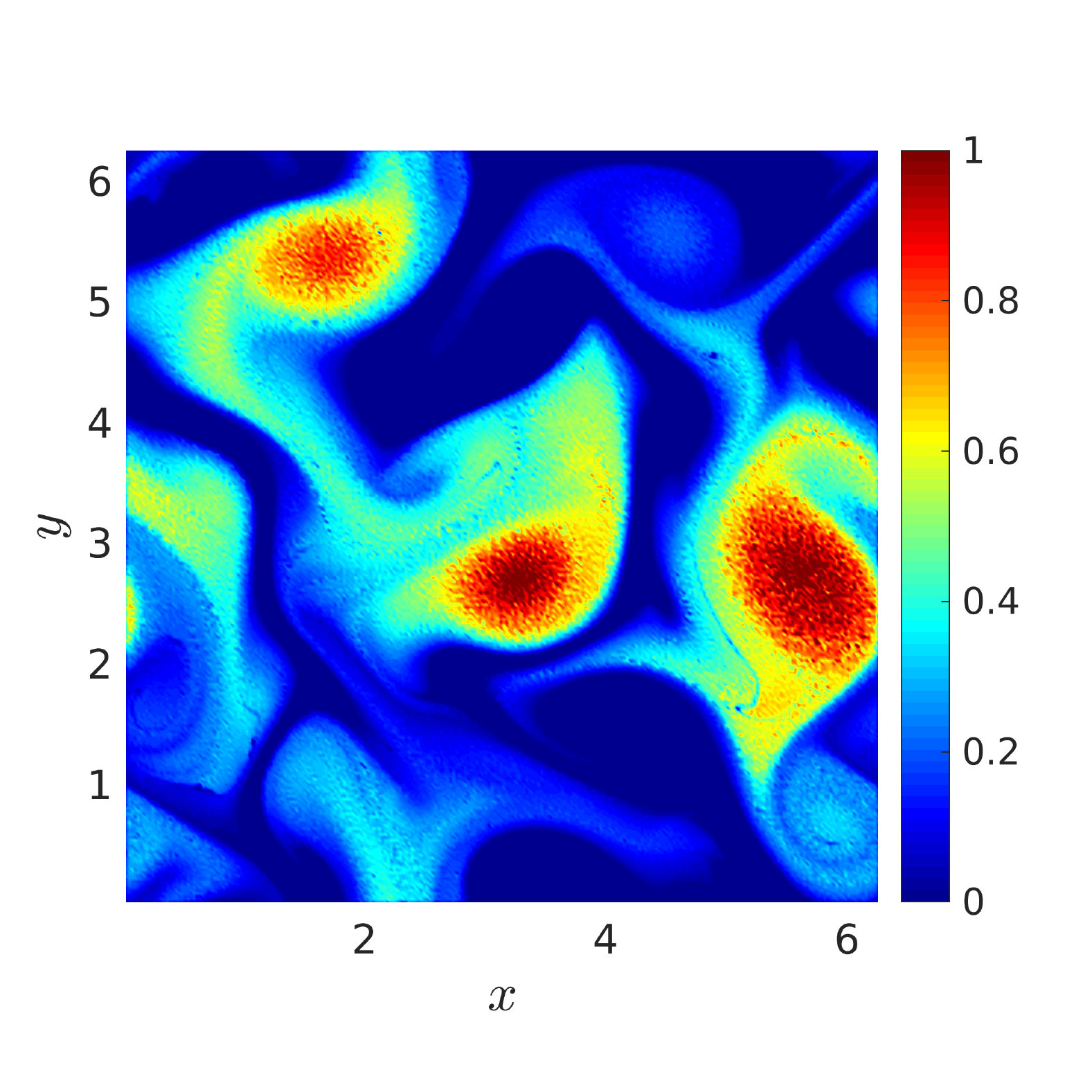}
\caption{Superposition vector $\mathfrak{s}$ (see Section \ref{sec:subpart} for the definition) of the two sparse vectors generated for the turbulence flow using two singular vectors ($r=2$) (left). Push forward of the left image $\mathfrak{s}$ by $L$ (right).}
\label{fig:turbulence_eigenvector_sparse}
\end{figure}

We consider a numerical solution to the Navier-Stokes equations in two dimensions with random-in-phase forcing, introduced in \cite{farazmand13}. Specifically, we consider the velocity field $u:\mathbb{T}^2\times \mathbb{R} \to \mathbb{R}^2$ which solves the equations
\begin{align}
\partial_t u + u \cdot \nabla u &= -\nabla p + \nu \Delta u + f \nonumber \\
\nabla \cdot u &= 0
\end{align}
where the toral domain $\mathbb{T}^2=[0, 2\pi] \times [0, 2 \pi]$; see \cite{farazmand13} for details.
We use the Ulam approximation $P$ of the transfer operator $\mathcal{P}$ constructed on a $256\times 256$ grid, for a flow time of 50 time units, exactly as in \cite{Hadjighasem_etal}.
Using the method of \cite{FSM10,F13}, we compute singular vectors of the matrix $L$ (in the notation of \cite{FSM10,F13}, which is a slight reweighting of the stochastic matrix $P$).
The second left and right singular vectors of $L$ are shown in Figure \ref{fig:turbulence_eigenvector_orig}.
The left and right images in Figure \ref{fig:turbulence_eigenvector_orig} may be compared directly with \cite[Figure 9(g) and Figure 11(a)]{Hadjighasem_etal}, respectively, where in the latter figures, only a binary plot of positive and negative values was displayed.
In Figure \ref{fig:turbulence_eigenvector_orig} there are three highlighted coherent regions:  two with extreme positive values (red), and one with extreme negative values (dark blue).
When applying Algorithm \ref{alg1} to the first two left singular vectors (the leading left singular vector is constant and the second left singular vector is shown in Figure \ref{fig:turbulence_eigenvector_orig} (left)), these three features are separated and highlighted in red;  see Figure \ref{fig:turbulence_eigenvector_sparse}.
The sparse basis vectors shown in Figure \ref{fig:turbulence_eigenvector_sparse} are considerably more informative than either (i) a simple partition based on the sign of the singular vectors or (ii) the result of a hard partitioning method such as k-means. The sparsity quantifiers in the case of $r=2$ are 31.6\% (subspace error), 38.4\% (absolute sparsity) and 54.6\% (relative sparsity).

\begin{figure}[hbt]
\centering
\includegraphics[scale=0.5]{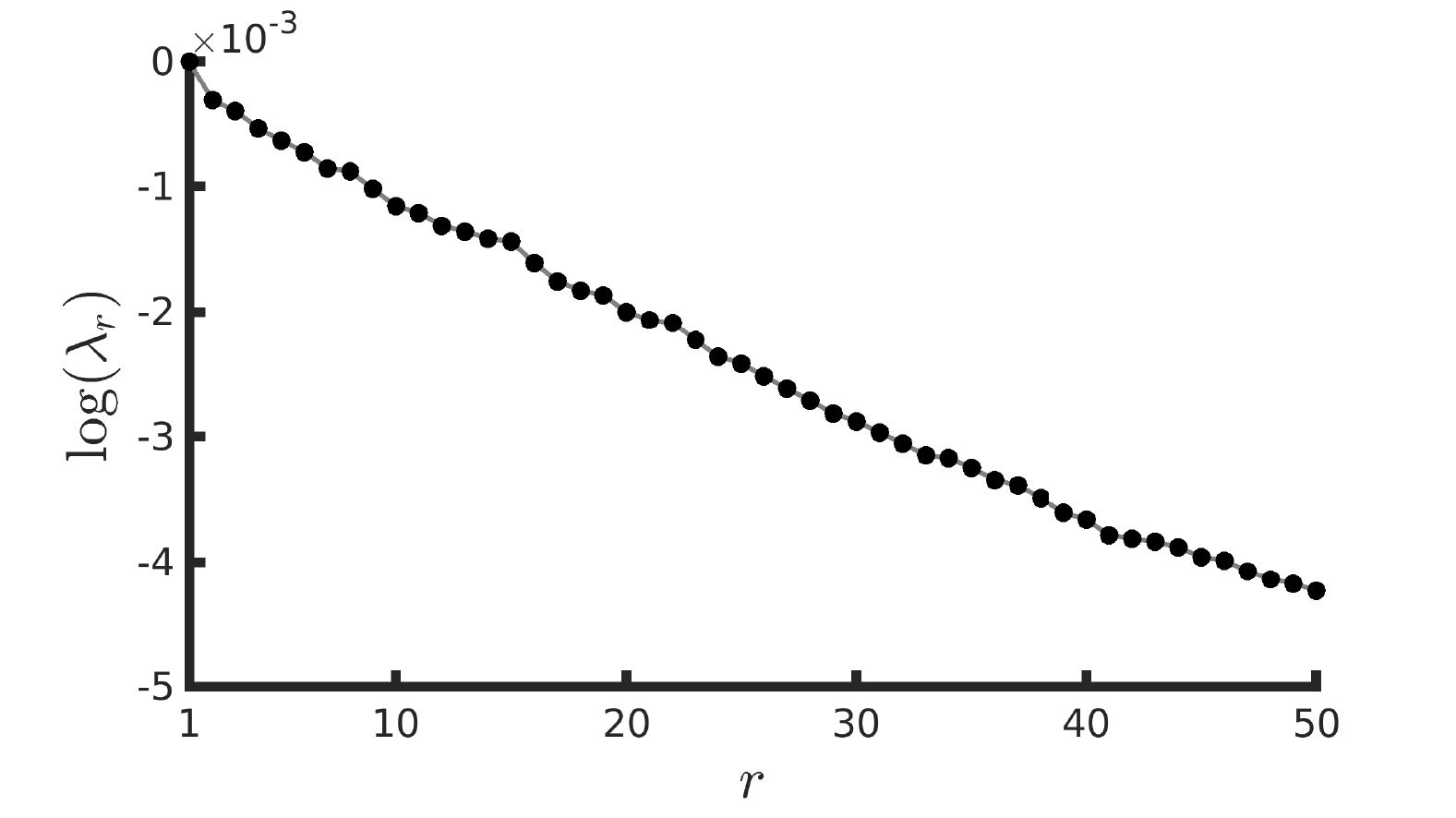}
\includegraphics[scale=0.5]{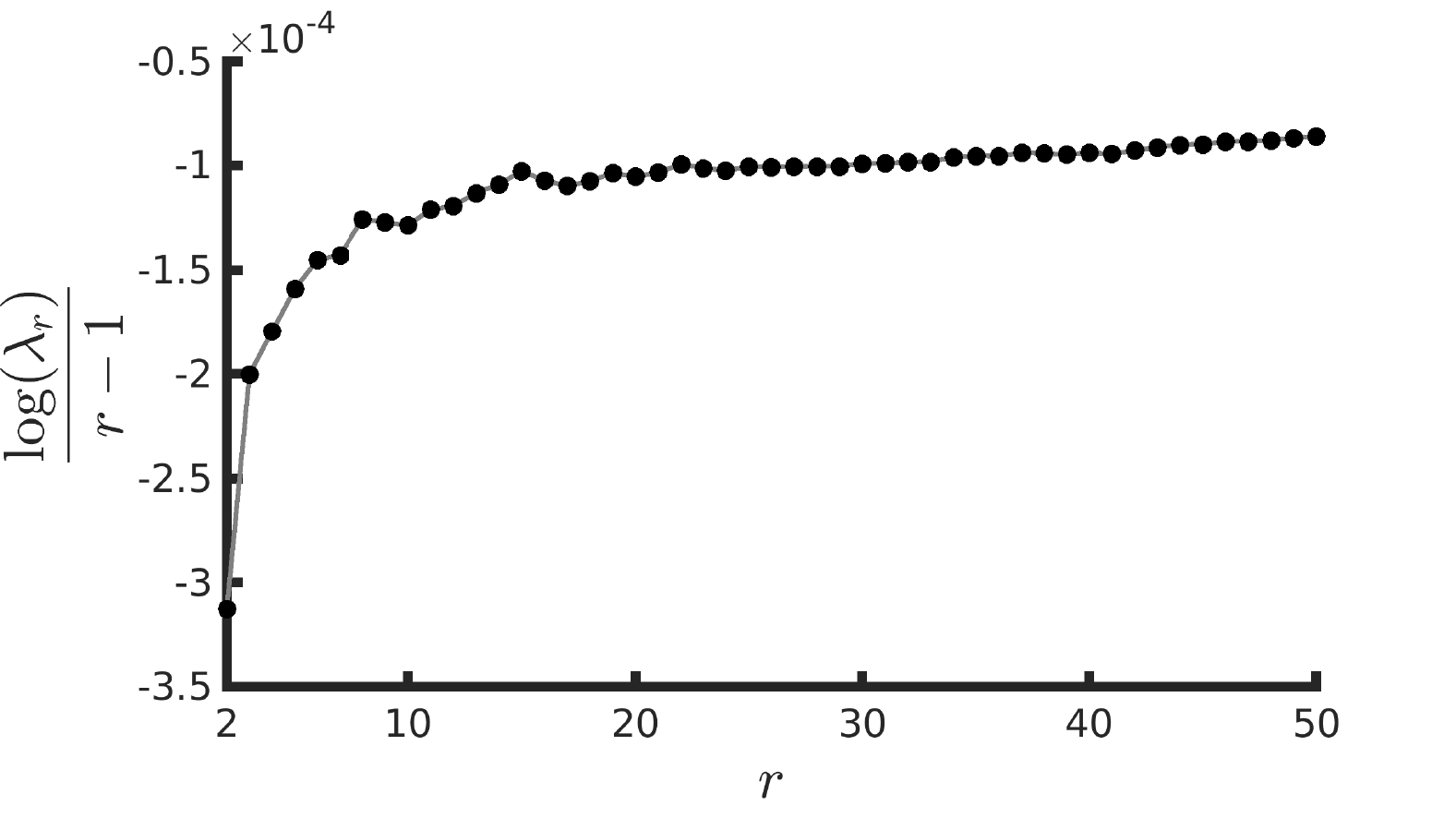}
\caption{Plot of $\log\lambda_r$ vs $r$ (left) and $\log\lambda_{r}/(r-1)$ vs $r$ (right) for the turbulence flow.}
\label{fig:turbulence_spectrum}
\end{figure}

Turbulence is often cited as a canonical example of dynamics operating over a broad range of spatial scales.
This fact is backed up by the lack of a clear eigengap for the transfer operator singular value spectrum;  see Figure \ref{fig:turbulence_spectrum} (left).
The Weyl rescaling in Figure \ref{fig:turbulence_spectrum} (right) shows some minor drops, which may provide weak eigengap information.
In situations where there is no clear eigengap, we propose to use information from the sparse vectors output by Algorithm \ref{alg1} to decide how many sparse vectors contain ``reliable'' coherent features.
The upper envelope (dotted black)
of Figure \ref{fig:turbulence_cumulative} shows $\Min(S^{(r)})$ vs. $r$.
There are several ``drops'' in this envelope, with the three largest drops at $r=5,30$, and 36.
In this turbulence example, we illustrate a refinement of this choice of $r$, simultaneously taking into account the ``reliability'' of the sparse vectors.




\subsubsection{A sparse vector heuristic for simultaneously selecting the number of input vectors $r$ and number of features $k$}
\label{sec:heuristic}
\begin{figure}[hbt]
\centering
\includegraphics[width=\textwidth]{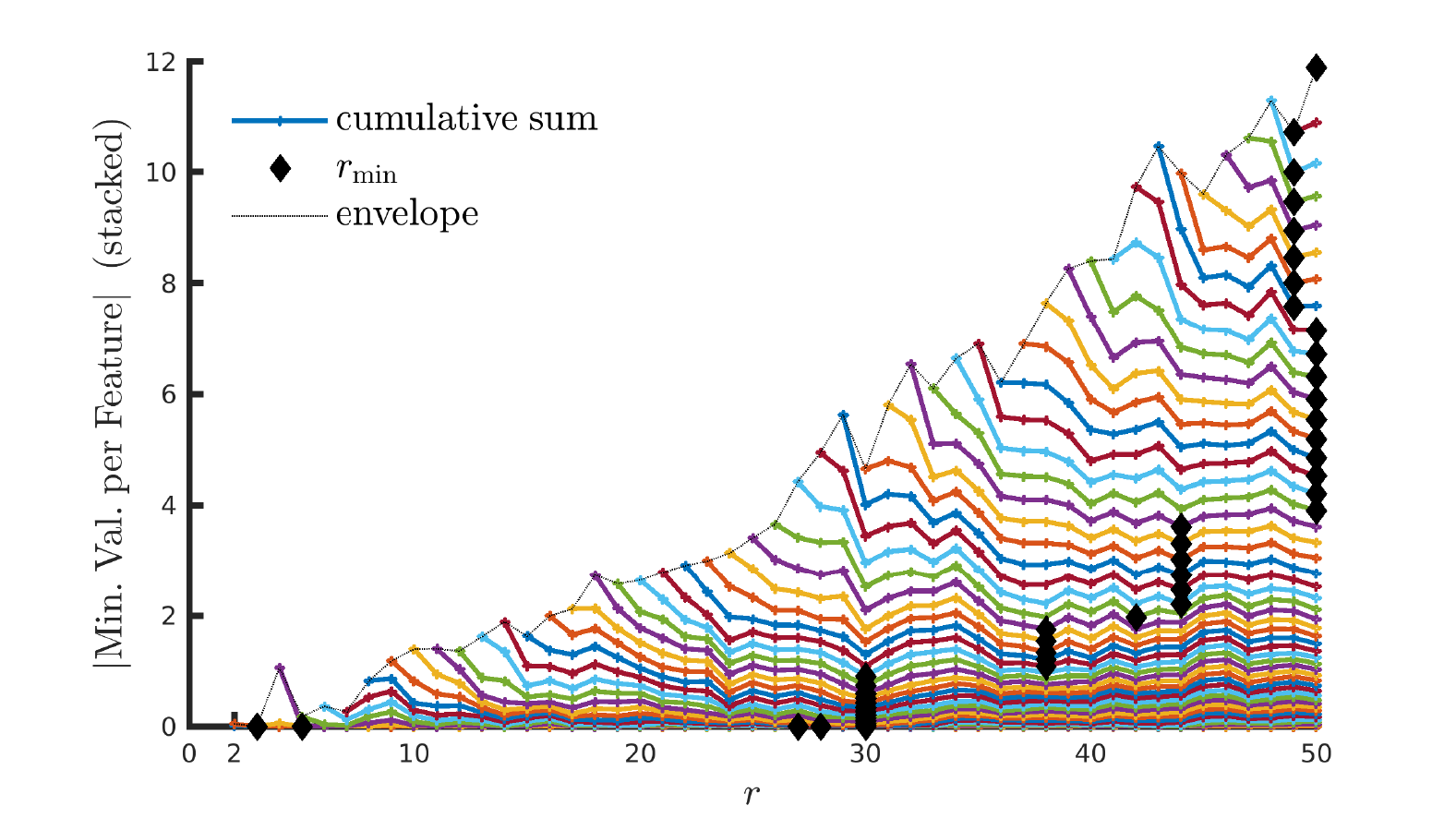}
\caption{``Stacked'' plot of sparse vector minimum value curves.} \label{fig:turbulence_cumulative}
\end{figure}

\begin{figure}[hbt]
\centering
\includegraphics[width=.6\textwidth]{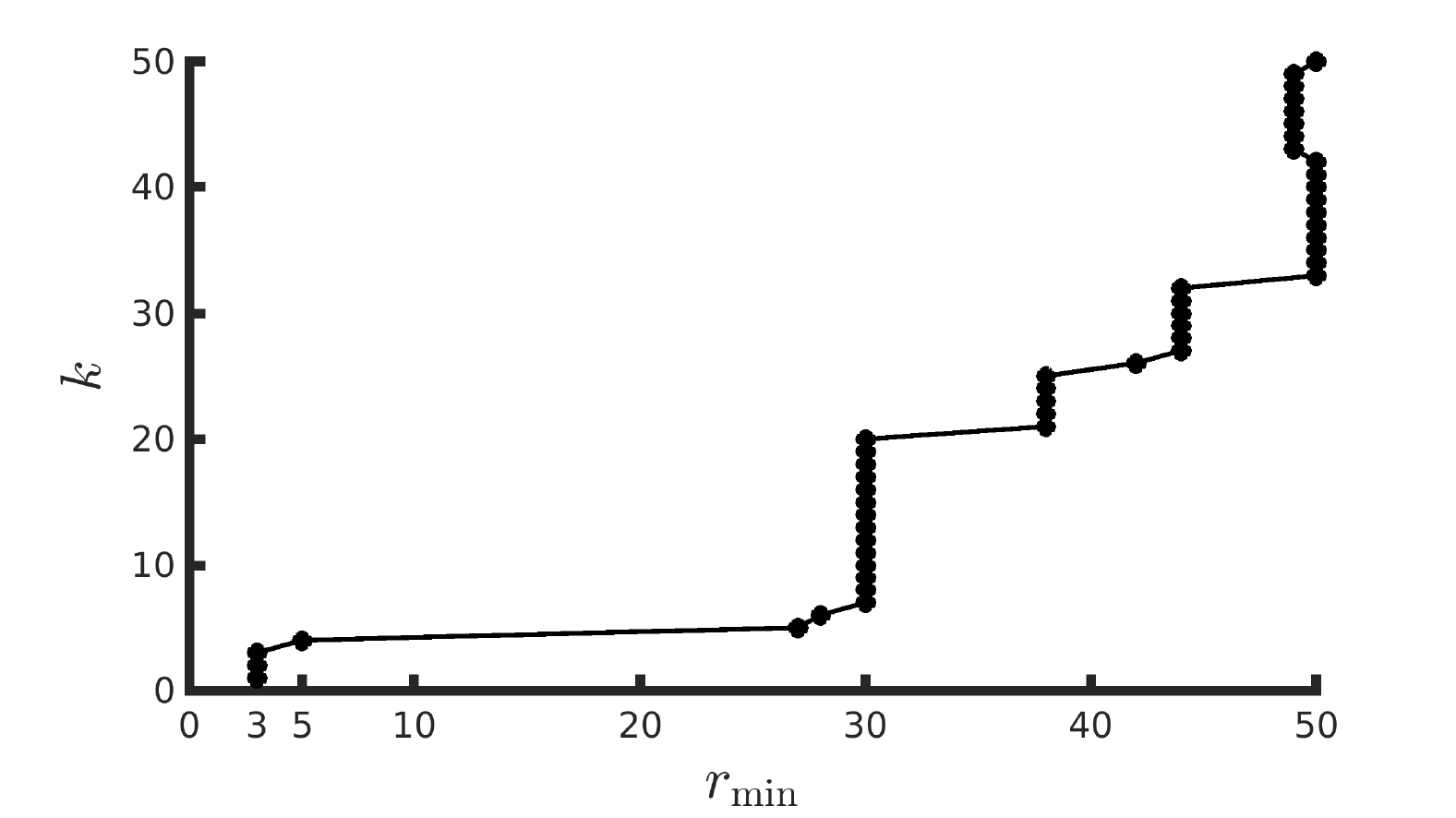}
\caption{Turbulence example: $r_{\mathrm{min}}(k)$ vs number of coherent features $k$.} \label{fig:turbulence_rmin}
\end{figure}

For some maximum number of columns $r_{\max}$ and each $r=2,\ldots,r_{\max}$, we apply Algorithm \ref{alg1} to the first $r$ columns of the eigenvector matrix $V$, and denote the resulting arrays $S^{(r)}$.
For each $r$ and $k$, let $\Min(S^{(r)},k):=\sum_{j=1}^k -\min_{1\le i \le p} S^{(r)}_{ij}$ be the sum of the minimum values of the first $k$ sparse vectors.
In Figure \ref{fig:turbulence_cumulative}, for each $k \le r$, we plot $\Min(S^{(r)},k)$ vs. $r$ as a coloured line, for each $2\le r\le r_{max}$.
The $k^{\rm th}$ coloured line shows the ``performance'' (as determined by the minimum value heuristic) of Algorithm \ref{alg1} at separating $k$ features from various $r$ input vectors.
The best performance for $k$ features is achieved when a single coloured line attains a minimum of $\Min(S^{(r)},k)$ across the various values of $r$.
For a given $k$, these minimal $r$ are denoted
$r_{\min}(k):=\argmin_{k \le r \le r_{\max}} \Min(S^{(r)},k)$ (when multiple $r$ values give the same $\Min(S^{(r)},k)$ value, we take $r_{\min}(k)$ as the smallest of these).
These ``optimal'' $(k,r)$ combinations are shown as black diamonds on Figure \ref{fig:turbulence_cumulative}.
Figure \ref{fig:turbulence_rmin} shows a plot of $r_{\min}(k)$ vs.\ $k$.

\begin{figure}[hbt]
\centering
\includegraphics[width=0.325\textwidth]{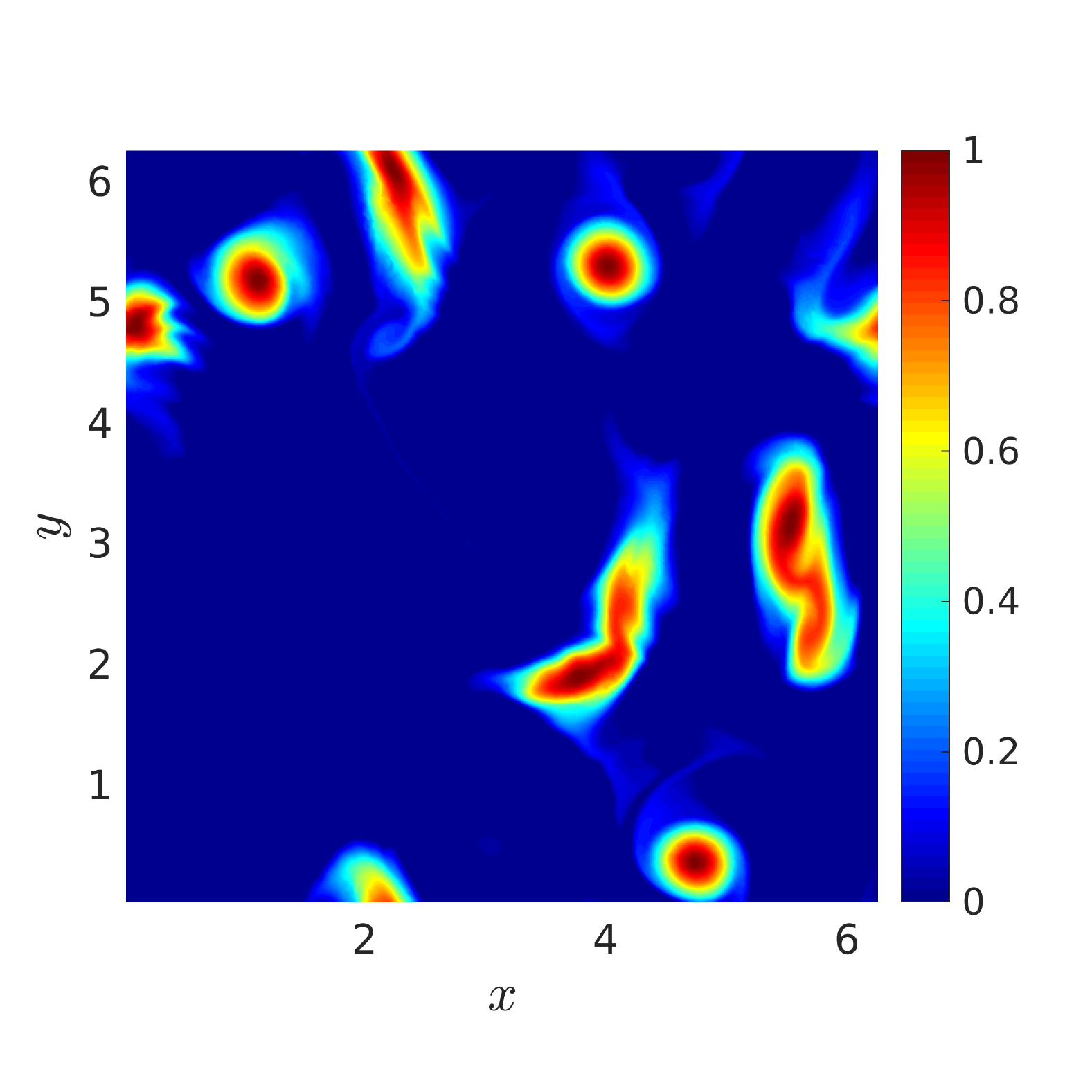} \hspace{0.1cm}
\includegraphics[width=0.325\textwidth]{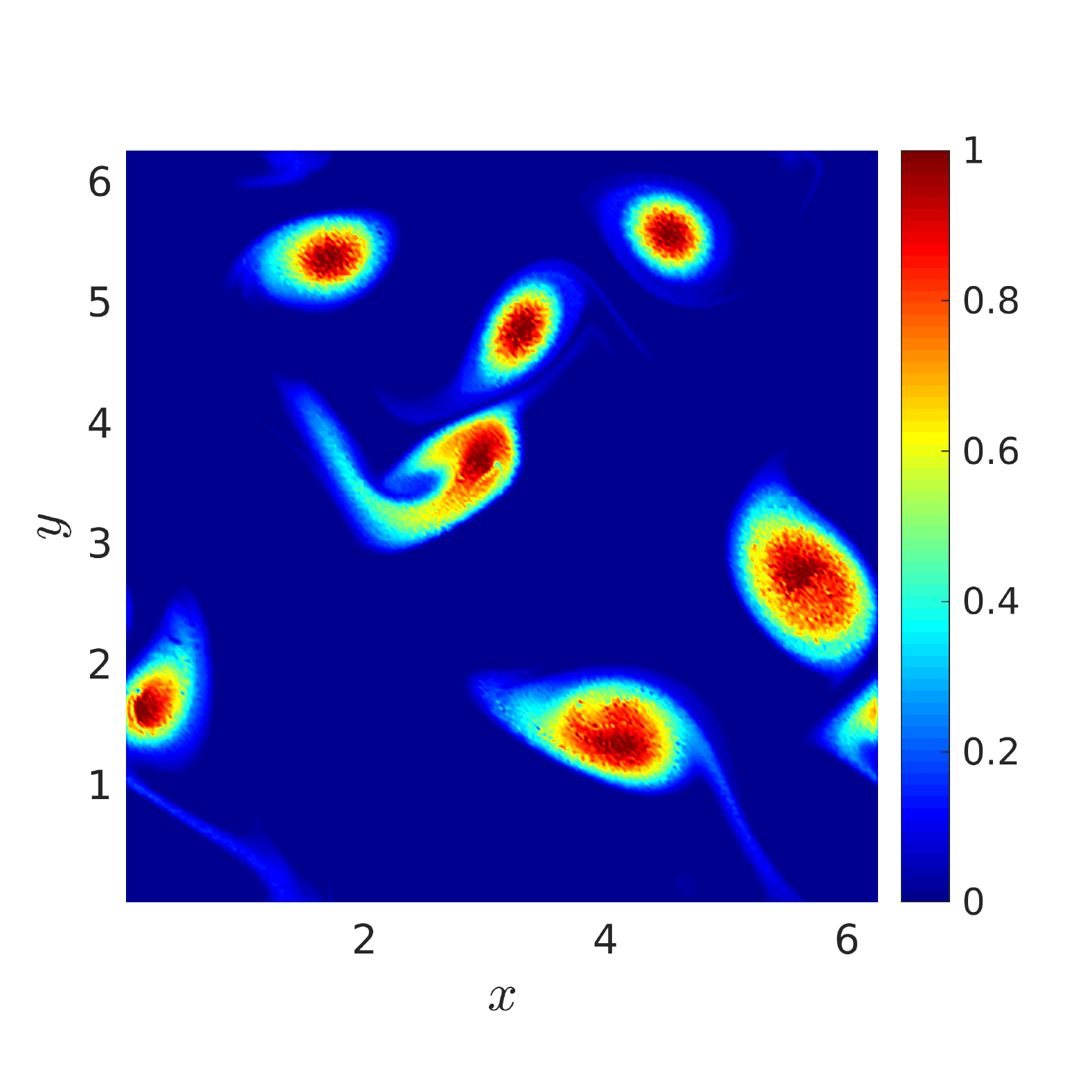}
\caption{Superposition vector $\mathfrak{s}$ of the 7 most reliable features at initial time (left) according to the minimum value reliability criterion ($r=30$) and its forward-time image (right).}
\label{fig:turbulence_sparse_r30_k7_sparse}
\end{figure}

\begin{figure}[hbt]
\centering
\includegraphics[width=0.325\textwidth]{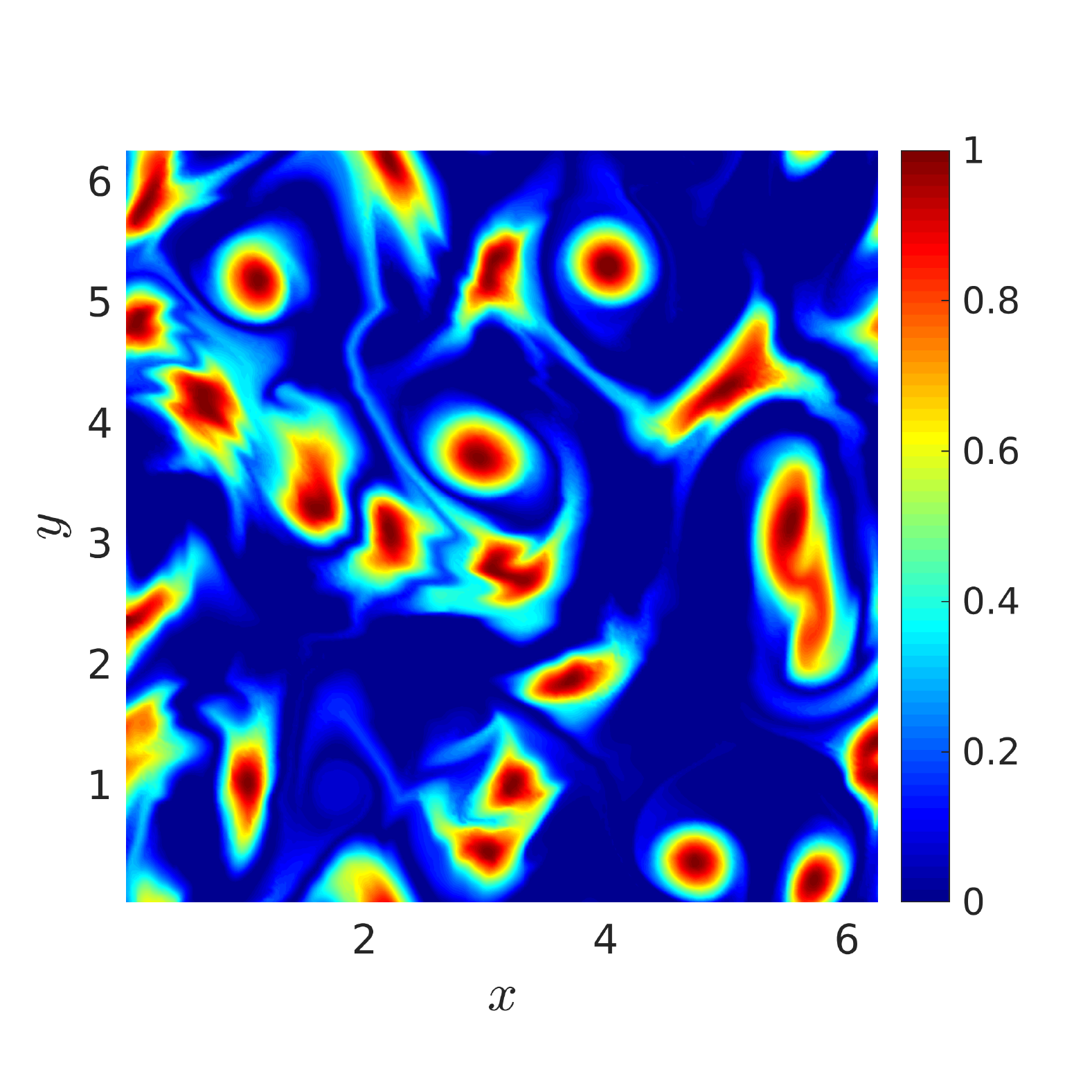} \hspace{0.1cm}
\includegraphics[width=0.325\textwidth]{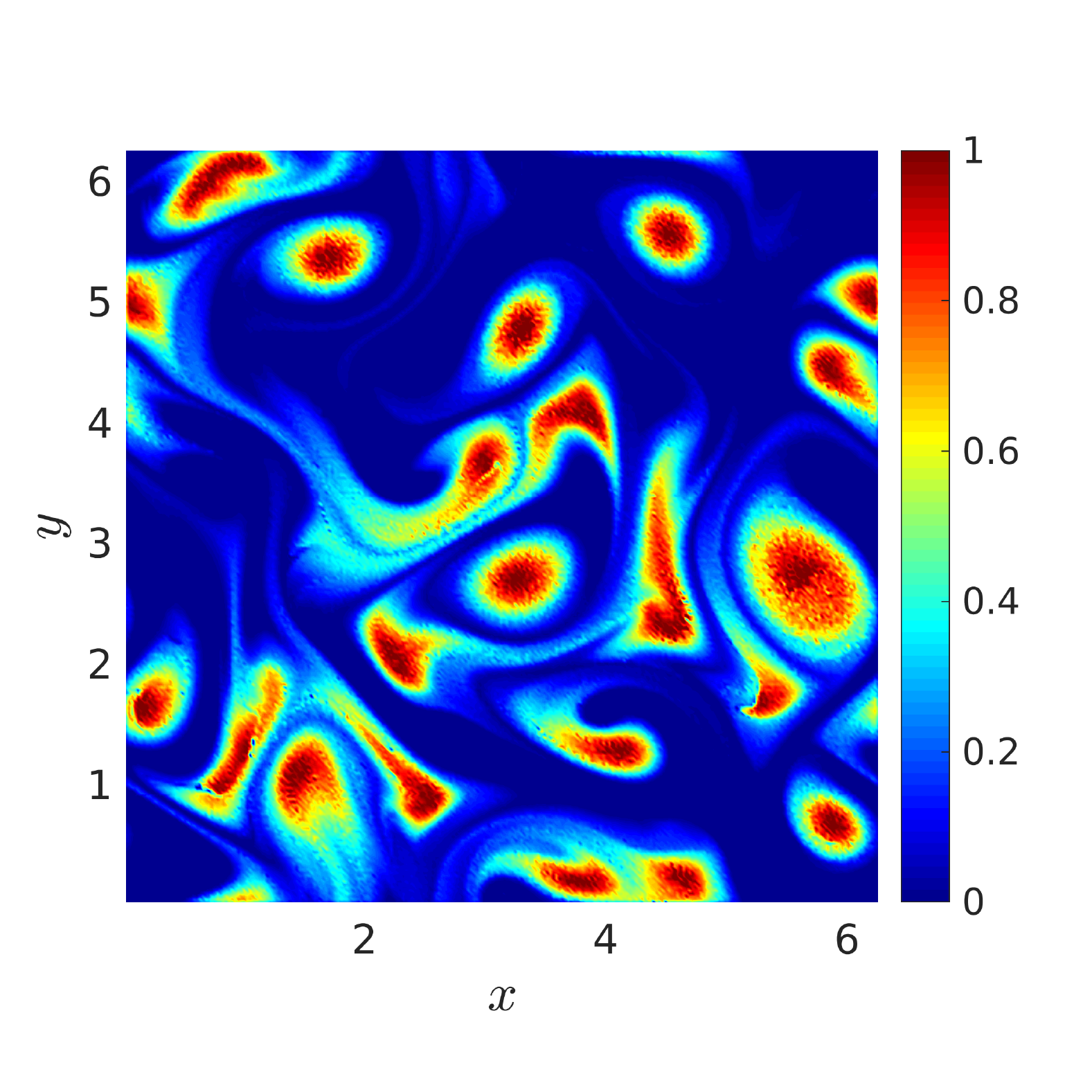}
\caption{Superposition vector $\mathfrak{s}$ of the 21 most reliable features at initial time (left) according to the minimum value reliability criterion ($r=38$) and its forward-time image (right).}
\label{fig:turbulence_sparse_r38_k21_sparse}
\end{figure}

Notice that there are $r_{\min}(k)$ takes on only certain $r$ values.
These common $r$ values (e.g.\ $r=30,38,\ldots$ in Figure \ref{fig:turbulence_cumulative}) could be interpreted as representing \emph{natural spatial scales} of the dynamics because one expects the eigenvectors (or singular vectors) input to Algorithm \ref{alg1} to become more oscillatory as $r$ is increased, highlighting increasingly smaller spatial features.
One may select any of these optimal $(k,r)$ combinations. To be conservative, in our experiments we select $k$ at the lower end of each ``vertical strip'' of black diamonds.
If we increase $k$ we often simply obtain more coherent sets.
For this example, we select $(k,r)=(7,30)$ and $(21,38)$ for further analysis;  the corresponding sparse superposition vectors and their forward-time images are shown in Figures \ref{fig:turbulence_sparse_r30_k7_sparse} and \ref{fig:turbulence_sparse_r38_k21_sparse}.
From these two figures one already sees 7 (resp.\ 21) clearly highlighted coherent sets in red that remain coherent after the flow time of 40 time units.
The span of the collection of sparse vectors in both the case of $r=30$ and that of $r=38$ produce reasonable approximations of the span of the corresponding leading eigenvectors with a subspace error of 10.4\% (resp. 9.8\%). The sparse vectors have an absolute sparsity of 7.7\% (resp. 7.1\%) and a relative sparsity of 26.8\% (resp. 25.2\%).

\subsubsection{Extraction of a subpartition from sparse vector output} \label{sec:threshcheeger}
One could apply a single global threshold to our selected sparse vectors according to
Algorithms \ref{thresh1} and \ref{thresh2}.
We found that because Figure \ref{fig:turbulence_sparse_r30_k7_sparse} is already a collection of vectors with mostly disjoint supports, both of these algorithms produced a very low threshold with little change to the sparse vectors.
On the other hand, the sparse vectors that are combined in Figure \ref{fig:turbulence_sparse_r38_k21_sparse} have ``overlaps'' of large values in small parts of phase space, causing Algorithms \ref{thresh1} and \ref{thresh2} to produce very high thresholds close to 1.
In place of these thresholding algorithms, one could easily apply a global threshold manually by inspecting Figure \ref{fig:turbulence_sparse_r30_k7_sparse} or \ref{fig:turbulence_sparse_r38_k21_sparse};  for example a threshold of around 0.7--0.75 would produce coherent sets with small boundary lengths relative to enclosed area at the initial and final times.


To end this section we apply thresholding according to the Cheeger ratio.
Such a thresholding has been applied to specific single vectors \cite{F15,FJ15,FK17,Hadjighasem_etal}, but because Algorithm \ref{alg1} has automatically separated \emph{all} features, we may apply this thresholding to either the superposition vector $\mathfrak{s}$ or to each individual vector sparse vector $s_j$, $j=1,\ldots,k$.
Given a vector (corresponding to a discrete representation of a real-valued function on a manifold $M$), we select a threshold by choosing the level set that minimises a quantity we call the \emph{scale-invariant Cheeger ratio}. Let $M$ be a $d$-dimensional manifold and let $T:M\to M$ be a volume-preserving transformation representing our dynamics. For simplicity of presentation, we assume that one application of $T$ represents the nonlinear dynamics over the full finite-time interval and that we compute the coherence via the Cheeger ratio at only the initial and final time;  see \cite{F15,FK17} for more detail and generalisations.
Let $\Gamma\subset M$ be a $d-1$ dimensional submanifold disconnecting $M$ into two pieces $M_1, M_2$.
The scale-invariant Cheeger ratio is
\begin{equation}
\label{eq:cheeger}
h^D(\Gamma):=\frac{\ell_{d-1}(\Gamma)+\ell_{d-1}(T(\Gamma))}{2(\min\{\ell_d(M_1),\ell_d(M_2)\})^{\frac{d-1}{d}}},
\end{equation}
where $\ell_{d-1}$ and $\ell_d$ are $(d-1)$- and $d$-dimensional volume measure, respectively.
In this example we treat each sparse vector $s_j$ separately, and each $\Gamma$ will be selected as a level set of $s_j$, $j=1,\ldots,k$.
In particular we choose the level set value so that (\ref{eq:cheeger}) is minimised.
In terms of numerics, having selected a $(k,r)$ combination, for each $1\le j\le k$, we linearly interpolate $s_{j}$ on a fine regular grid on $M$ to obtain a function $u_j: M\to \mathbb{R}$.
For threshold values $\tau$ ranging from $\min_x u_j(x)$ to $\max_x u_j(x)$
we calculate $h^D(\Gamma_\tau)$, where $\Gamma_\tau$ is the level set $\{x\in M:u_j(x)=\tau\}$, and select the $\tau$ yielding the minimum value.
We use MATLAB's \verb"contourc" command to calculate level sets; \verb"isosurface" is the corresponding command for three dimensions. This takes approximately two minutes for 21 interpolated sparse vectors in $\mathbb{R}^{65536}$  on an Intel Core i7-6700 processor with 16GB RAM.
Note that in contrast to Algorithm \ref{thresh1} or \ref{thresh2}, we calculate a \emph{separate threshold} for each vector before applying a hard thresholding.

\begin{figure}[hbt]
\centering
\includegraphics[width=0.75\textwidth]{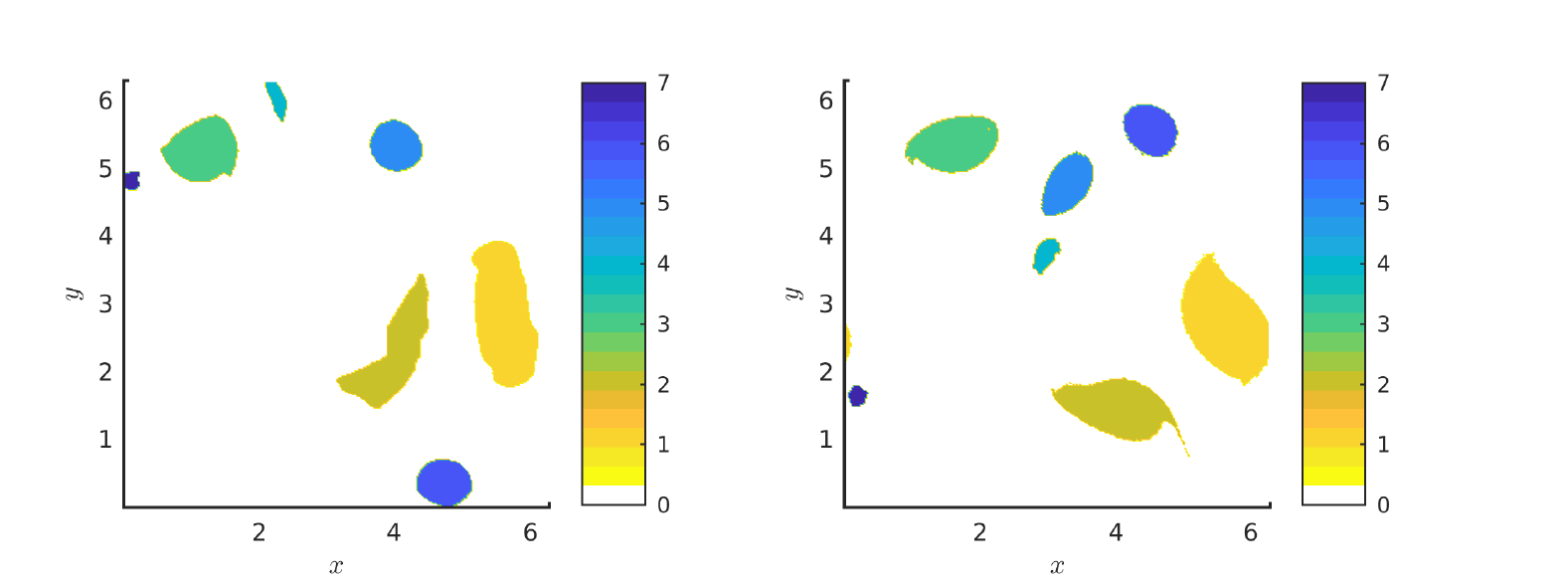}
\caption{Partition $(k,r)=(7,30)$ of phase space showing 7 coherent sets at initial- (left) and forward- time (right) obtained by thresholding.}
\label{fig:turbulence_sparse_r30_k7_thresh}
\end{figure}

We apply this thresholding procedure to the $(k,r)$ combination (7,30) whose superposition of sparse vectors is shown in Figure \ref{fig:turbulence_sparse_r30_k7_sparse}); the result is the partition in Figure \ref{fig:turbulence_sparse_r30_k7_thresh}.
All features obtained are highly coherent, though two features are perhaps overly small.
The scale-invariant dynamic Cheeger ratio for these two small sets has a very flat minimum over a large range of thresholds and so these two sets could be further enlarged without changing appreciably the Cheeger ratio (and therefore the finite-time coherence).

\begin{figure}[hbt]
\centering
\includegraphics[width=0.75\textwidth]{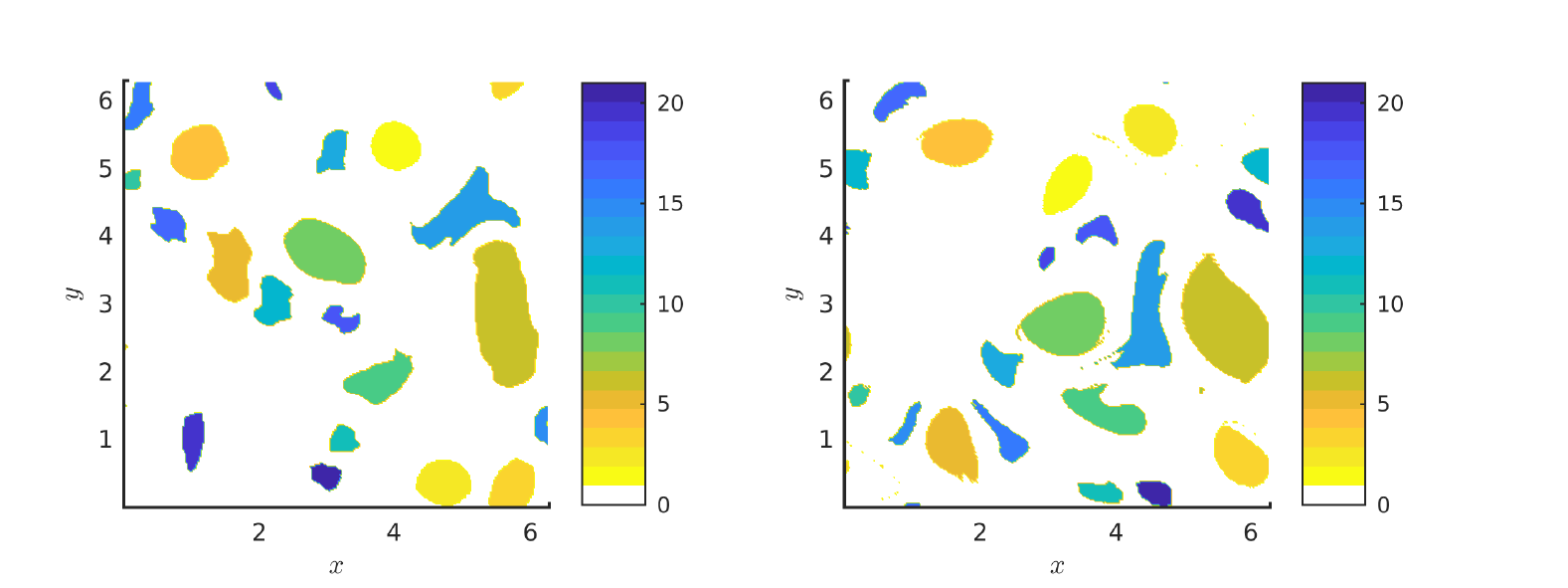}
\caption{Partition $(k,r)=(21,38)$ of phase space showing 21 coherent sets at initial- (left) and forward- time (right) obtained by thresholding.}
\label{fig:turbulence_sparse_r38_k21_thresh}
\end{figure}

We also apply this thresholding procedure to the $(k,r)$ combination (21,38),
whose superposition vector $\mathfrak{s}$ is shown in Figure \ref{fig:turbulence_sparse_r38_k21_sparse}.
The resulting thresholded sets are displayed in Figure \ref{fig:turbulence_sparse_r38_k21_thresh}, coloured according to their reliability (a reliability index of 1 is the most reliable (coherent)).
We retain all the coherent sets identified in Figure \ref{fig:turbulence_sparse_r30_k7_thresh}, and gain another 15 coherent sets.
In this case, while some of the additional coherent sets are relatively small, many of these sets cannot be significantly increased without introducing filamentation;  that is, they do not have flat minima for their Cheeger ratios.


\subsubsection{Extraction of a partition with k-means}
\begin{figure}
\centering
\includegraphics[width=0.75\textwidth]{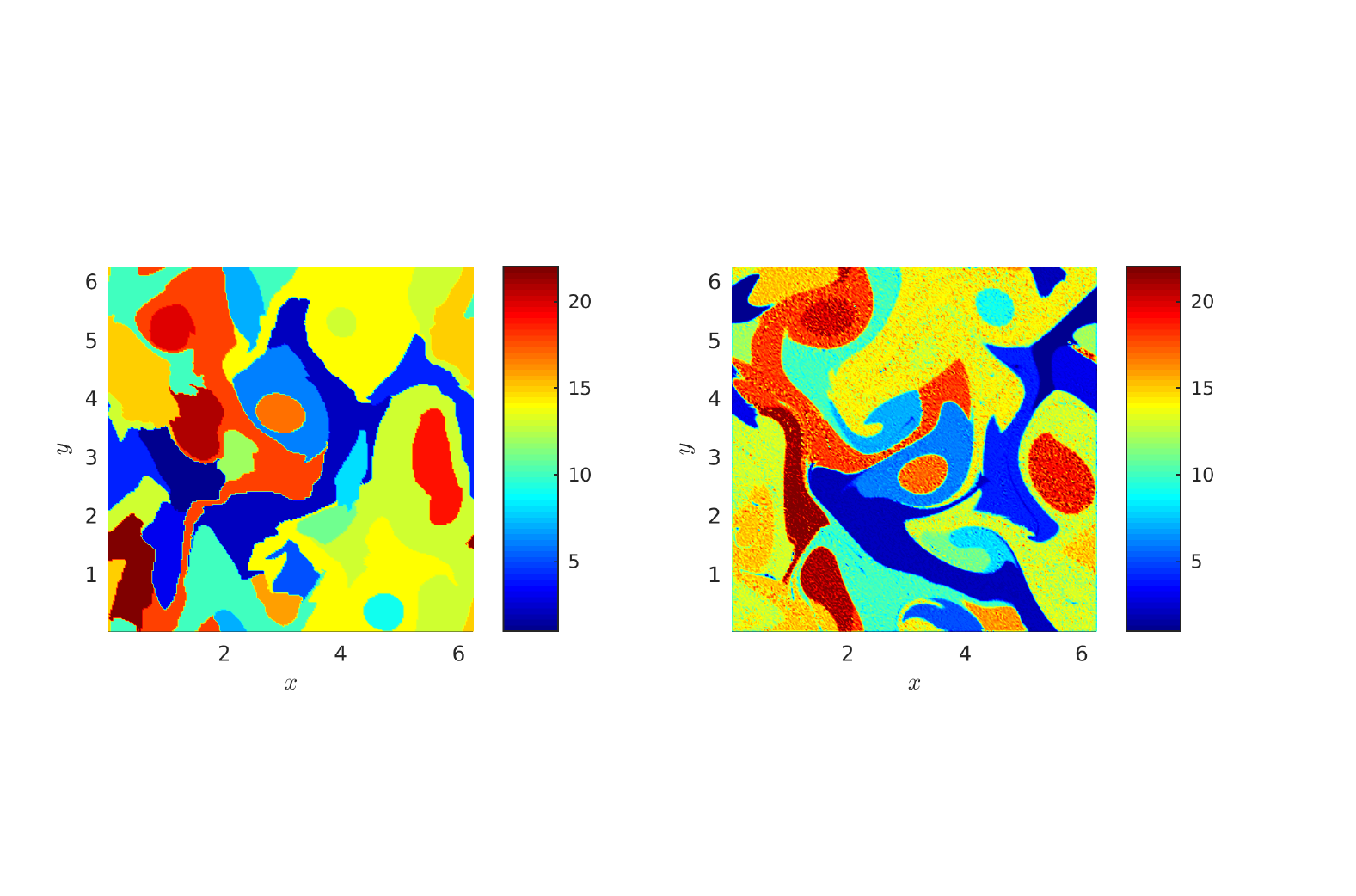}
\caption{Initial and forward-time plot of the k-means partition of state space into 22 clusters ($r=38$).}
\label{fig:turbulence_kmeans_r38_clust22}
\end{figure}

The use of sparse vectors via Algorithm \ref{alg1} clearly outperforms k-means clustering of the singular vectors.
As an example, we apply k-means to the $(k,r)$ combination (21,38), namely using 38 eigenvectors and requesting $k=21$ clusters; see Figure \ref{fig:turbulence_kmeans_r38_clust22}.
We add one extra cluster to give k-means a chance to find the complement of 21 coherent sets as a 22nd coherent set
The k-means clustering correctly detects some coherent sets, notably many of those in Figure \ref{fig:turbulence_sparse_r30_k7_thresh}, but also produces many more highly filamented sets, and misses several very coherent features visible in Figure \ref{fig:turbulence_sparse_r38_k21_thresh}.
To a large extent, k-means fails in difficult examples such as turbulence because k-means \emph{imposes} a partition of the phase space into coherent sets, when a partition is not reasonable.
The superposition of the sparse vectors in Figure \ref{fig:turbulence_sparse_r38_k21_sparse} highlights the coherent sets in a much more precise way than k-means can, and we can further refine Figure \ref{fig:turbulence_sparse_r38_k21_sparse} by thresholding, as in Figure \ref{fig:turbulence_sparse_r38_k21_thresh}.
Increasing $k$ causes k-means to add more filamented features, while decreasing $k$ causes k-means to miss some of the coherent features in Figures \ref{fig:turbulence_sparse_r38_k21_sparse} and \ref{fig:turbulence_sparse_r38_k21_thresh}.
We also experimented with the use of the k-means sum of squared errors (the objective that k-means attempts to minimise) to determine the number of clusters, but we failed to find any ``kinks'' to identify any values for $k$, as described for example in \cite{sugar99}.


\subsection{North Atlantic surface currents}\label{sec:northatlantic}
As our final example, we consider a real-world dataset, namely geostrophic AVISO
\cite{aviso} velocity fields.
We consider the 90-day period 15 January 2015 to 15 April 2015;  the velocity fields are updated every 2 hours.
We consider a synthetic initial $291\times 224$ Mercator grid of tracers on the domain $\mathcal{M}=[-65^\circ,-35^\circ]\times[30^\circ,50^\circ]$.
After removing trajectories initialised on land, we numerically integrate the remaining 56205 trajectories, linearly interpolating the velocity field between the AVISO-provided velocity fields spaced one day apart.
Using trajectory positions\footnote{In particular, derivatives of the velocity field or flow map need not be estimated.} every 3 days, we compute the dynamic Laplacian using the ``adaptive transfer operator approach'' of \cite{FJ18} on alpha complexes\footnote{With the lowest alpha such that the alpha complexes are connected.} with zero Neumann boundary conditions;  see \cite{FJ18} for details on implementation.

\begin{figure}[hbt]
\centering
\includegraphics[width=0.6\textwidth]{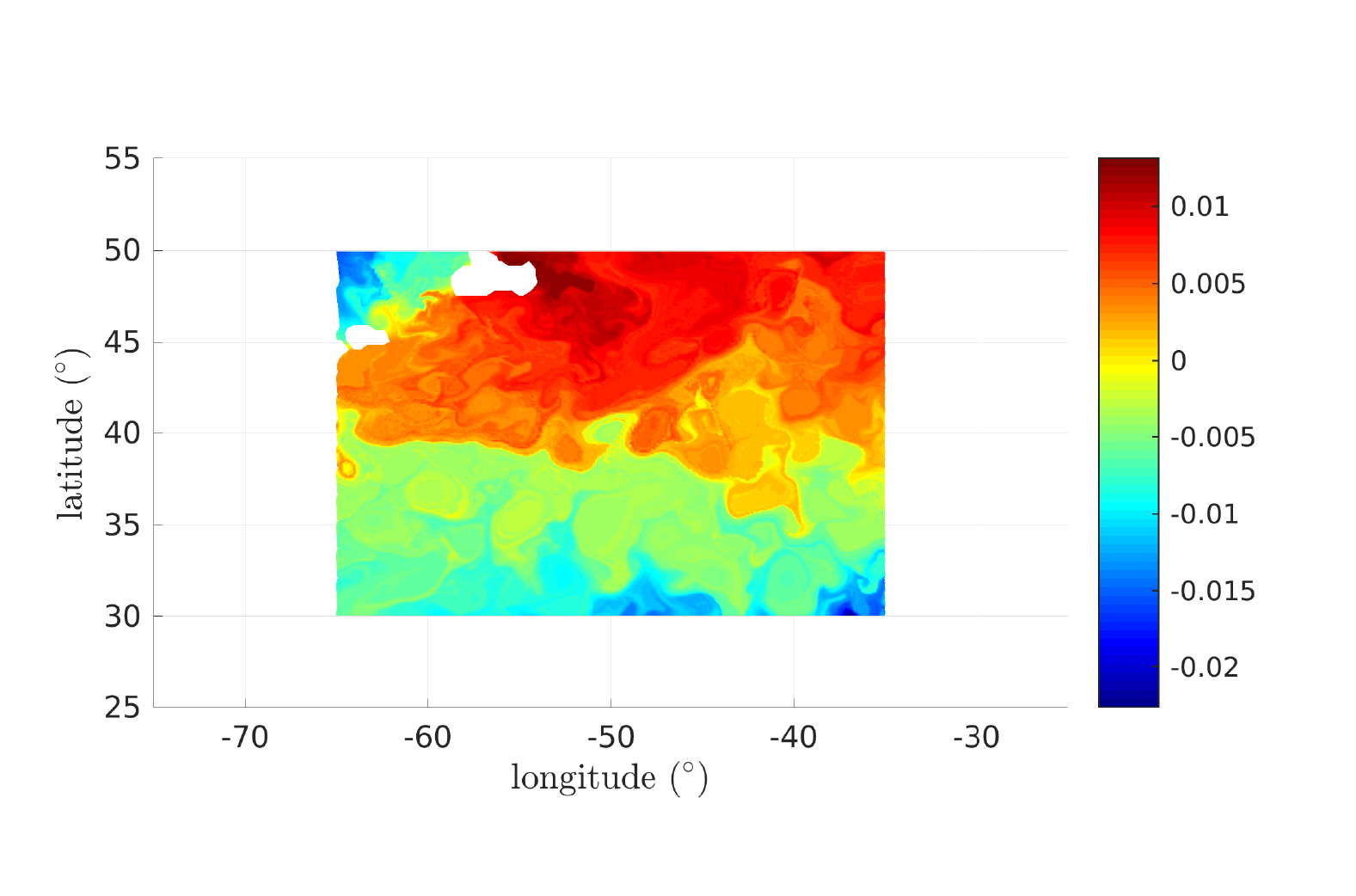}
\includegraphics[width=0.6\textwidth]{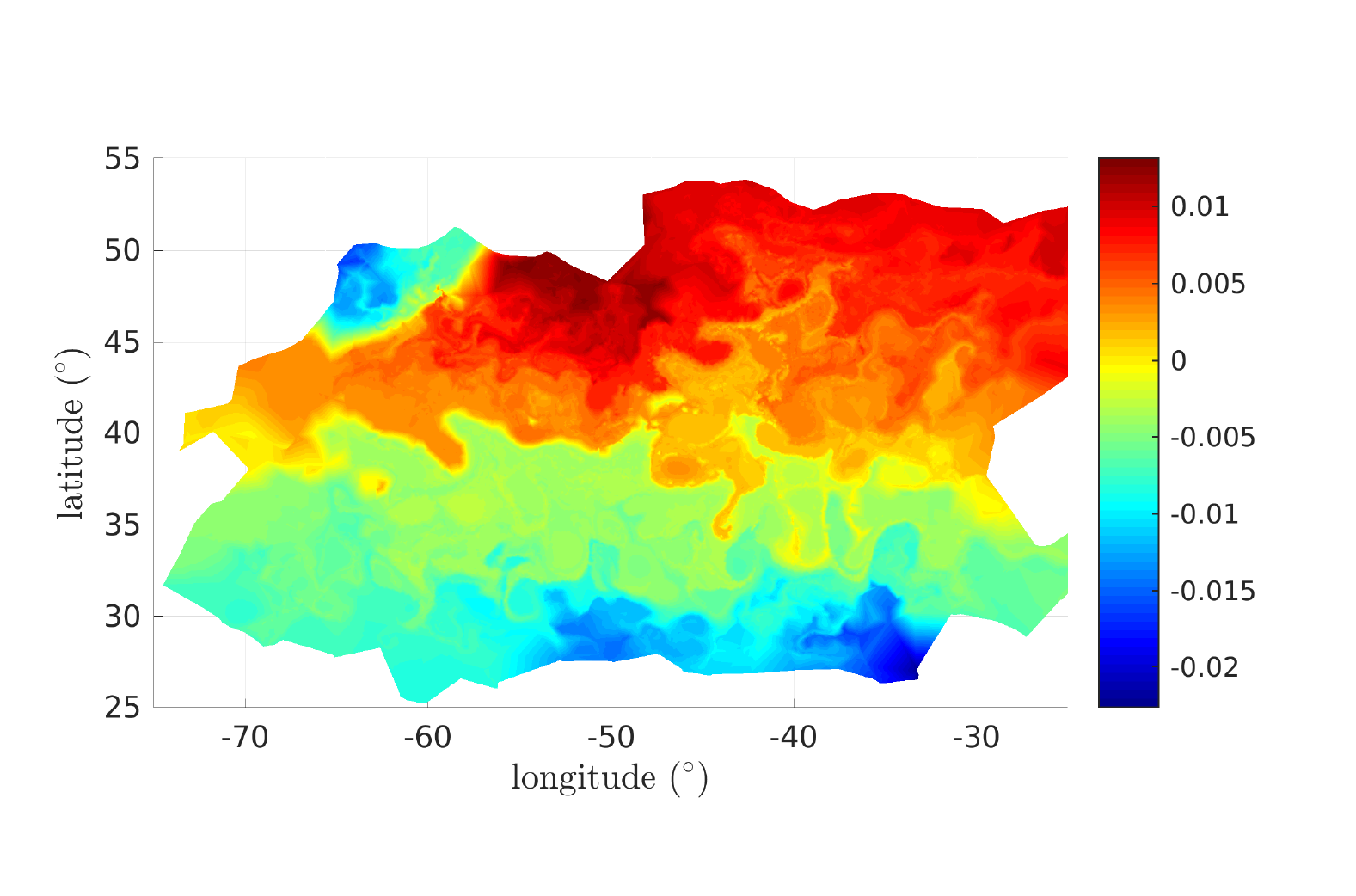}
\caption{Third eigenfunction at 15 January 2015 (upper) and its forward-time image at 15 April 2015 (lower).}
\label{fig:nat_eigenvector3}
\end{figure}

The second eigenfunction highlights the Gulf of St Lawrence in the extreme northwest of the 15 January 2015 domain as a coherent set evolving distinctly from the rest of the domain (not shown), and the third eigenfunction separates the domain at the Gulf Stream;  see the sharp red/green interface in Figure \ref{fig:nat_eigenvector3}.
One may directly compare Figure \ref{fig:nat_eigenvector3} (lower) with \cite[Figure 9]{balasuriya_etal}.
The latter figure shows a satellite-based map of sea surface temperature at 5 April 2015.
The Gulf Stream is seen as a very strong temperature gradient emanating from the east coast of the USA.
There is a very close match with the strong red/green interface (and also the red/orange interface between -50 and -40 longitude) in the eigenvector in Figure \ref{fig:nat_eigenvector3} (lower).
This demonstrates that the strong temperature gradient of the surface ocean is highly correlated with coherent transport of water parcels on the surface.

\begin{figure}[hbt]
\centering
\includegraphics[scale=0.5]{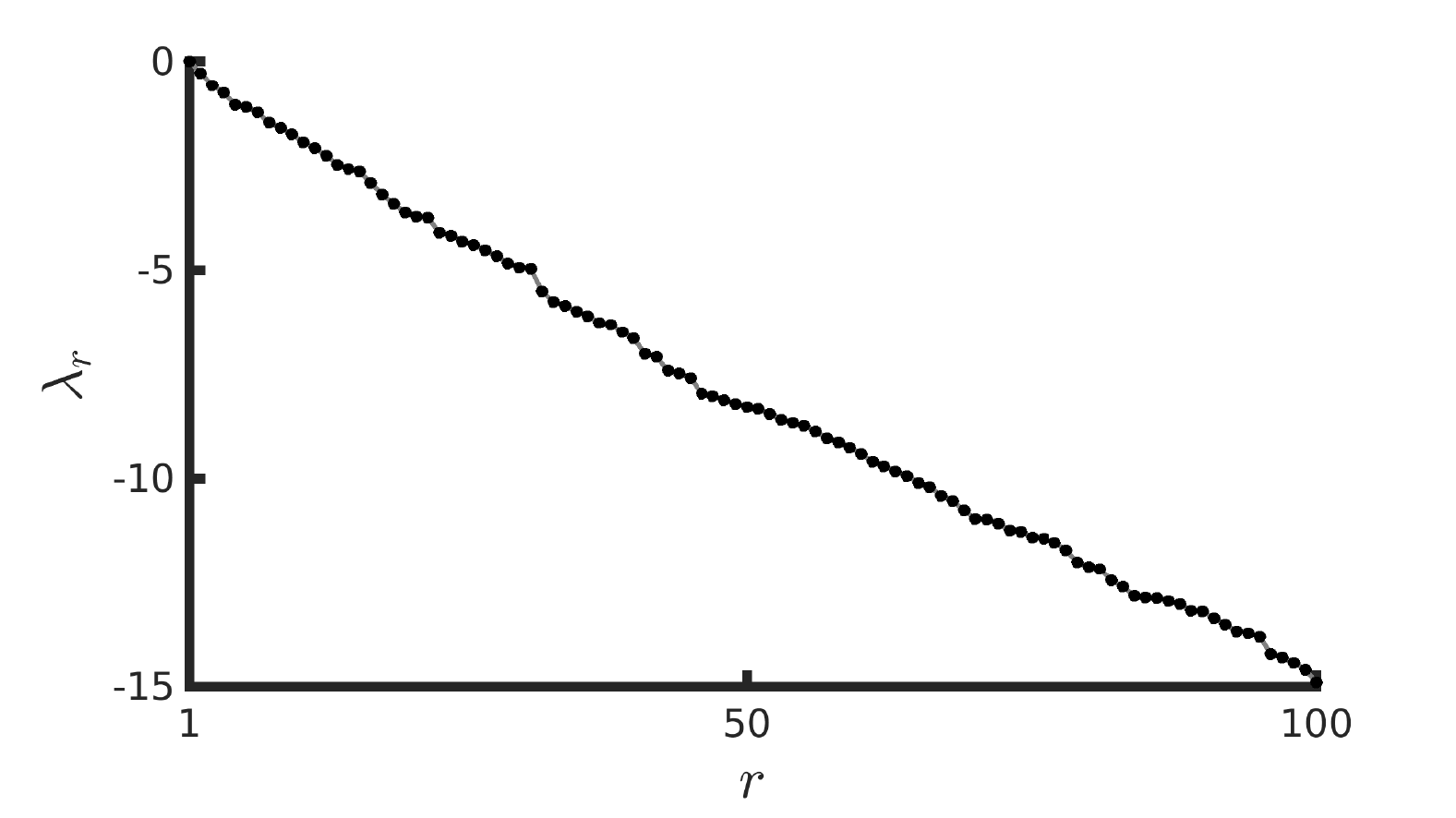}
\includegraphics[scale=0.5]{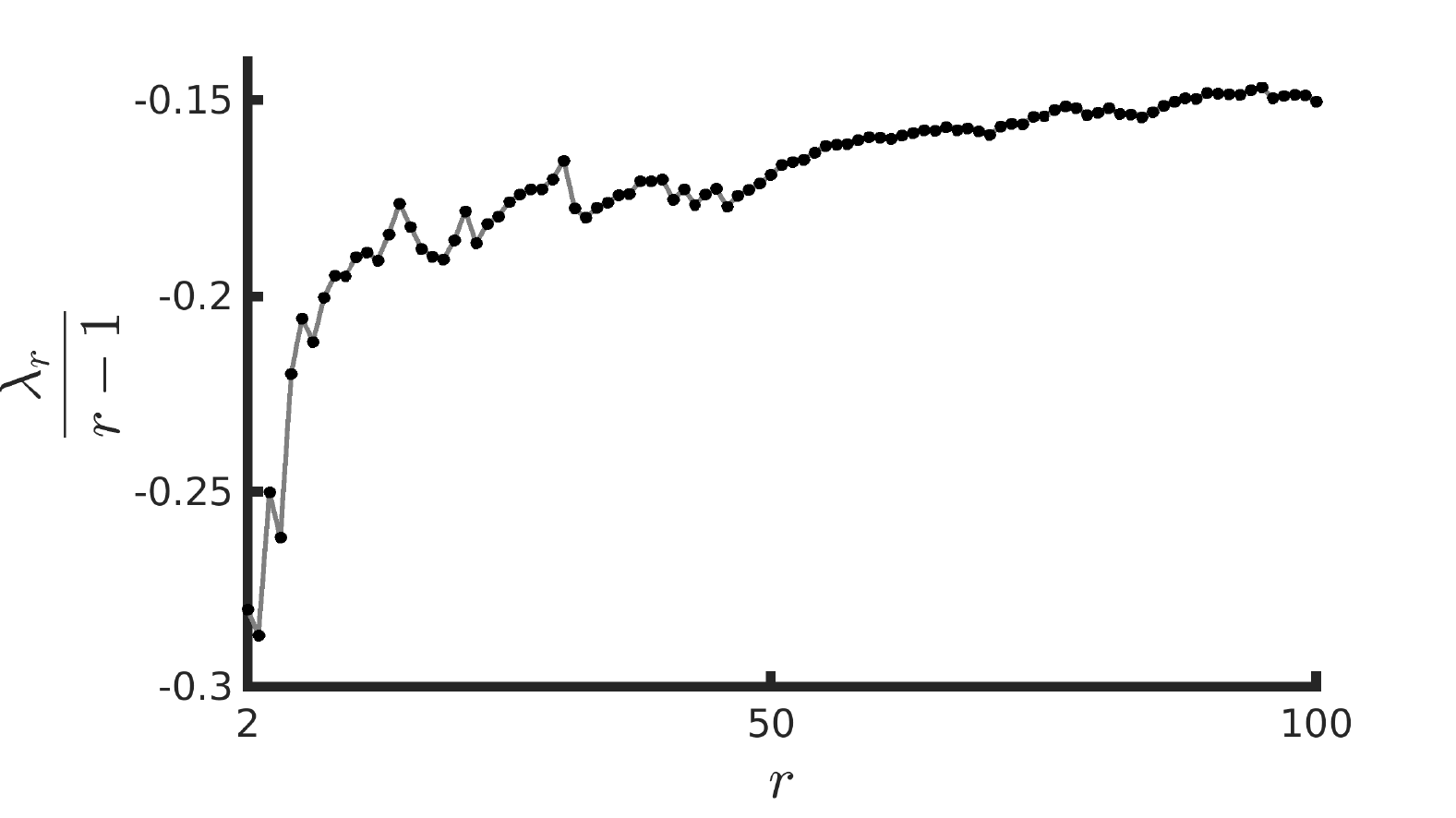}
\caption{Spectrum of the dynamic Laplacian for the North Atlantic dataset. Plot of $\lambda_r$ vs $r$ (left) and $\lambda_{r}/(r-1)$ vs $r$ (right).}
\label{fig:nat_spectrum}
\end{figure}

\begin{figure}[hbt]
\centering
\includegraphics[width=\textwidth]{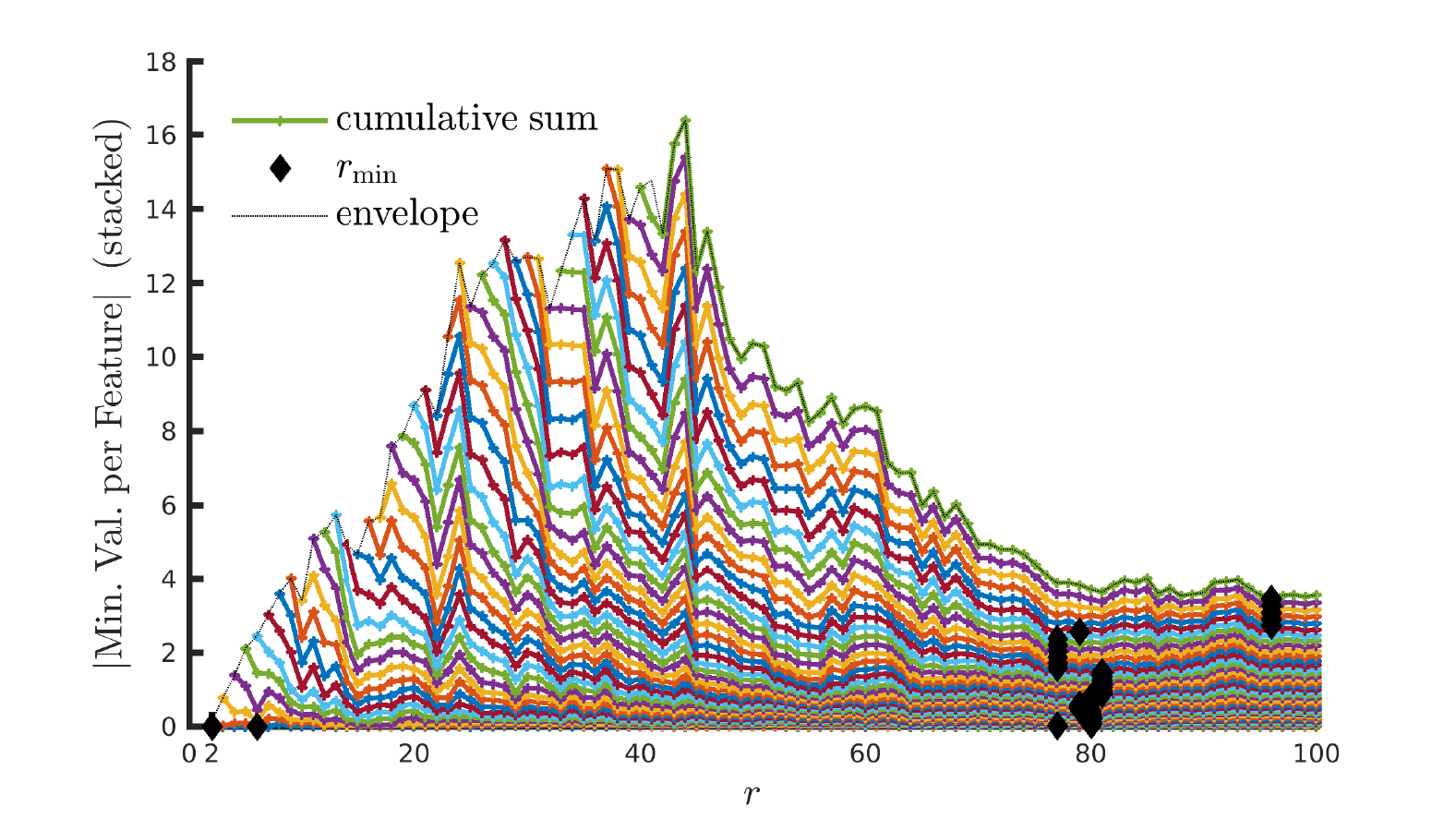}
\caption{North Atlantic negative proportion stacked bar chart} \label{fig:north_atlantic_cumulative}
\end{figure}

\begin{figure}[hbt]
\centering
\includegraphics[width=0.6\textwidth]{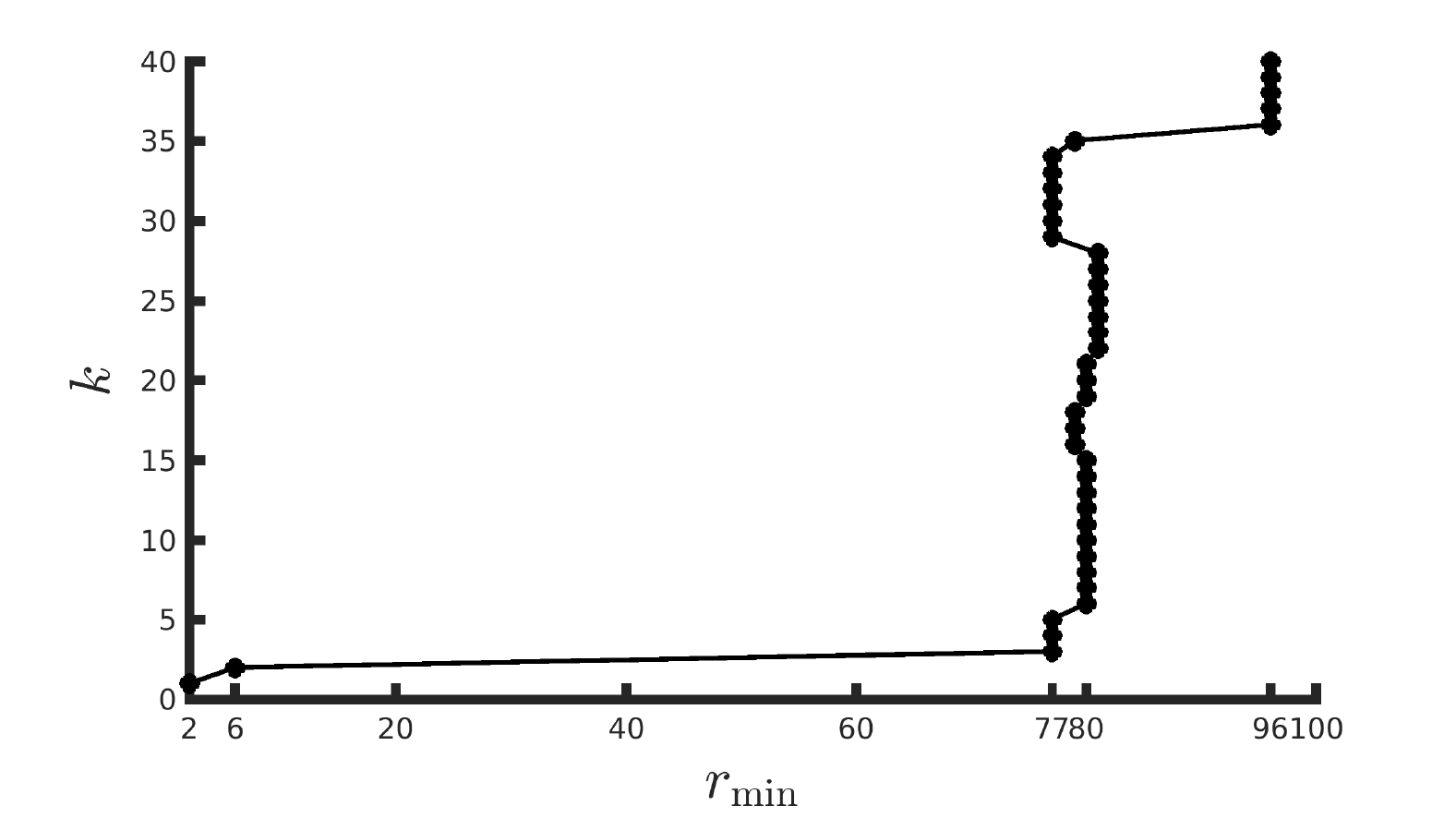}
\caption{$r_{\mathrm{min}}$ vs number of coherent features.} \label{fig:north_atlantic_rmin}
\end{figure}

The spectrum for this dataset, shown in Figure \ref{fig:nat_spectrum} (left) reveals little separation of scales;  the largest drops in the Weyl rescaling (Figure \ref{fig:nat_spectrum}; right) are at $r=3$ and $r=30$.
Therefore, we again select values for $r$ and $k$ using the heuristic described in Section \ref{sec:heuristic}.
We construct the cumulative minimum value plot (Figure \ref{fig:north_atlantic_cumulative}) and plot $r_{\min}(k)$ in Figure \ref{fig:north_atlantic_rmin}.
Our initial choice of $r_{\max}=50$ as our maximum number of eigenvectors was possibly too small to identify many of the coherent sets, because $r_{\min}(k)=50$ for almost all values of $2\le k\le 50$.
Increasing $r_{\max}$ to 100, the coloured curves for $k=1,\ldots,40$ all start to plateau before $r_{\max}=100$, and so we consider this $r_{\max}$ sufficiently large.
We only plot the coloured curves for $k=1,\ldots,40$, to simplify the figure.

\begin{figure}[hbt]
\centering
\includegraphics[width=0.9\textwidth]{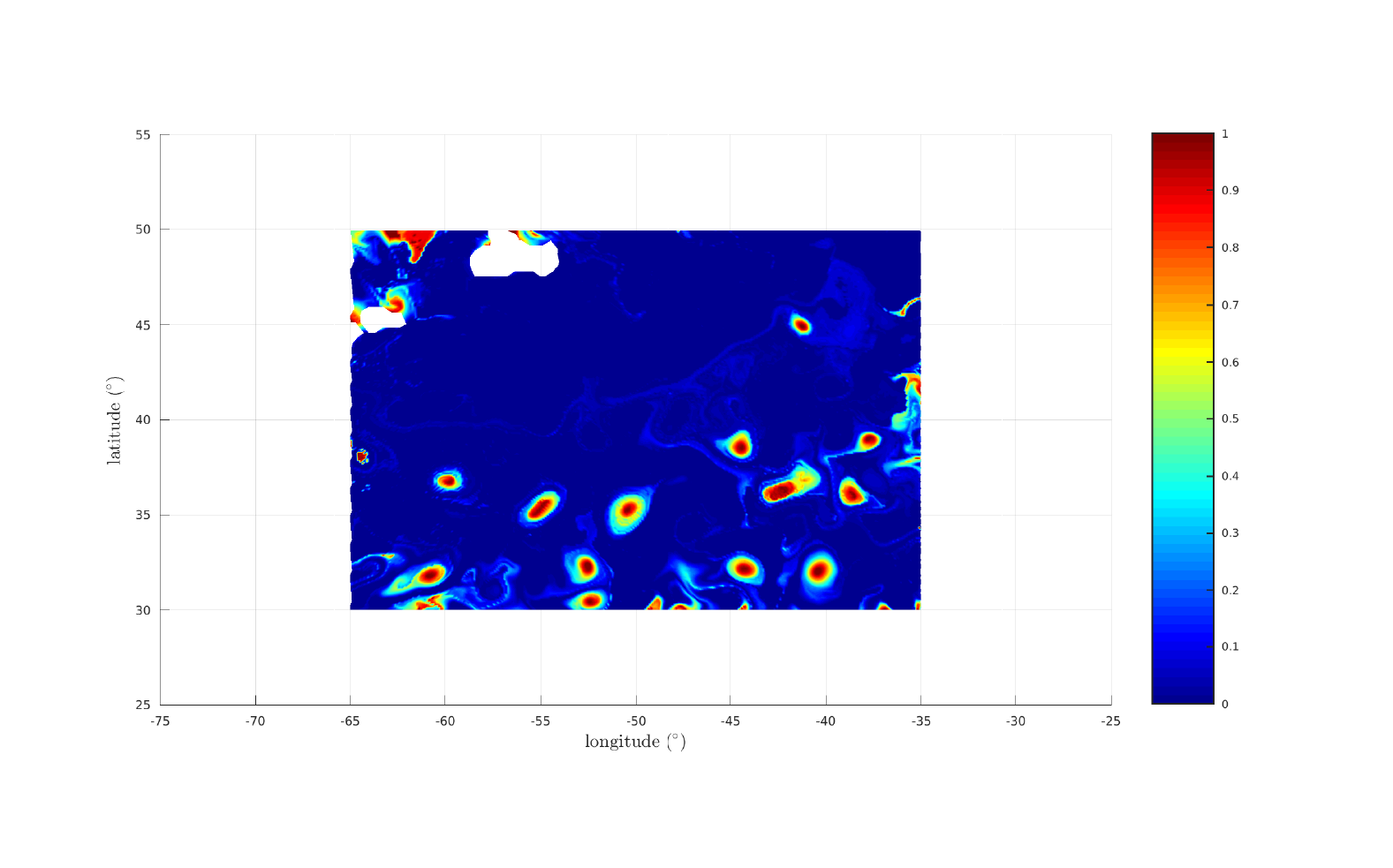}
\includegraphics[width=1\textwidth]{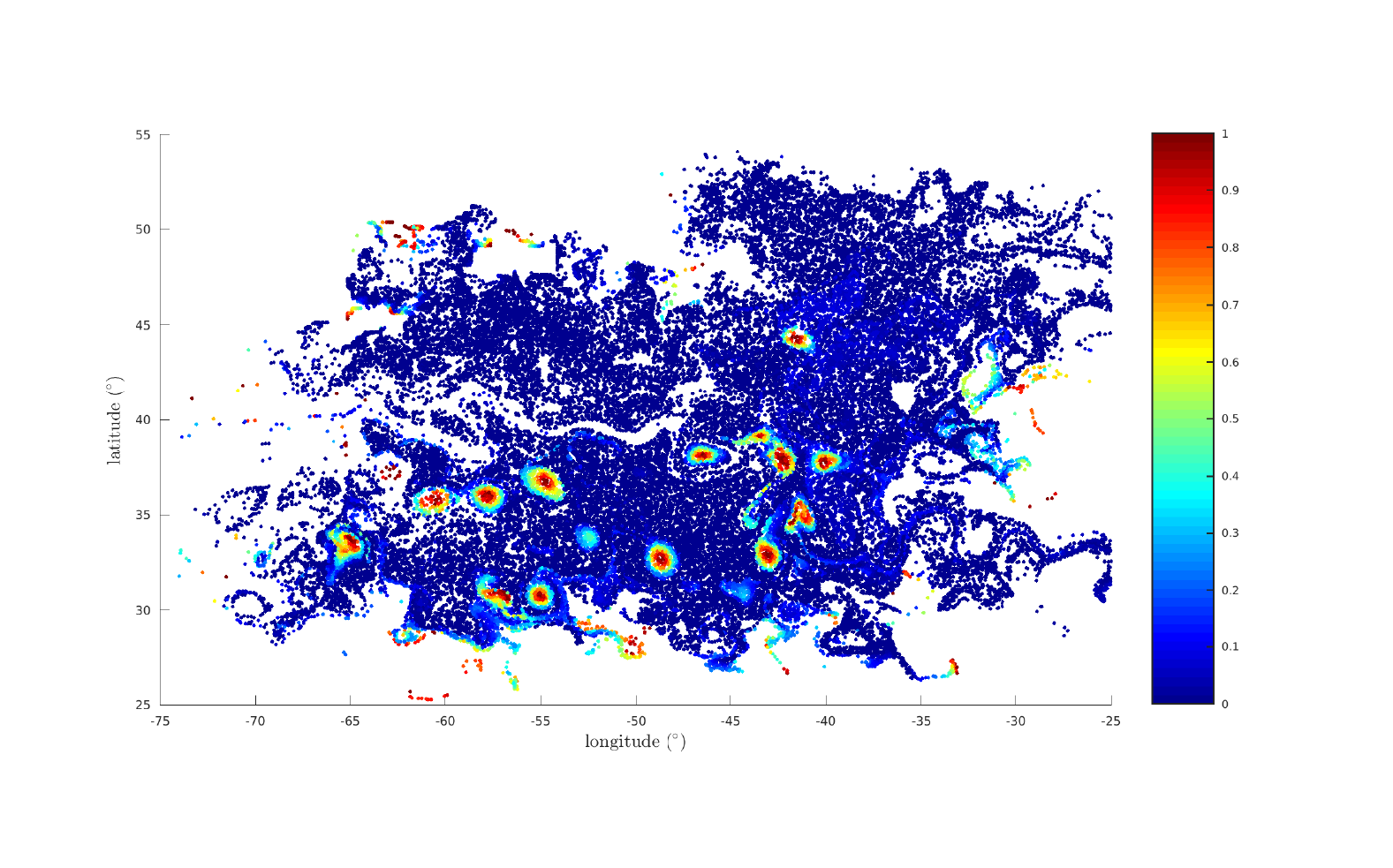}
\caption{Superposition vector $\mathfrak{s}$ of the 36 most reliable features for the North Atlantic dataset using Algorithm \ref{alg1} with $r=96$ at 15 January 2015 (upper) and its forward-time image at 15 April 2015 plotting only final trajectory locations (lower).}
\label{fig:nat_sparse_r96_k36}
\end{figure}

\begin{figure}[hbt]
\centering
\includegraphics[width=0.3\textwidth]{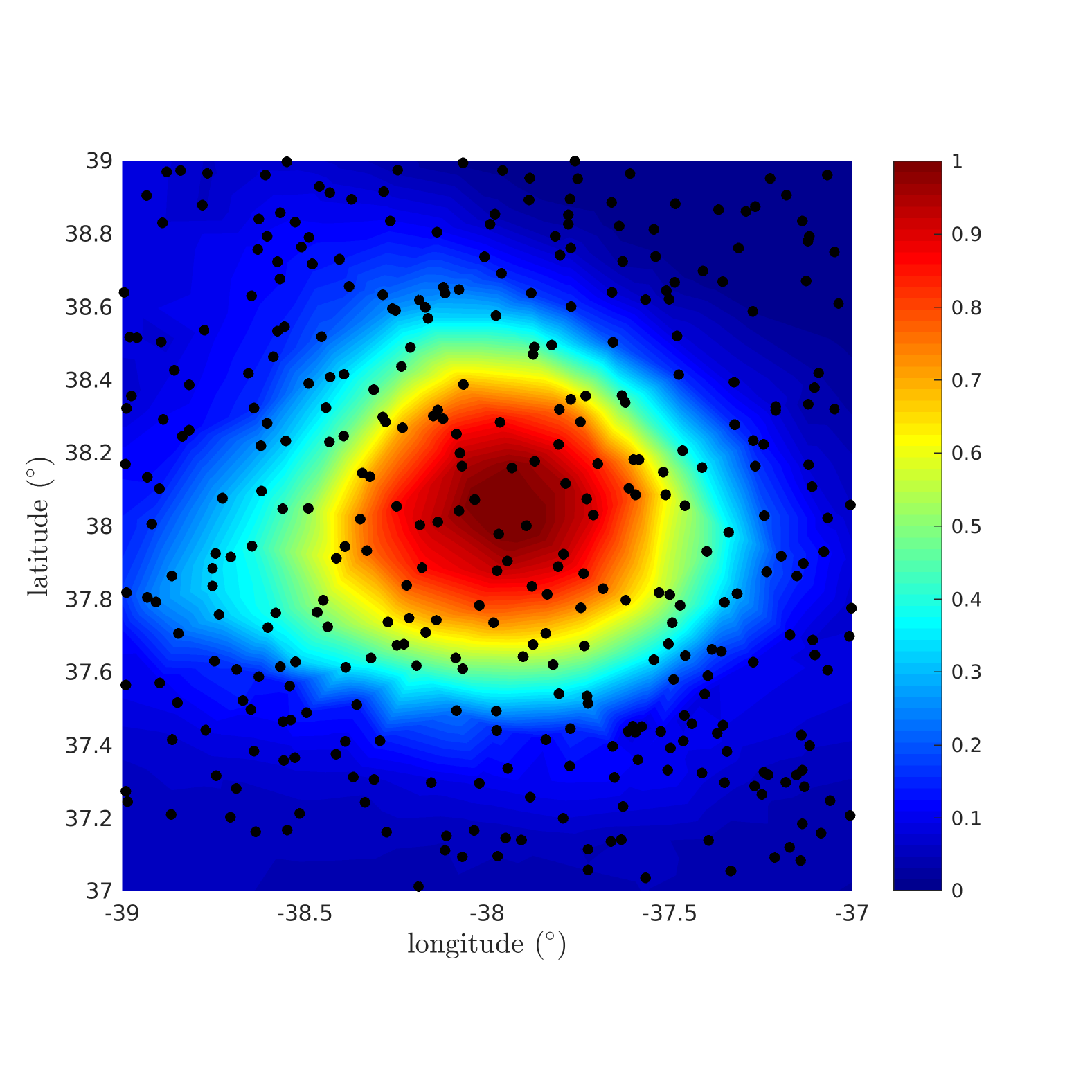}
\includegraphics[width=0.3\textwidth]{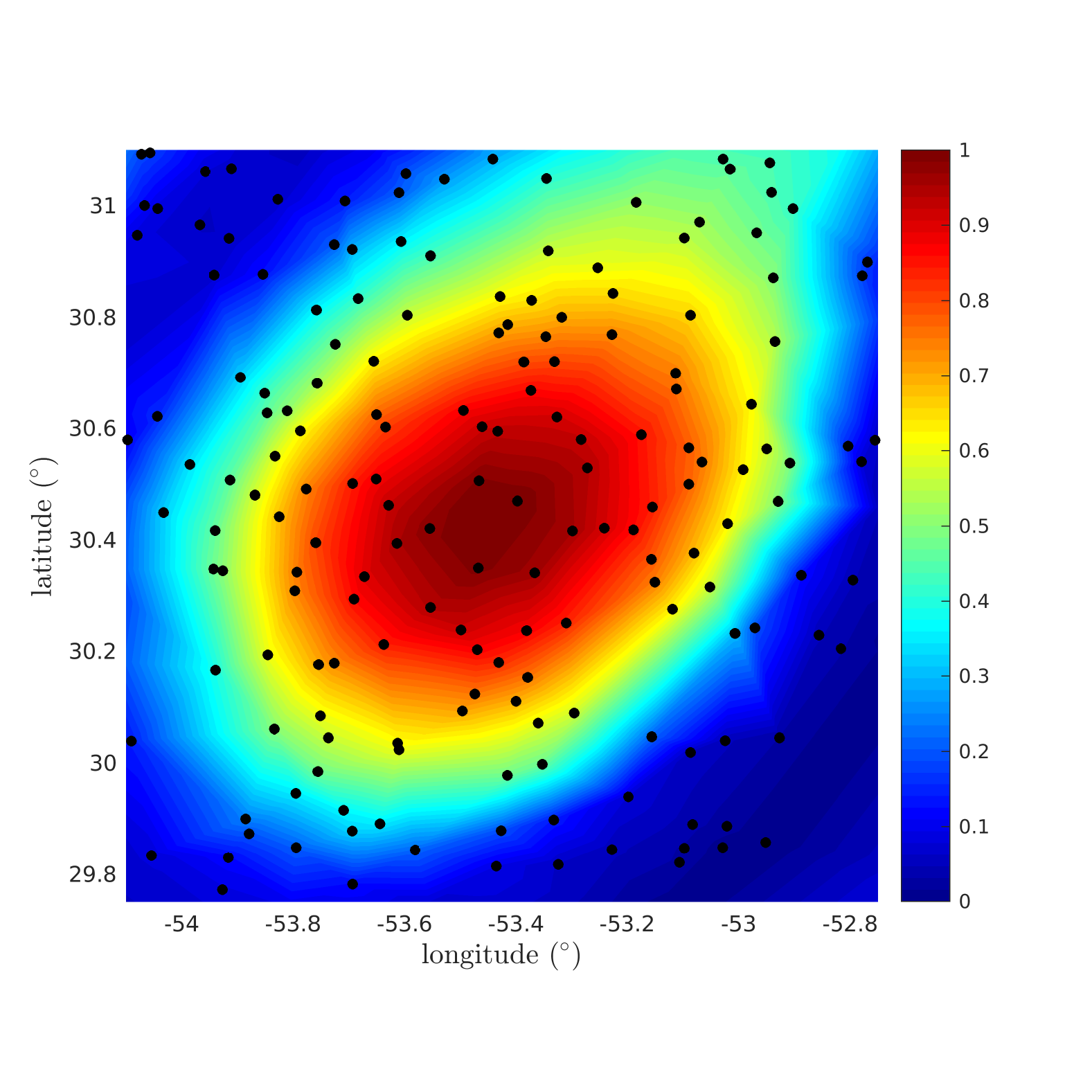}
\caption{Zoom of two sparse vectors produced by Algorithm \ref{alg1} with $(k,r)=(36,96)$, plotted at 28 February 2015.  The black dots are locations of trajectory data used to compute the dynamic Laplace operator.}
\label{fig:eddyzoom}
\end{figure}

As discussed in Subsection \ref{sec:heuristic}, we want to consider the first $k$ value in a run of $k$ values with the same $r_{\min}(k)$. We take $r=96,k=36$. We repeated the experiment with various $k$ and $r_{\min}(k)$ pairs with $r_{\min}(k)$ around 77--81, but the resulting regions were similar to a subset of those for $r=96,k=36$. The span of these 96 sparse vectors produced an approximation of the original subspace with an error of 3.9\%. The absolute sparsity is 3\% and the relative sparsity is 19.1\%.
We plot the sum of the resulting sparse vectors in Figure \ref{fig:nat_sparse_r96_k36}, together with their forward time image. We obtain many coherent features, as well as some partial coherent features at the boundary.
A zoom of two of the sparse vectors highlighting two different eddies at 28 February 2015 is shown in Figure \ref{fig:eddyzoom}.
The black dots in Figure \ref{fig:eddyzoom} indicate trajectory locations used in the construction of the dynamic Laplacian.
Note that we resolve high-resolution estimates of the eddy shapes using the finite-element approach of \cite{FJ18} from lower-resolution trajectory data.

\begin{figure}[hbt]
\centering
\includegraphics[width=0.49\textwidth]{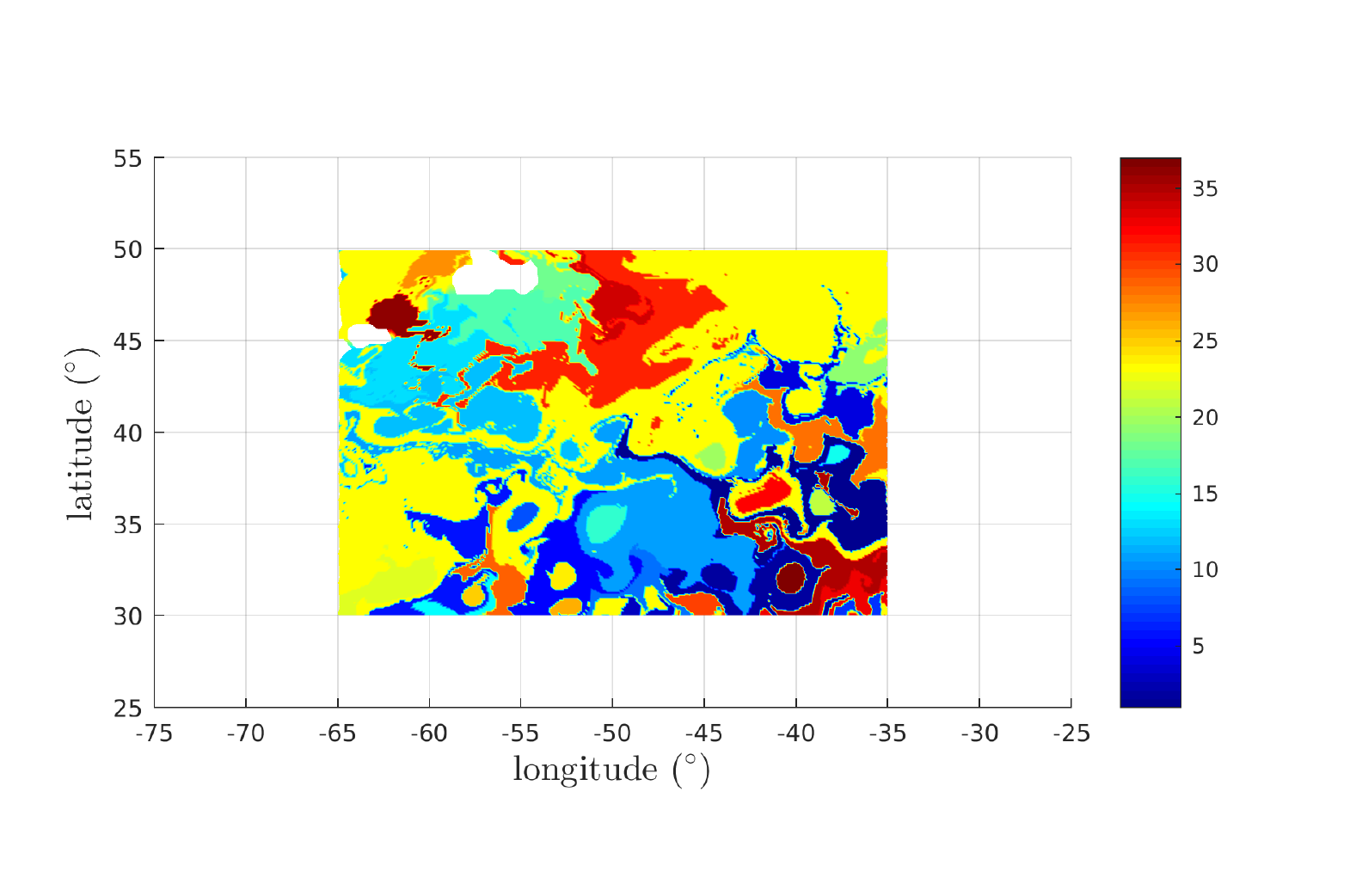}
\includegraphics[width=0.49\textwidth]{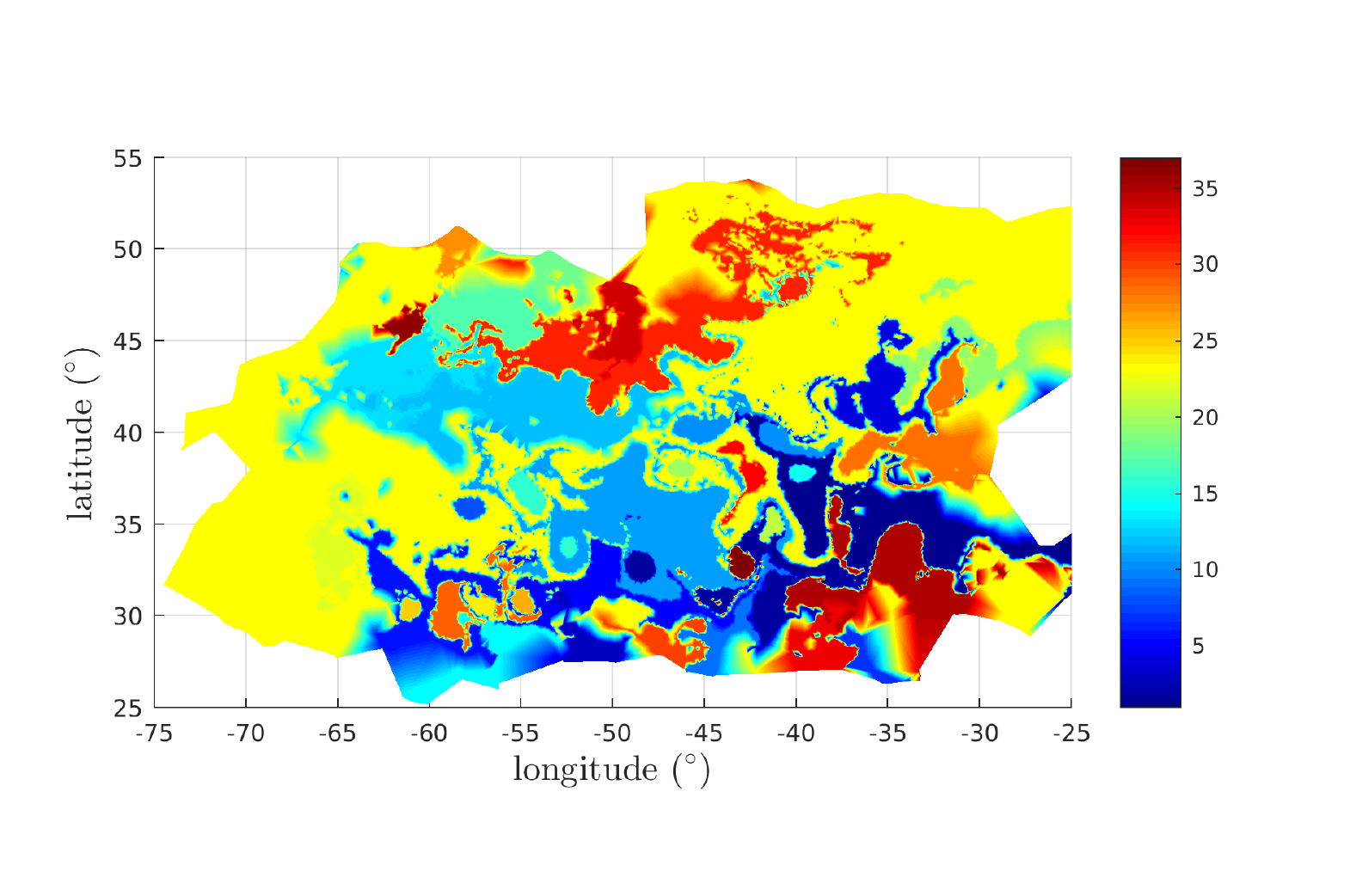}
\caption{Initial (left) and forward-time (right) plot of the k-means partition of state space for the North Atlantic dataset into 37 clusters using $r=96$ eigenvectors.}
\label{fig:nat_kmeans_r96_clust37}
\end{figure}

The eddy field in Figure \ref{fig:nat_sparse_r96_k36} can be compared directly with a backward finite-time Lyapunov exponent field in \cite[Figure 10 (upper right)]{balasuriya_etal}.
The FTLE field suggests a Gulf Stream feature as in Figure \ref{fig:nat_eigenvector3} (lower) and some eddies of similar size as Figure \ref{fig:nat_sparse_r96_k36};  in both cases, the FTLE field is not cleanly identifying these figures.
Figure 4 \cite{rypina_etal} shows fields of encounter volume, exact diffusivity, long-time diffusivity, and diffusive timescale at 11 January 2015 on a domain approximately -75 to -45 degrees longitude by 35 to 42.5 degrees latitude.
These four (spatially smaller) fields can be directly compared with the corresponding part of the larger domain in Figure \ref{fig:nat_sparse_r96_k36} (upper).
The three eddies highlighted in Figure \ref{fig:nat_sparse_r96_k36} (upper) in the subdomain of \cite[Figure 4]{rypina_etal} appear in the latter figure  to varying extents, but not as sharply as in Figure \ref{fig:nat_sparse_r96_k36}.
Other features one might find in \cite[Figure 4]{rypina_etal} do not appear in Figure \ref{fig:nat_sparse_r96_k36} (upper), possibly because they are less coherent than the many other features highlighted.
The partition produced by k-means again proved ineffective;
see Figure \ref{fig:nat_kmeans_r96_clust37}.



\section{Discussion}
\label{sec:conclusion}

We proposed a simple, powerful post-processing of the output of spectral clustering to separate individual features or clusters.
Spectral clustering outputs a truncated eigenspectrum of some Markov- or Laplace-type matrix or operator, with associated eigenbasis.
Features or clusters are encoded in the vectors of this eigenbasis, but may be mixed together within several eigenvectors.
We adapted sparse PCA methodologies to identify an approximate sparse basis for the given eigenbasis.
The use of sparsity  (i) efficiently and reliably separates individual features and (ii) produces a natural cutoff to identify those data points that do not belong to any cluster.

Spectral clustering often relies on the existence of an eigengap to determine (i) the number of clusters or features and (ii) where to truncate the eigenspectrum.
In many situations, there is no clear eigengap.
We proposed a refinement of the standard eigengap procedure in Section \ref{sec:eigengap} that takes into account the spectral structure of Markov and Laplace-type operators, to correctly scale the eigenspectrum to allow fair comparison of gaps along the spectrum.
We further introduced vector-based heuristics in Section \ref{sec:negmass} based on the output of our sparse eigenbasis approximation (SEBA) algorithm to suggest a good choice of the size of eigenbasis $r$.
As described in Section \ref{sec:scales}, these vector-based heuristics determine natural \emph{spatial scales}, while the eigenvalues at eigengaps determine natural \emph{temporal scales}.
A rejection criterion to help determine suitable numbers of features or clusters $k$ was proposed in Section \ref{sec:reliability}.
We described a method for simultaneously selecting $k$ and $r$ in Section \ref{sec:heuristic}.

We believe that our techniques are a useful toolkit for general spectral clustering problems.
One specific motivation was the identification of almost-invariant and coherent sets in time-dependent nonlinear dynamical systems;  in particular, simply and reliably extracting many such sets and determining corresponding natural spatial and temporal scales.
We applied our techniques to three example nonlinear systems.

Firstly, the Bickley jet, which possesses a modest number of strongly coherent sets that were easily separated and extracted using eigenfunctions of the dynamic Laplace operator.
In this example, other post-processing techniques such as k-means clustering also perform very well.
Secondly, turbulent, randomly forced Navier-Stokes flow on a 2-torus, with coherent sets at multiple spatial scales and with no clear temporal scale separation.
We identified good combinations for the ``number of coherent sets'' and the ``eigenbasis size'' and produced well-separated coherent sets through the superposition vector $\mathfrak{s}$.
In contrast, post-processing methods such as k-means performed poorly, as expected for the situation where much of the phase space should not be classified as coherent.
Thirdly, we analysed trajectory data in the North Atlantic Ocean derived from satellite altimetry.
We found a large separation of natural spatial scales with a jump from scales commensurate with the entire ocean basin, to features around one degree of longitude across.
Leading eigenvectors identified the Gulf of St Lawrence as a coherent set strongly dynamically disconnected from the rest of the domain, and the Gulf Stream as a persistent large-scale coherent transport barrier.
No intermediate scales were identified;  the next natural spatial scale occurred at the scale of mesoscale eddies, at which we identified numerous eddies.

\section{Acknowledgements}
The authors would like to thank Anthony J.\ Roberts for suggesting sparse PCA as a possibly useful methodology at a workshop in Blackheath in 2017, Mohammad Farazmand and Alireza Hadjighasem for providing the transition matrix used in Section \ref{sec:turbulence}, and Irina Rypina for providing the AVISO trajectory data used in Section \ref{sec:northatlantic}.
GF is partially supported by an ARC Discovery Project, CPR  is supported by an Australian Government Research Training Program Scholarship, and KS is supported by an ARC Discovery Project.
\appendix






\section{Code} \label{sec:code}
\subsection{MATLAB code for Algorithm \ref{alg1}}
\label{code:alg1}
\begin{lstlisting}
function S=SEBA(V)

%V is pxr matrix (r vectors of length p as columns, assumed orthonormal)
%S is pxr matrix with columns approximately spanning the column space of V

[V,~]=qr(V,0); %Enforce orthonormality
[p,r]=size(V);
mu=0.99/sqrt(p);
S=zeros(size(V));
Rnew=speye(r);  %Initialise rotation
R=0;

while norm(Rnew-R)>1e-14 %tolerance may be decreased
    R=Rnew;
    Z=V*R';
    %Thresholding to solve sparse approximation problem
    for i=1:r
        S(:,i)=sign(Z(:,i)).*max(abs(Z(:,i))-mu,0);
        S(:,i)=S(:,i)/norm(S(:,i));
    end
    %Polar decomposition to solve Procrustes problem
    [P,~,Q]=svd(S'*V,0);
    Rnew=P*Q';
end

%Choose correct parity of vectors and scale so largest value is 1.
for i=1:r
    S(:,i)=S(:,i)*sign(sum(S(:,i)));
    S(:,i)=S(:,i)/max(S(:,i));
end

%Sort so that most reliable vectors appear first (starting with column 1 in S).
[~, I] = sort(min(S), 'descend');
S = S(:, I);
\end{lstlisting}

\subsection{MATLAB code for Algorithm \ref{thresh1}}
\label{code:thresh1}

\begin{lstlisting}
function [S,A,taupu] = subpartition_unity(S)
%SUBPARTITION_UNITY Apply hard thresholding to obtain a sub-partition of 
% unity from the columns of S. Also create a maxmimum-likelihood 
% sub-partition.

S=max(S,0);                         % Take non-negative part
S_descend=sort(S,2,'descend');      % Sort each row in descending order
S_sum=cumsum(S_descend,2);          % Must be <=1 for partition of unity
taupu=max([S_descend(S_sum>1);0]);  % Largest element where row sum > 1 
S(S<=taupu)=0;                      % Apply hard thresholding

[M,A]=max(S,[],2);  % Find column with largest element, by row
A(M==0)=0;          % Suppress zero rows
\end{lstlisting}

\subsection{MATLAB code for Algorithm \ref{thresh2}}
\label{code:thresh2}
\begin{lstlisting}[breaklines=false]
function [S,A,taudp] = disjoint_support(S)
%DISJOINT_SUPPORT Apply hard thresholding to obtain disjoint supports for 
% the columns of S. Also create a maximum-likelihood sub-partition.

S=max(S,0);                     % Take non-negative part of S
S_descend=sort(S,2,'descend');  % Sort to find largest element in each row
taudp=max(S_descend(:, 2));     % Take largest non-leading element
S(S<=taudp)=0;                  % Apply hard thresholding
					
[M,A]=max(S,[],2);  % Find column with largest element, by row
A(M==0)=0;          % Suppress zero rows
\end{lstlisting}


\subsection{MATLAB code for producing a ``natural spatial scales'' plot of $(k,r)$ combinations}
\label{app:scales}
\begin{lstlisting}
function SpatialScalesPlot(V)
%V is pxR matrix (r vectors of length p as columns, assumed orthonormal)

R=size(V,2);
column_minima=nan(R-1,R);
for i=1:R-1
    S=SEBA(V(:,1:i+1));   % Column minima are in descending order
    column_minima(i,1:i+1)=-min(S);
end
minimum_sums=cumsum(column_minima,2);

[~,best_r]=min(minimum_sums);
cla reset;
plot(best_r(1:R)+1,(1:R),['o','-'],'Color','k','MarkerSize',4,...
    'MarkerFaceColor','k');% Plot best r for each k=1,..,r
\end{lstlisting}

\subsection{MATLAB code for producing a ``stacked'' sparse vector minimum value plot}
Note that the first seven lines coincide with the code in Section \ref{app:scales}, so these lines need not be recomputed if one is creating both plots.
\begin{lstlisting}
function MinValStackedPlot(V,kmax)
%V is pxR matrix (r vectors of length p as columns, assumed orthonormal)
%kmax is number of lines to output, assumed <= R

R=size(V,2);
column_minima=nan(R-1,R);
for i=1:R-1
    S=SEBA(V(:,1:i+1));   % Column minima are in descending order
    column_minima(i,1:i+1)=-min(S);
end
minimum_sums=cumsum(column_minima,2);

[best_sum,best_r]=min(minimum_sums);
cla reset; hold on
for j=1:kmax
    iStart=max(j-1,1);
    plot(iStart+1:R,minimum_sums(iStart:R-1,j),'o-',...
        'MarkerSize',2,'Linewidth',1.5);
end
plot(best_r(1:kmax)+1,best_sum(1:kmax),'kd','MarkerSize',4,...
    'MarkerFaceColor','k');% Plot best r for each k
plot([2,2:kmax],minimum_sums([1,R:R:R*(kmax-1)]),'k:',...
    'MarkerSize',2.5) % Plot upper envelope
hold off
\end{lstlisting}

\section{Proof that the expression \eqref{eq:weighted_spcart} can be minimised as outlined in Section \ref{sec:weights}}
\label{app:proof}
%

We wish to solve problem \eqref{eq:diagonal_spcart}, which is equivalent to problem \eqref{eq:weighted_spcart}. We can apply the same alternating optimisation scheme as in Section \ref{sec:sparsebasis}. With fixed $S'$, problem \eqref{eq:diagonal_spcart} is the same as problem \eqref{eq:spcart} with $V$ replaced with $V'$ and $S$ replaced with $S'$. This can be solved using Step \ref{itm:setS} of Algorithm \ref{alg1}.

With fixed $R$, we use the property of the (weighted) Frobenius norm that $\|A\|_{F,\nu}=\|AR\|_{F,\nu}$ for every matrix $A \in \mathbb{R}^{p \times r}$ and orthogonal $R \in \mathbb{R}^{r \times r}$. To see this, note that $\|A\|_{F,\nu}^2=\Trace(A^\top D_\nu A)$, so
\begin{align*}
\|A\|^2_{F,\nu}&=\Trace(A^\top D_\nu A(RR^\top))=\Trace(R^\top A^\top D_\nu AR)=\|AR\|^2_{F,\nu}.
\end{align*}

This lets us reformulate problem \eqref{eq:diagonal_spcart}, for fixed $R$, as
\begin{align}
\argmin_{S'\in\mathfrak{U}^{p,r}}\frac{1}{2}\|V'R^\top-S'\|_F^2 + \mu\|D_\nu^{\frac{1}{2}}S'\|_{1,1}
&= \argmin_{S'\in\mathfrak{U}^{p,r}}\sum_{j=1}^r \frac{1}{2}\|Z_j-S'_j\|_2^2 + \mu\|D_\nu^{\frac{1}{2}}S'_j\|_1, \label{eq:columns}
\end{align}
where $Z_j$ is the $j$th column of $V'R^\top$ and $S'_j$ is the $j$th column of $S'$. Since $\|Z_j\|_2=\|S_j'\|_2=1$, we have $\|Z_j-S'_j\|_2=\|Z_j\|_2+\|S'_j\|_2-2Z_j^\top S'_j=2-2Z_j^\top S'_j$, so problem \eqref{eq:columns} is equivalent to $\argmax_{S'\in\mathfrak{U}^{p,1}}Z_j^\top S'_j-\mu\|D_\nu^{1/2}S'_j\|_1$. By a slight modification of equations (5)-(7) in \cite{journee10a}, taking $p=1$, $A=Z_j^\top$ and $x=1$ and replacing $\gamma \|z\|_1$ with $\sum_{i=1}^n \mu\nu_i^{1/2}|z_i|$, it can be shown that our problem \eqref{eq:columns} is solved for each $j$ by setting 
$S'_j=C'_{\mu,\nu}((V'R^\top)_j)/\|C'_{\mu,\nu}((V'R^\top)_j)\|_2$, exactly as in Step \ref{itm:setR} of Algorithm \ref{alg1} with $V$ replaced with $V'$ and $C_\mu$ replaced with $C'_{\mu,\nu}$.

Lastly, we prove that $C'_{\mu,\nu}(D_\nu^{1/2}z)$ is nonzero for every $z \in \mathbb{R}^p$ with $\|z\|_{2,\nu}=1$, if and only if $\mu<1/\sqrt{\|\nu\|_1}$.
We start by showing that if $\mu<1/\sqrt{\|\nu\|_1}$ then $C'_{\mu,\nu}(D_\nu^{1/2}z)$ is nonzero for every $z \in \mathbb{R}^p$ with $\|z\|_{2,\nu}=1$.
To see this, choose any such $z$, then
\begin{align*}
1=\|z\|_{2,\nu}=\sqrt{\sum_{i=1}^p \nu_i z_i^2}\le \sqrt{\sum_{i=1}^p \nu_i}\max_{1\le j\le p}z_j=\sqrt{\|\nu\|_1} \max_{1\le j \le p}z_j.
\end{align*}
That is, there is some index $j$ with $z_j\ge 1/\sqrt{\|\nu\|_1}$. Recalling that $\nu_i>0$ for every $i=1,\ldots,p$, $\mu<1/\sqrt{\|\nu\|_1}$ implies $|\nu_j^{1/2}z_j|-\mu\nu_j^{1/2}>0$, i.e. the $j$th element of $C'_{\mu,\nu}(D_\nu^{1/2} z)$ is nonzero as required. To prove the reverse direction, note that the vector $c \in \mathbb{R}^p$ with every element equal to $1/\sqrt{\|\nu\|_1}$ has $\|c\|_{2,\nu}=1$. But for $\mu\ge 1/\sqrt{\|\nu\|_1}$, $|\nu_i^{1/2}c_i|-\mu\nu_i^{1/2}\le 0$ for every $i=1,\ldots,p$, so $C'_{\mu,\nu}(D_\nu^{1/2}c)=0$, and the proof is complete.

\end{document}